\documentclass{article}
\pdfoutput=1
\usepackage{graphicx}
\usepackage{curves}
\usepackage{epic}
\usepackage{eepic}
\usepackage{amsmath}
\usepackage{amsfonts}
\usepackage{amsthm}
\usepackage{amssymb}

\theoremstyle{plain}
\newtheorem{Th}{Theorem}[section]
\newtheorem{Lem}[Th]{Lemma}
\newtheorem{Cor}[Th]{Corollary}
\newtheorem{Prop}[Th]{Proposition}

\theoremstyle{definition}
\newtheorem{Def}[Th]{Definition}
\newtheorem*{Not}{Notation}
\newtheorem{Ex}[Th]{Example}
\newtheorem{Assum}[Th]{Assumption}

\theoremstyle{remark}
\newtheorem*{Rem}{Remark}

\def\thm{\begin{Th}}
\def\endthm{\end{Th}}
\def\lemma{\begin{Lem}}
\def\endlemma{\end{Lem}}
\def\cor{\begin{Cor}}
\def\endcor{\end{Cor}}
\def\prop{\begin{Prop}}
\def\endprop{\end{Prop}}
\def\definition{\begin{Def}}
\def\enddefinition{\end{Def}}
\def\remark{\begin{Rem}}
\def\endremark{\end{Rem}}
\def\example{\begin{Ex}}
\def\endexample{\end{Ex}}
\def\demo{\begin{proof}}
\def\enddemo{\end{proof}}

\def\notation{\begin{Not}}
\def\endnotation{\end{Not}}
\def\assumption{\begin{Assum}}
\def\endassumption{\end{Assum}}

\def\A{\mathcal{A}}
\def\B{\mathcal{B}}
\def\C{\mathcal{C}}
\def\D{\mathcal{D}}
\def\E{\mathcal{E}}
\def\F{\mathcal{F}}

\def\H{\mathcal{H}}
\def\M{\mathcal{M}}

\def\P{\mathcal{P}}
\def\Q{\mathcal{Q}}

\def\X{\mathcal{X}}

\def\I{\mathcal{I}}
\def\N{\mathcal{N}}
\def\W{\mathcal{W}}
\def\hE{\widehat{\mathcal{E}}}

\def\BbR{\mathbb{R}}
\def\BbN{\mathbb{N}}

\def\BbZ{\mathbb{Z}}

\def\a{\alpha}
\def\b{\beta}
\def\c{\gamma}
\def\d{\delta}

\def\vp{\varphi}
\def\s{\sigma}

\def\GG{\Gamma}

\def\SS{\Sigma}

\def\inte#1{{\rm int}(#1)}
\def\diam#1{{\rm diam}(#1)}
\def\word#1#2{{#1}_1\ldots{#1}_{#2}}
\def\sd#1#2{#1\backslash#2}
\def\word#1#2{{#1}_1\ldots{#1}_{#2}}

\def\ol#1{\overline{#1}}

\def\norm#1{||#1||}

\def\W{\mathcal{W}}
\def\hE{\widehat{\E}}

\def\wG{\widetilde{G}}

\def\Con{{\rm Con}}
\def\wA{\widetilde{A}}

\usepackage{euscript}

\begin{document}
\begin{center}
{\bf\Large Conductive homogeneity of locally symmetric polygon-based self-similar sets}\\
by\\
{\large J. Kigami \& Y. Ota\\
Graduate School of Informatics\\
Kyoto University}
\end{center}

\begin{abstract}
We provide a rich family of self-similar sets, called locally symmetric polygon-based self-similar sets, as examples of metric spaces having conductive homogeneity, which was introduced in \cite{Ki22} as a sufficient condition for the construction of counterparts of ``Sobolev spaces'' on compact metric spaces. In particular, our results imply the existence of ``Brownian motions'' on our family of self-similar sets at the same time. Unlike the known examples like the Sierpinski carpet~\cite{BB1}, unconstrained carpet~\cite{CaoQiu} and the Octa-carpet~\cite{AndrewsIV, CaHaHu}, our examples may have no global symmetries, i.e. the group of isometries is trivial, as in Figure~\ref{JgonIntro}.
\end{abstract}

\tableofcontents

\section{Introduction}\label{INT}

\begin{figure}[ht]
\centering
\includegraphics[width = 300pt]{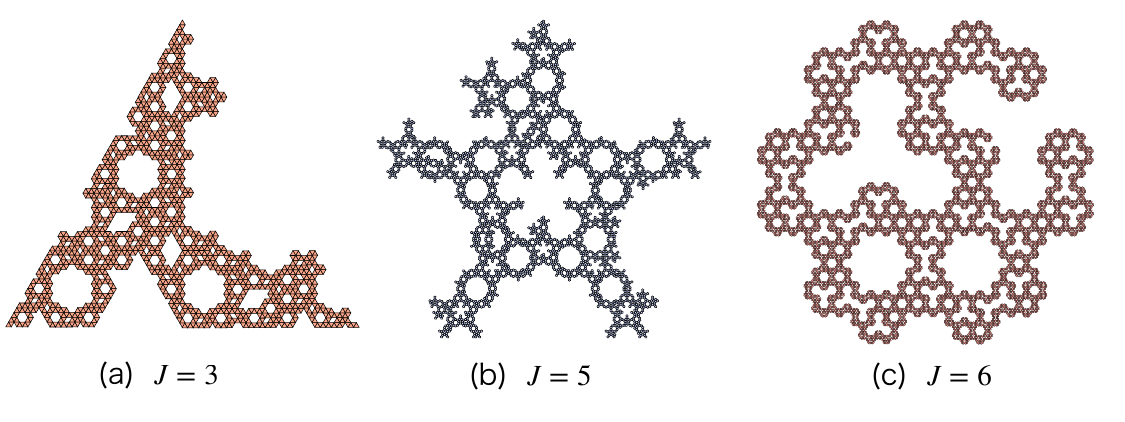}
\vspace{-10pt}
\caption{Locally symmetric $J$-gon-based self-similar sets}\label{JgonIntro}
\end{figure}

The main objective of this paper is to introduce a new class of self-similar sets and to study when self-similar sets in the new class have conductive homogeneity. The new class is called locally symmetric polygon-based self-similar sets, or to be more exact, $G$-symmetric $J$-gon-based self-similar systems where $J \ge 3$ and $G$ is a subgroup of the symmetric group $D_J$ of the regular $J$-gon. We will explain the roles of the local symmetry and the group $G$ later in this introduction. The notion of conductive homogeneity was introduced in \cite{Ki22} as a sufficient condition for the construction of counterparts of ``$(1, p)$-Sobolev spaces'' on a given compact metric space through certain scaling limits of discrete $p$-energies on graphs approximating the original space. In Figure~\ref{JgonIntro}, we present samples of locally symmetric $J$-gon-based self-similar sets having conductive homogeneity.
More details on the self-similar sets in Figure~\ref{JgonIntro}-(a), (b), and (c) will appear in Examples~\ref{GPS.ex10}, \ref{WNG.ex10}, and \ref{WNG.ex10}, respectively.\par

In fact, a variety of analyses on metric spaces, connected with constructions of ``Sobolev spaces'' in one way or another, has been developed since the late 1980s, partly because the notion of fractal had emerged as models of shapes in nature. In general, fractals do not carry ``differentiable'' structures, so that conventional analysis based on differentiation can not be applied. \par
The mainstream of such analyses is based on the idea of upper gradients, which is a generalisation of local Lipschitz constants of Lipschitz continuous functions, as a substitute for derivatives. This idea was initially explored by Haj{\l}asz\cite{Haj1}, Cheeger\cite{Cheeger} and Shanmugalingam\cite{Shanm1} in the 1990s. See \cite{HeiKoShTy}, which gives a panoramic view of this direction of the study of Sobolev spaces on metric spaces since the 1990s. However, the mainstream theory has been known not to work for some self-similar sets, for example, the Vicsek set and the Sierpinski carpet by Kajino-Murugan~\cite{KajMur2, KajMur3}. See \cite[Introduction]{Ki22} for detailed accounts.\par
The notion of conductive homogeneity has been introduced to widen the construction of ``Sobolev spaces'' to the cases as above where the mainstream theory does not apply. The naive idea comes from the following elementary observation about the Sobolev spaces on the unit interval $I = [0, 1]$. Let $f: [0, 1] \to \BbR$ be continuous and let $p > 1$. Define
\[
\E^m_p(f) = \sum_{i = 0}^{2^m - 1}\bigg|f\bigg(\frac{i + 1}{2^m}\bigg) - f\bigg(\frac i{2^m}\bigg)\bigg|^p.
\]
Then $f \in W^{(1, p)}([0, 1])$ if and only if $(2^{p - 1})^m\E_p^m(f)$ converges as $m \to \infty$. Moreover, if $f \in W^{1, p}([0, 1])$, then as $m \to \infty$,
\[
(2^{p - 1})^m\E_p^m(f) \to \int_0^1 |\nabla{f}|^pdx.
\]
To make this fact a prototype, our scenario to define a counterpart of $(1, p)$-Sobolev space on a compact metric space is as follows. Let $(X, d)$ be a compact metric space and let $\{(T_m, E_m)\}_{m \ge 0}$, where $T_m$ is the collection of vertices and $E_m$ is the collection of edges, be a sequence of discrete graphs approximating $(X, d)$. Define a discrete $p$-energy of a function $f: X \to \BbR$ by
\[
\E_p^m(f) = \frac 12\sum_{(x, y) \in E_m} |f(x) - f(y)|^p.
\]
Then the problem is to find a ``nice'' scaling constant $\s_p$ such that the space
\[
\{f|\text{$(\s_p)^m\E^m_p(f)$ ``converges'' somehow as $m \to \infty$}\}
\]
is a nice function space like the $(1, p)$-Sobolev space $W^{1, p}(\Omega)$ for a domain $\Omega$ of an Euclidean space. Roughly speaking, the notion of conductive homogeneity, or to be exact, $p$-conductive homogeneity, is a sufficient condition for this scenario to work well.  In \cite{Ki22}, $p$-conductive homogeneity is shown to enable us to construct a function space $W^p$ having properties analogous to $W^{1, p}(\Omega)$ in the case $p > \dim_{AR}(X, d)$, which is a key geometric constant of the metric space $(X, d)$ called the Ahlfors regular conformal dimension. See Definition~\ref{CHC.def100} for its exact definition. Moreover, see Section~\ref{CHC} for the precise definition and the consequences of conductive homogeneity.\par
\begin{figure}[ht]
\centering
\includegraphics[width = 300pt]{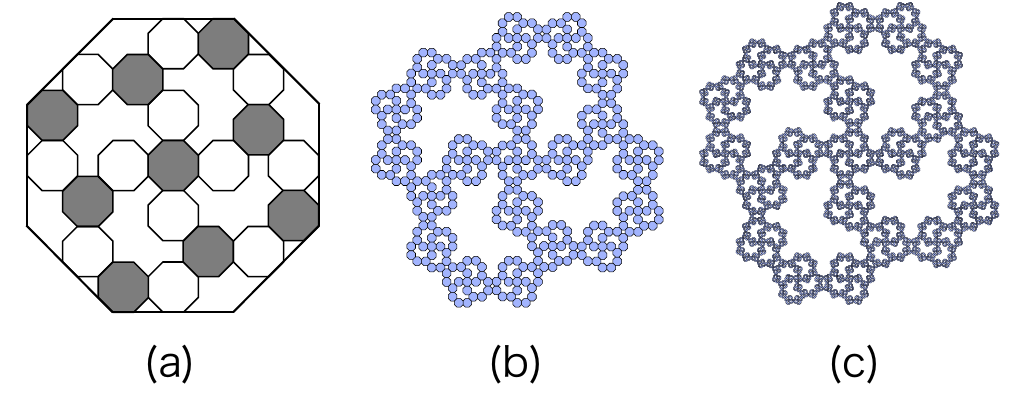}
\caption{$J = 8$, $G = Rot_4$}\label{OctaRot4}
\end{figure}

The main purpose of this paper is to find a rich family of metric spaces having conductive homogeneity, or more precisely, $p$-conductive homogeneity for any $p > \dim_{AR}(X, d)$. As mentioned at the beginning, the family of interest is $G$-symmetric $J$-gon-based self-similar systems. The group $G$ represents the global symmetry of the associated self-similar set. Sufficient conditions for conductive homogeneity will be described in terms of $G$ and its action on the self-similar set. The group $G$ can be trivial, i.e. it can consist only of the identity map. In fact, this is the case for all the examples in Figure~\ref{JgonIntro}. In such a case, conductive homogeneity is determined by the action of $G$ along with the configuration of contracted copies of the regular $J$-gon. See Figure~\ref{OctaRot4}-(a) for an example of such a configuration.\par
\begin{figure}[ht]
\centering
\includegraphics[width = 300pt]{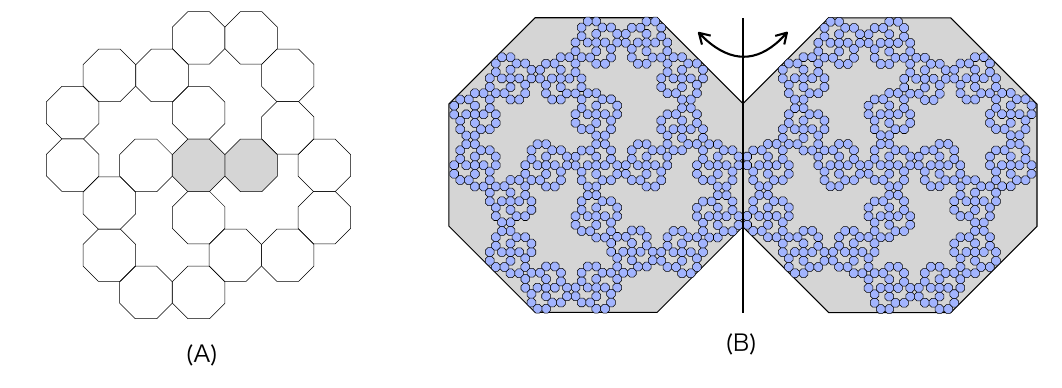}
\caption{Local Symmetry}\label{OctaLS}
\end{figure}

Besides the group $G$, one of the most important features of $G$-symmetric $J$-gon-based self-similar systems is the local symmetry. Roughly speaking, the local symmetry means that if two cells, i.e. contracted $J$-gons in the first step of the iteration of contractions, share a boundary segment, then the copies of self-similar sets in these cells are symmetric under the reflection in the shared boundary line.  For the octagonal example given in Figure~\ref{OctaRot4}, the group $G$ is $Rot_4$, where $Rot_k$ for $k \in \BbN$ is defined as
\[
\text{the group generated by the rotation about the origin $0$ by the angle $2\pi/k$}.
\]
The first step of the generation of the self-similar set is shown as the collection of small octagons in Figure~\ref{OctaRot4}-(a) and Figure~\ref{OctaLS}-(A), which are essentially the same except for the colouring. (Please forget what the colours dark grey and white mean for the moment until Example~\ref{GPS.ex20}, where more detailed accounts of this example will be given.) In Figure~\ref{OctaLS}-(A), we have two light grey cells sharing a boundary segment. The enlargement of this part is Figure~\ref{OctaLS}-(B) and the two small copies of the self-similar set are symmetric under the reflection in the shared boundary segment. See the property (C4) below for a more exact description of the local symmetry.

Since our examples are planar, it is easy to see that $2 > \dim_{AR}(X, d)$, so that they are $2$-conductive homogeneous under any of sufficient conditions for $p$-conductive homogeneity for all $p > \dim_{AR}(X, d)$. Moreover, $2$-conductive homogeneity implies the existence of non-trivial self-similar Dirichlet forms. See Theorem~\ref{CHC.thm30} for details. Consequently, our results in this paper provide new examples of self-similar sets having a non-trivial diffusion process, which could be called Brownian motion. Since the construction of Brownian motions on the generalised Sierpinski carpets by Barlow-Bass\cite{BB1}, there has been no construction of ``Brownian motions'' on infinite-ramified self-similar sets for three decades except unconstrained carpets in \cite{CaoQiu} and the Octa-carpet, see Figure~\ref{GPS.ex13}, and its generalisation in \cite{AndrewsIV} and \cite{CaHaHu}. Note that these examples have the full symmetry of the polygons, namely the square for unconstrained carpets and the regular Octagon for the Octa-carpet. Our Brownian motions constructed through $2$-conductive homogeneity provide new examples of those on regular polygon-based infinitely-ramified self-similar sets which do not necessarily inherit the full symmetry from the regular polygon.\par

For simplicity, we are going to explain exact definitions and results in the case of hexagon-based self-similar sets. Let $Q_*$ be a regular hexagon whose centre is the origin of $\BbR^2$. We assume that $Q_*$ includes its interior and the boundary. Let $\{p_i\}_{i = 0, 1, \ldots, 5}$ be the verticies of $Q_*$ and let $b_i$ be the line segment between $p_i$ and $p_{i + 1}$. In Figure~\ref{Q56}, the hexagon $Q_*$ is illustrated as $Q_*^{(6)}$ in the right-hand side. Moreover, let 
\[
D_6 = \{g| g \in O(2), g(Q_*) = Q_*\},
\]
where $O(2)$ is the orthogonal group of $\BbR^2$. Now, we define our class of self-similar sets.

\definition
Let $S$ be a finite subset and let $G$ be a subgroup of $D_6$. Then a family of contractions $\{f_s\}_{s \in S}$ is called a $G$-symmetric hexagon-based self-similar system if the following conditions are satisfied.\\
(A1)\,\,There exist $r \in (0, 1)$ and  $\{(\vp_s, c_s)\}_{s \in S} \subseteq D_6 \times \BbR^2$ such that
\[
f_s(x) = r\vp_s(x) + c_s
\]
for any $x \in \BbR^2$ and $s \in S$. Moreover, $f_s(Q_*) \subseteq Q_*$ for any $s \in S$. \\
(A2)\,\,For any $i \in \BbZ_J$, there exists $(s, j) \in S \times \BbZ_J$ such that $f_s(b_j) \subseteq b_i$.\\
(A3)\,\,For any $g \in G$, there exists $g_*: S \to S$ such that $(f_{g_*(s)})^{-1}\circ{g}\circ{f_s} \in G$ for any $s \in S$.\\
(A4)\,\,If $f_s(Q_*) \cap f_{s'}(Q_*) \neq \emptyset$ for $s \neq s' \in S$, then 
\[
f_s(Q_*) \cap f_{s'}(Q_*) = f_s(b_i) = f_{s'}(b_j)
\]
for some $i, j \in \{0, 1, \ldots, 5\}$ and
\[
(f_{s'})^{-1}{\circ}R_{s, s'}\circ{f_s} \in G,
\]
where $R_{s, s'}$ is the reflection in the line including $f_s(b_i)$.\\
(A5)\,\,
Define 
\[
E_1^{\ell} = \{(s, s')| s, s' \in S, \text{(A4)-(a) holds}\}.
\]
Then $(S, E_1^{\ell})$ is a connected graph.
\enddefinition

By the classical result on a family of contractions presented in Proposition~\ref{GPS.prop10}, there exists a non-empty compact subset $K$ of $Q_*$ satisfying
\[
K = \bigcup_{s \in S}f_s(K).
\]
This self-similar set $K$ is the main object of our study in this paper. 

\notation
For $w = w_1w_2\ldots{w_m} \in S^m$, we define $f_w = f_{w_1}\circ\ldots{\circ}f_{w_m}$ and $K_w = f_w(K)$.
\endnotation

In Theorem~\ref{GPS.thm10}, each of the conditions (A2), (A3), (A4) and (A5) is shown to imply the following conditions (C2), (C3), (C4) and (C5) respectively:\\
(C2)\,\,Non-degeneracy of $K$, i.e. $K \cap b_j \neq \emptyset$ for any $j = 0, \ldots, 5$,\\
(C3)\,\,$K$ is $G$-symmetric, i.e. For any $g \in G$, $g(K) = K$ for any $g \in G$,\\
(C4)\,\,$K$ is locally symmetric, i.e. if $f_w(Q_*) \cap f_v(Q_*) = f_w(b_i) = f_v(b_j)$ for some $w, v \in S^m$ and $i, j \in \{0, 1, \ldots, 5\}$, then $K_w$ and $K_v$ are symmetric in the line segment $f_w(b_i) = f_v(b_j)$.\\
(C5)\,\, $(K, d_*)$ is connected, where $d_*$ is the restriction of the Euclidean metric on $K$.\\
Our main interest of this paper is when $(K, d_*)$ has $p$-conductive homogeneity for any $p > \dim_{AR}(K, d_*)$. The first key is the action of $G$ on $\BbZ_6 = \{0, 1, \ldots, 5\}$. The action is defined by
\[
g(b_i) = b_{g(i)}
\]
for $i \in \BbZ_J$. This action shows how the group $G$ moves the boundary segments of $Q_*$. The second key is the essential boundary segments $(\BbZ_6)^e$ defined by
\begin{multline*}
(\BbZ_6)^e = \{g(i)| g \in G, \text{there exist $n \ge 1$, $w, v \in S^n$ and $i, j \in \BbZ_6$}\\
\text{such that $f_w(b_i) = f_v(b_j)$.}\}
\end{multline*}
This set represents boundary segments that actually appear as an intersection of two small copies of $Q_*$. Note that it is $G$-invariant, i.e. $G((\BbZ_6)^e) = (\BbZ_6)^e$.\par
To clarify our notations, we give definitions of two subgroups $D_3$ and $D_3^V$ of $D_6$ by
\begin{align*}
D_3 &= Rot_3 \cup \{g_0, g_1, g_2\}\\
D_3^V &= Rot_3 \cup \{\tilde{g}_0, \tilde{g}_1, \tilde{g}_3\},
\end{align*}
where $g_i$ is the reflection in the straight line $p_ip_{i + 3}$ for $i = 0, 1, 2$ and $\tilde{g_i}$ is the reflection in the straight line passing through the midpoints of $b_i$ and $b_{i + 3}$ for $i = 0, 1, 2$. Two groups $D_3$ and $D_3^V$ are isomorphic but their actions on $\BbZ_6$ are different, namely, $D_3^V$ acts transitively on $\BbZ_6$ while $D_3$ does not.\par
\begin{figure}[ht]
\centering
\includegraphics[width = 300pt]{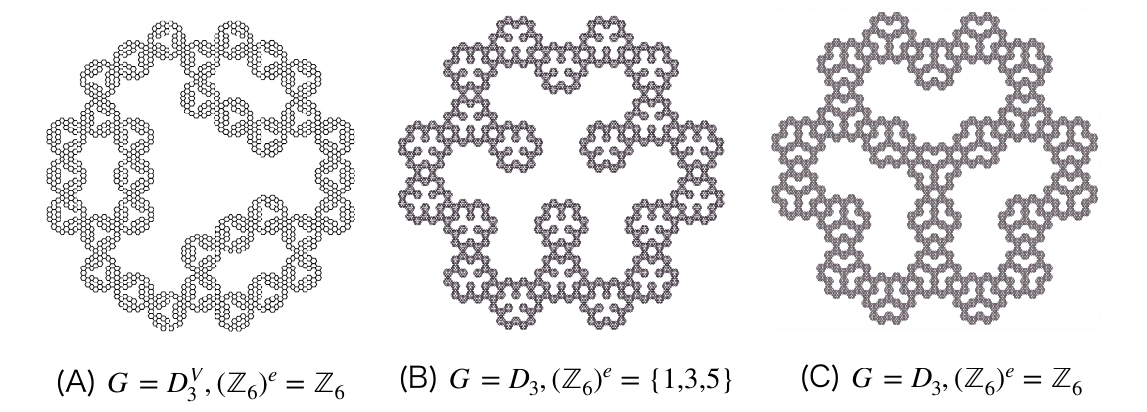}
\vspace{-10pt}
\caption{Hexagonal cases}\label{Hexaintro}
\end{figure}

In Figure~\ref{Hexaintro}, we present three different combinations of $(G, (\BbZ_6)^e)$. More details on these examples (A), (B) and (C) in Figure~\ref{Hexaintro} will be examined in Examples~\ref{RNJ.ex10}, \ref{COP.ex10} and \ref{MIS.ex00}, respectively.\par
The next theorem gives a sufficient condition for conductive homogeneity in terms of $G$ and $(\BbZ_6)^e$.

\thm[= Theorems~\ref{MTH.thm10} and \ref{MTH.thm30}]\label{INT.thm00}
If $G$ acts on $(\BbZ_6)^e$ transitively, i.e. $(\BbZ_6)^e = \{g(i)|g \in G\}$ for some $i \in (\BbZ_6)^e$, then $(K, d_*)$ is $p$-conductively homogeneous for any $p > \dim_{AR}(K, d_*)$.
\endthm

Note that the groups $Rot_6$ and $D_3^V$ act on $\BbZ_6$ transitively. Consequently, $(\BbZ_6)^e = \BbZ_6$ and Theorem~\ref{INT.thm00} applies for those groups. This includes the example in Figure~\ref{Hexaintro}-(A).

In the case when $G = D_3$, the essential boundary segments $(\BbZ_6)^e$ may not be $\BbZ_6$ any longer. In fact, $G = D_3$ and $(\BbZ)^e = \{1, 3, 5\}$ for the example in Figure~\ref{Hexaintro}-(B). However, since $D_3$ acts transitively on $\{1, 3, 5\}$, this example still has $p$-conductive homogeneity for $p > \dim_{AR}(K, d_*)$ by Theorem~\ref{INT.thm00}. \par
On the contrary, $G = D_3$ and $(\BbZ_6)^e = \BbZ_6$ for the example in Figure~\ref{Hexaintro}-(C), and hence the action of $G$ on $(\BbZ_6)^e$ is not transitive. In such a case and the case when $G$ is trivial as in Figure~\ref{JgonIntro}-(C), we need to explore deeper structures of the graph $(T_1, E_1^{\ell})$. Since it is a little too much for the introduction, we leave such cases to later sections where we obtain Theorems~\ref{MIS.thm20} and \ref{WNG.thm20} as our sufficient conditions for conductive homogeneity.\par
For general $J$ of $J$-gon, we have analogous results. Additionally, in the triangle-based cases, the situation is special, i.e. the following theorem holds:

\thm[= Theorem~\ref{MTH.thm10}]\label{INT.thm30}
For a triangle-based case, no matter what $G$ is, $(K, d_*)$ is $p$-conductively homogenous for any $p > \dim_{AR}(K, d_*)$.
\endthm

The organisation of this paper is as follows. In Section~\ref{GPS}, we introduce the notion of $G$-symmetric $J$-gon-based self-similar systems, which is the geometric foundation of this paper, and show its basic properties, in particular, the local symmetry in Theorem~\ref{GPS.thm10}. In Section~\ref{RPB}, we explore further geometric properties of a $G$-symmetric $J$-gon-based self-similar system to verify those required to apply the theory of conductive homogeneity in \cite[Section~2.1]{Ki22}. Section~\ref{ICP} is devoted to studying the topological property called isolated contact points of cells, whose existence makes the theory a little more complicated. Additionally, we introduce the essential boundary segments $(\BbZ_J)^e$ in this section. In Section~\ref{CHC}, we introduce the theory of conductive homogeneity, more precisely, the definition and consequences studied in \cite{Ki22}. In particular, we present an equivalent condition of conductive homogeneity in Theorem~\ref{CHC.thm40}. The equivalent condition plays an essential role in showing conductive homogeneity in later sections. Starting from Section~\ref{CHP}, we present our results on sufficient conditions for conductive homogeneity. To begin with, Section~\ref{CHP} treats sufficient conditions concerning only the group $G$ and its action on $(\BbZ_J)^e$. In the following two sections, Sections~\ref{FBT} and \ref{SBT}, we give two theorems, Theorems~\ref{BAS.thm10} and \ref{BAS.thm20}, which make bridges between Theorem~\ref{CHC.thm40} and the results on conductive homogeneity in this paper. The former theorem describes a relation between conductive homogeneity and the combinatorics of families of paths in a self-similar set, whose detailed structure is given in the latter theorem. Section~\ref{PTM} is devoted to proofs of theorems in Section~\ref{CHP}. In Section~\ref{ENJ}, we consider the case when $J$ is even and $G = D_{J/2}$ and $Rot_{J/2}$. In this case, we need to study detailed structures of the combinatorics of small copies of the original $J$-gon by contraction mappings. The main result is given in Theorem~\ref{MIS.thm20}, which will be proven in Section~\ref{PTE} after the necessary preparations given in Section~\ref{BCM}. In Section~\ref{WNG}, we study the case where $G$ is trivial, i.e. $G$ consists only of the identity. In Proposition~\ref{WNG.prop10}, we show that $G$ is trivial if and only if there exists a ``folding map''. Then we present a sufficient condition for conductive homogeneity in this case in Theorem~\ref{WNG.thm20}, which is proven in Section~\ref{POT}.\par
The results in this paper, except those in Sections~\ref{WNG} and \ref{POT}, are essentially contained in the second author's doctoral dissertation \cite{Ota1}. However, the statements and proofs look different due to the later progress of the study.

\section{Regular polygon-based self-similar sets}\label{GPS}

In this section, we introduce a class of regular polygon-based self-similar sets, which form the geometric framework of this paper. In what follows, $J$ is a natural number with $J \ge 3$ and $Q_*^{(J)}$ represents a unit regular $J$-gon. More precisely, it is defined as follows.

\definition[Regular $J$-gon]\label{GPS.def10}
Let $J \in \BbN$ with $J \ge 3$ and let
\[
p_i = \frac{1}{\cos{(\frac{\pi}J})}\bigg(\cos{\Big(\frac {\pi}{J} + \frac{2\pi}Ji - \frac{\pi}2\Big)}, \sin{\Big(\frac{\pi}{J} + \frac{2\pi}Ji - \frac{\pi}2\Big)}\bigg).
\]
for $i \in \BbZ_{J}$, where $\BbZ_J = \BbZ/J\BbZ = \{0, \ldots, J - 1\}$ equipped with a distance $\d_J(i, j)$ for $i, j \in \BbZ_J$ defined as
\[
\d_J(i, j) = \min\big\{|k - l|\,\big| k, l \in \BbZ, k \equiv i\!\!\!\mod J,\, l \equiv j\!\!\!\mod J\big\}.
\]
If no confusion may occur, we omit the subscript $J$ in $\d_J$ and use $\d$ to denote $\d_J$.
Define $b_i$ as the line segment connecting $p_{i - 1}$ and $p_i$ for $i \in \BbZ_J$. Throughout this paper, we identify the collection of boundary segments, $\{b_i| i \in \BbZ_J\}$, with $\BbZ_J$ in the natural manner. Moreover, define
\[
Q_*^{(J)} = \{tx| t \in [0, 1], x \in b_i\,\,\text{for some $i \in \BbZ_{J}$}\}.
\]
If no confusion may occur, we use $Q_*$ to denote $Q_*^{(J)}$.
\enddefinition

\begin{figure}[ht]
\centering
\includegraphics[width = 300pt]{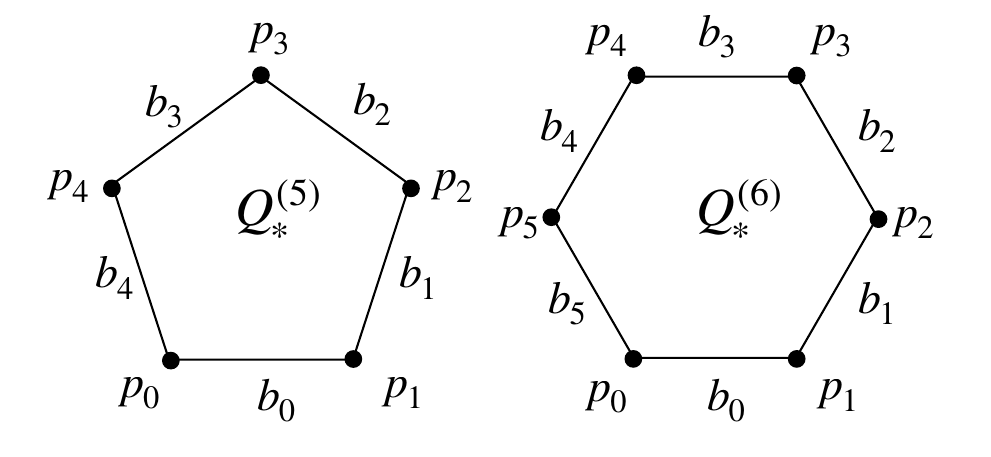}
\caption{$Q_*^{(5)}$ and $Q^{(6)}_*$}\label{Q56}
\end{figure}

Namely $Q_*^{(J)}$ is the regular $J$-gon whose vertices are $p_0, \ldots, p_{J - 1}$. See Figure~\ref{Q56} for the pentagonal and hexagonal cases.\par
In the definition of $p_i$, the normalizing factor $1/\cos{(\pi/J)}$ is applied to make the incircle of $Q_*^{(J)}$ the circle with the center $(0, 0)$ and radius $1$. Moreover, it follows that $Q_*^{(nJ)} \subseteq Q_*^{(J)}$ for any $n \in \BbN$.\par

\definition\label{GPS.def20}
Define $R_{\theta}$ as the reflection in the line $\{(t\cos{\theta}, t\sin{\theta})| t \in \BbR\}$ and define $\Theta_{\theta}$ as the rotation about the origin $0$ by an angle $\theta$.
\enddefinition

\definition\label{GPS.def30}
Define
\[
D_J = \{g| g \in O(2), g(Q_*^{(J)}) = Q_*^{(J)}\}.
\]
and 
\[
Rot_J = \{\Theta_{\frac{2\pi}Ji}| i \in \BbZ_J\}.
\]
\enddefinition

It is easy to see that
\[
D_{J} = Rot_J \cup \begin{cases}
\{R_{\frac{\pi}{J}i - \frac {\pi}2}| i \in \BbZ_{J}\}\quad&\text{if $J$ is even},\\
\{R_{\frac{\pi}{J} + \frac{2\pi}Ji - \frac{\pi}2}| i \in \BbZ_{J}\}\quad&\text{if $J$ is odd}.
\end{cases}
\]

\definition\label{GPS.def40}
Define $\rho: \BbZ_{J} \to \BbR/2\pi$ by 
\[
\rho(i) = \frac{2\pi}Ji
\]
for any $i \in \BbZ_{J}$.
\enddefinition

\lemma\label{GPS.lemma00}
 $R_{\rho(i)}$ is the reflection in the line parallel to $b_i$. Moreover, define $D_J^*$ as the subgroup of $O(2)$ generated by $D_J$ and $\{R_{\rho(i)}\}_{i \in \BbZ_J}$. Then
\[
D_J^* = \begin{cases}
D_J\quad&\text{if $J$ is even,}\\
D_{2J}\quad&\text{if $J$ is odd.}
\end{cases}
\]
\endlemma

The next definition gives the geometric framework of self-similar sets studied in this paper. 

\definition\label{GPS.def50}
Let $S$ be a finite set and let $\{f_s\}_{s \in S}$ be a family of contractions on $\BbR^2$. Let $G$ be a subgroup of $D_J$. The triple $(S, \{f_s\}_{s \in S}, G)$ is called a $G$-symmetric $J$-gon-based self-similar system, ``$(J, G)$-s.s. system'' for short, if and only if the following assumptions (A1) through (A5) are satisfied:\\
(A1)\,\,There exist $r \in (0, 1)$ and  $\{(\vp_s, c_s)\}_{s \in S} \subseteq D_J^* \times \BbR^2$ such that
\[
f_s(x) = r\vp_s(x) + c_s
\]
for any $x \in \BbR^2$ and $s \in S$. Moreover, $f_s(Q_*) \subseteq Q_*$ for any $s \in S$. \\
(A2)\,\,For any $i \in \BbZ_J$, there exists $(s, j) \in S \times \BbZ_J$ such that $f_s(b_j) \subseteq b_i$.\\
(A3)\,\,For any $g \in G$, there exists $g_*: S \to S$ such that $(f_{g_*(s)})^{-1}\circ{g}\circ{f_s} \in G$ for any $s \in S$.\\
(A4)\,\,If $f_s(Q_*) \cap f_{s'}(Q_*) \neq \emptyset$ for $s \neq s' \in S$, then either (a) or (b) below holds.\par\vspace{5pt}\noindent
(a)\,\,$f_s(Q_*) \cap f_{s'}(Q_*) = f_s(b_i) = f_{s'}(b_j)$ for some $i, j \in \BbZ_{J}$ and
\[
(\vp_{s'})^{-1}{\circ}\vp_s{\circ}R_{\rho(i)} \in G
\]
\noindent(b)\,\,$f_s(Q_*) \cap f_{s'}(Q_*) =f_s(p_i) = f_{s'}(p_j)$ for some $i, j \in \BbZ_{J}$.\par\vspace{5pt}\noindent
(A5)\,\,
Define 
\[
E_1^{\ell} = \{(s, s')| s, s' \in S, \text{(A4)-(a) holds}\}.
\]
Then $(S, E_1^{\ell})$ is a connected graph.
\enddefinition

\remark
For simplicity, we sometimes omit the subscript $*$ in $g_*$ and use $g$ to denote $g_*$.
\endremark

Examples of $(J, G)$-s.s.\,systems will be presented later. See Examples~\ref{GPS.ex10} and \ref{GPS.ex20} for example.

\remark
By \eqref{GPS.eq20}, the condition (A4)-(a) is equivalent to the following statement:\\
$f_s(Q_*) \cap f_{s'}(Q_*) = f_s(b_i) = f_{s'}(b_j)$ for some $i, j \in \BbZ_J$ and 
\[
(f_{s'})^{-1}{\circ}R_{s, s'}{\circ}f_s \in G,
\]
where $R_{s, s'}$ is the reflection in the line segment $f_s(Q_*) \cap f_{s'}(Q_*)$.
\endremark
\remark
By (A5), we see that $\cup_{s \in S} f_s(Q_*)$ is a connected subset of $Q_*$.
\endremark

Hereafter in this section, we always assume that $(S, \{f_s\}_{s \in S}, G)$ is a $(J, G)$-s.s. system.

\prop\label{GPS.prop10}
There exists a unique non-empty compact subset $K$ of $Q_*$ such that
\[
K = \bigcup_{s \in S}f_s(K),
\]
which is called the self-similar set associated with $(S, \{f_s\}_{s \in S})$. Moreover, let $d_*$ be the normalized restriction of the Euclidean metric to $K$, i.e.
\[
d_* (x, y) = \frac{|x - y|}{\sup_{a, b \in K} |a - b|}
\]
for $x, y \in K$. Then
\[
\dim_H(K, d_*) = -\frac{\log{\#(S)}}{\log r},
\]
where $\#(A)$ is the number of elements of a set $A$ and $\dim_H(K, d_*)$ be the Hausdorff dimension of $(K, d_*)$.
\endprop

\demo
The existence and uniqueness of $K$ follow from \cite[Theorem~1.1.4]{AOF}. Let $O = \inte{Q_*}$, where $\inte{A}$ is the interior of a set $A \subseteq \BbR^2$ with respect to the Euclidean metric. By (A4), it follows that if $s \neq s'$, then $f_s(O) \cap f_{s'}(O) = \emptyset$, and $f_s(O) \subseteq O$ for any $s \in S$. Thus, the open set condition holds. Hence we have the desired results by \cite[Corollary~1.5.9]{AOF}. 
\enddemo

The next definition is a collection of basic notations on self-similar sets.

\definition\label{GPS.def60}
Define 
\[
T_m = S^m = \{\word sn| s_1, \ldots, s_m \in S\}, 
\]
and
\[
T = \bigcup_{m \ge 0} T_m,
\]
where $T_0 = \{\phi\}$. For $w \in T$, define $|w|$ as the unique number $n$ satisfying $w \in T_n$. For $m \ge 1$ and $\word sm \in T_m$, define
\[
f_{\word sm} = f_{s_1}\circ\cdots\circ{f_{s_m}}, \quad\vp_{\word sm} = \vp_{s_1}\circ\cdots\circ\vp_{s_m}
\]
\[
Q_{\word sm} = f_{\word sm}(Q_*^{(J)}), \quad K_{\word sm} = f_{\word sm}(K), \quad c_{\word sm} = f_{\word sm}(0)
\]
and
\[
b_i(\word sm) = f_{\word sm}(b_i).
\]
In case $w = \phi \in T_0$, $f_{\phi}$ is defined as the identity map. Moreover, for $A \subseteq T_n$ and $m \ge 0$, we define 
\[
S^m(A) = \{w\word sm| w \in A, s_1, \ldots, s_m \in S\}.
\]
For simplicity, we write $S(A) = S^1(A)$ and $S^m(w) = S^m(\{w\})$ for $w \in T$.
\enddefinition

\remark
It follows that
\[
f_w(x) = r^{|w|}\vp_w(x) + c_w
\]
for any $w \in T$ and $x \in \BbR^2$, where $c_w = f_w(0)$. Moreover, since $d_*$ is normalized, we see that $\diam{K_w, d_*} = r^{|w|}$ for any $w \in T$, where $\diam{A, d}$ is a diameter of a subset $A$ of a metric space $(X, d)$.
\endremark

\remark
Let $E = \{(w, v)| w, v \in T, w \in S(v)\,\,\text{or}\,\, v \in S(w)\}$ and set $\A: T \times T \to \{0, 1\}$ by
\[
\A(w, v) = \begin{cases} 1\quad&\text{if $(w, v) \in E$},\\
0\quad&\text{otherwise}.
\end{cases}
\]
Then $(T, \A)$ is a tree, i.e. a connected non-directed graph with no loop and $E$ is the collection of its edges.
\endremark

We give some additional properties of the self-similar set $K$.

\prop\label{GPS.prop15}
Define $\SS = S^{\BbN} = \{w_1w_2\ldots | \text{$w_i \in S$ for any $i \in \BbN$}\}$. Then for any $\omega = \omega_1\omega_2\ldots \in \SS$, the intersection $\bigcap_{m \ge 1} K_{\word{\omega}m}$ consists of a single point. Moreover, let $\chi(\omega)$ be the single point. Then $\chi: \SS \to K$ is continuous, surjective and
\[
\chi\circ{\s_w^*} = f_w\circ{\chi}
\]
for any $w \in T$, where $\s_w^*: \SS \to \SS$ is defined by $\s_w^*(\omega) = w\omega$ for any $\omega \in \SS$. In particular, the family $\{K_w\}_{w \in T}$ is a partition of $K$ parametrized by the rooted tree $(T, \A, \phi)$ in the sense of {\rm \cite[Definition~2.2.1]{GAMS}} and {\rm\cite[Definition~2.3]{Ki22}}. Moreover,
\begin{equation}\label{GPS.eq10}
\sup_{x \in K}\#(\chi^{-1}(x)) \le 6.
\end{equation}
\endprop

\demo
The first part of the statements follows from \cite[Theorem~1.2.3]{AOF}. The fact that $\{K_w\}_{w \in T}$ is a partition is straightforward from the definition of a partition. Since 
\begin{equation}\label{GPS.eq15}
\#(\{w| w \in T_n, x \in Q_w\}) \le 6,
\end{equation}
where the value $6$ may be attained when $J = 3$, for any $x \in Q_*$ and $n \ge 0$, we have  \eqref{GPS.eq10}.
\enddemo

\definition\label{GPS.def70}
Define
\begin{multline*}
G_*(S, \{f_s\}_{s \in S}) = \{g | g \in D_J, \text{for any $w \in T$, there exists $v \in T$}\\\text{ such that $|v| = |w|$ and $g(Q_w) = Q_v$}\},
\end{multline*}
For simplicity, if no confusion may occur, we write $G_*$ in place of $G_*(S, \{f_s\}_{s \in S})$.
\enddefinition

\remark
It is easy to see that $g \in G_*$ satisfies the condition (A3). Therefore, replacing $G$ with $G_*$, we still have all the properties (A1) through (A5).
\endremark

The rest of this section is devoted to showing the geometric consequences of the conditions (A1) through (A5). In particular, Theorems~\ref{BAS.thm00} and \ref{GPS.thm10} give counterparts of (A1) through (A5) in the levels of $T_n$ and $K$ respectively.

\thm\label{BAS.thm00}
Let $(S, \{f_s\}_{s \in S}, G)$ be a $(J, G)$-s.s.\,system. For any $n \ge 1$, the following statements {\rm (B2)} through {\rm (B5)} are true.\\
{\rm (B2)}\,\,For any $i \in \BbZ_{J}$, there exist $w \in T_n$ and $j \in \BbZ_{J}$ such that $b_j(w) \subseteq b_i$,\\
{\rm (B3)}\,\,For any $g \in G$ and $n \ge 1$, there exists a map $g_*: T_n \to T_n$ such that $(f_{g_*(w)})^{-1}{\circ}g{\circ}f_w \in G$ for any $w \in T_n$. In particular $g(Q_w) = Q_{g_*(w)}$.\\
{\rm (B4)}\,\,Let $w \neq v \in T_n$. If $Q_w \cap Q_{v} \neq \emptyset$, then either {\rm (a)} or {\rm(b)} below holds.\\
{\rm (a)}\,\,$Q_w \cap Q_{v} = b_i(w) = b_{j}(v)$ for some $i, j \in \BbZ_{J}$. Moreover, let $R_{w, v}$ be the reflection in the line segment $b_i(w)$ amd let $g = (f_v)^{-1}{\circ}R_{w, v}{\circ}f_w$. Then $g \in G$ and $R_{w, v}(Q_{wu}) = Q_{vg_*(u)}$ for any $u \in T$.\\
{\rm (b)}\,\,$Q_w \cap Q_v = f_w(p_i) = f_v(p_j)$ for some $i, j \in \BbZ_J$.\\
{\rm (B5)}\,\,For $n \ge 1$, define
\[
E_n^{\ell} = \{(w, v)| w, v \in T_n, \text{{\rm (B4)-(a)} holds.}\}
\]
Then $(T_n, E_n^{\ell})$ is a connected graph.
\endthm

\remark
As in the case of $g_*: S \to S$ in (A3), we omit the subscript $*$ in $g_*$ in (B3) and use $g$ to denote $g_*$.
\endremark

\remark
By (B5), $\cup_{w \in T_m} Q_w$ is a connected subset of $Q_*$ for any $m \ge 1$.
\endremark

\thm\label{GPS.thm10}
Let $(S, \{f_s\}_{s \in S}, G)$ be a $(J, G)$-s.s.\,system. 
The following statements {\rm (C2)}, {\rm (C3)}, {\rm (C4)} and {\rm (C5)} are true:\\
{\rm (C2)}\,\,For any $i \in \BbZ_J$, $K \cap b_i \neq \emptyset$.\\
{\rm (C3)}\,\,For any $g \in G_*$ and $w \in T$, $g(K_w) = K_{g_*(w)}$.\\
{\rm (C4)}\,\,If $w, v \in T_n$ for some $n \ge 1$ and $Q_w \cap Q_v = b_i(w) = b_j(v)$ for some $i, j \in \BbZ_J$, then $K_w \cap K_v \neq \emptyset$  and $R_{w, v}(K_{wu}) = K_{vg_*(u)}$ for any $u \in T$, where $g = (f_v)^{-1}{\circ}R_{w, v}{\circ}f_w \in G$.\\
{\rm (C5)}\,\,$K$ is connected.
\endthm

Before we give proofs of these theorems, we present examples of $(J, G)$-s.s.\,systems. In each of the following examples, we show a configuration of small $J$-gons in the most left figure. We define $S$ as the collection of those small $J$-gons in the configuration. Consequently, it illustrates the first step of iteration, i.e. $\bigcup_{s \in S} Q_s$. Moreover, the contraction ratio $r$ and the centres $\{c_s\}_{s \in S}$ are uniquely determined by the configuration. 

\begin{figure}[ht]
\centering
\includegraphics[width = 300pt]{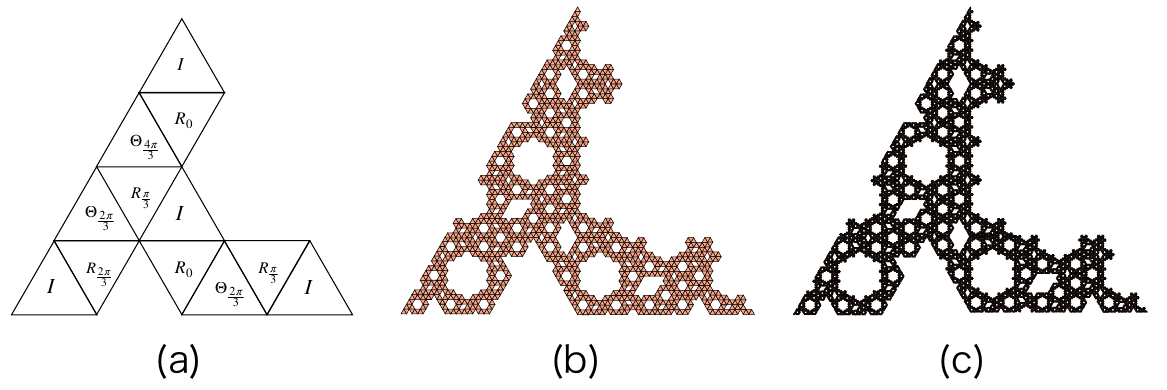}
\caption{$J = 3$, $G = \{I\}$}\label{tri00}
\end{figure}

\example\label{GPS.ex10}
Let $J = 3$. Define $S$ as the collection of small triangles in Figure~\ref{tri00}-(a). In this case, it is easy to see that the contraction ratio $r = \frac 14$ by the configuration of $s \in S$. The map  $\vp_s \in D_3^*$ is indicated in the corresponding $s \in S$ in Figure~\ref{tri00}-(a). Then $(S, \{f_s\}_{s \in S}, G)$ is a $(3, G)$-s.s.\,system with $G = \{I\}$. The middle figure (b) represents the second stage $\bigcup_{w \in T_2} Q_w$ and the right one (c) represents the corresponding self-similar set $K$. We will revisit this example as Example~\ref{BAS.ex10}, where the $p$-conductive homogeneity of $K$ is shown for any $p > \dim_{AR}(K, d_*)$.

\endexample

\begin{figure}[ht]
\centering
\includegraphics[width = 300pt]{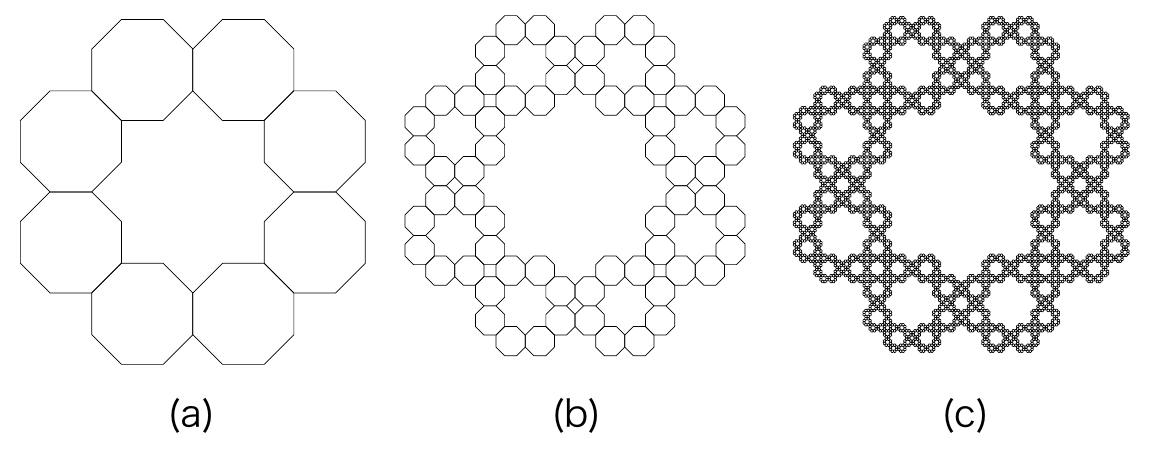}
\caption{Octa-carpet: $J = 8$, $G = D_4$}\label{OctaC}
\end{figure}

\example[Octa-carpet]\label{GPS.ex13}
Let $J = 8$ and let $S$ be the collection of small octagons in Figure~\ref{OctaC}-(a). The contraction ratio $r$ and the centres $\{c_s\}_{s \in S}$ are determined by the configuration of the octagons $s \in S$. Moreover, define $\vp_s = I$ for any $s \in S$. Then $(S, \{f_s\}_{s \in S}, G)$ is a $(8, G)$-s.s.\,systems with $G = D_8$. The associated self-similar set $K$ is called the Octa-gasket in \cite{BeHeStr} or Octa-carpet in \cite{CaHaHu}. We will revisit this example as Example~\ref{BAS.ex15}, where the $p$-conductively homogeneity of $K$ is shown for any $p > \dim_{AR}(K, d_*)$.

\endexample

\begin{figure}[ht]
\centering
\includegraphics[width = 350pt]{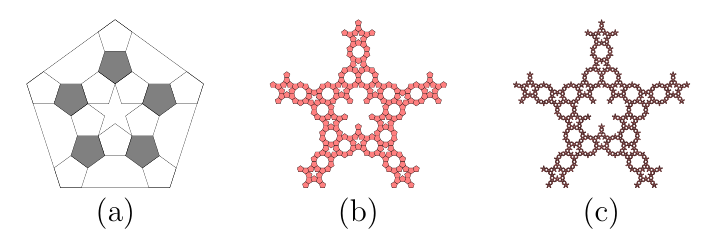}
\caption{$J = 5$, $G = D_5$}\label{PentaS}
\end{figure}

\example\label{GPS.ex15}
Let $J = 5$. Define $S$ as the collection of small white and grey pentagons in Figure~\ref{PentaS}. Define
\[
\vp_s = \begin{cases}
I\quad&\text{if the corresponding octagon $s$ is grey,}\\
R_0\quad&\text{if the corresponding octagon $s$ is white.}
\end{cases}
\]
As in the former examples, the contraction ratio $r$  and the centres $\{c_s\}_{s \in S}$ are determined by the configuration of the pentagons $s \in S$ in Figure~\ref{PentaS}-(a). Then $(S, \{f_s\}_{s \in S}, G)$ is a $(5, G)$-s.s.\,system with $G = D_5$. The middle figure (b) represents the second stage $\bigcup_{w \in T_2} Q_w$ and the right one (c) represents the corresponding self-similar set $K$. We will revisit this example as Example~\ref{BAS.ex20}, where the $p$-conductive homogeneity of $K$ is shown for any $p > \dim_{AR}(K, d_*)$.

\endexample

\example\label{GPS.ex20}
Let $J = 8$. Define $S$ as the collection of small grey and white octagons in Figure~\ref{OctaRot4}-(a). Define 
\[
\vp_s = \begin{cases}
I\quad&\text{if the corresponding octagon $s$ is grey,}\\
R_0\quad&\text{if the corresponding octagon $s$ is white.}
\end{cases}
\]
Then $(S, \{f_s\}_{s \in S}, G)$ is a $(8, G)$-s.s.\,system with $G = Rot_4$. As in the previous example, the middle figure (b) represents the second stage $\bigcup_{w \in T_2} Q_w$ and the right one (c) represents the corresponding self-similar set $K$. We will revisit this example as Example~\ref{MIS.ex20}, where $K$ is shown to be $p$-conductively homogeneous for any $p > \dim_{AR}(K, d_*)$.
\endexample

The rest of this section is devoted to proofs of Theorems~\ref{BAS.thm10} and \ref{GPS.thm10} along with basic properties of a $(J, G)$-s.s.\,systems.

\lemma\label{GPS.lemma10}
Let $w, v \in T_n$ for some $n \ge 1$ and let $g \in D_J$. Then the following three conditions are equivalent:\\
{\rm (1)}\,\,$g(Q_w) = Q_v$,\\
{\rm (2)}\,\,$(\vp_v)^{-1}\circ{g}\circ\vp_w \in D_J$ and $g(c_w) = c_{v}$.\\
{\rm (3)}\,\,$(f_v)^{-1}{\circ}g{\circ}f_w \in D_J$.\\
Moreover, if either of the above conditions holds, then 
\begin{equation}\label{GPS.eq100}
(f_v)^{-1}{\circ}g{\circ}f_w = (\vp_v)^{-1}{\circ}g{\circ}\vp_w.
\end{equation}
\endlemma

\demo
Since $Q_w = f_w(Q)$ and $Q_v = f_v(Q)$, the equivalence between (1) and (3) is obvious. Set $n = |w|$. Note that
\[
g(Q_w) = r^ng{\circ}\vp_w(Q_*) +g(c_w)
\]
and
\[
Q_v = r^n\vp_v(Q_*) + c_v.
\]
It follows that (1) is equivalent to $g{\circ}\vp_w(Q_*) = \vp_v(Q_*)$ and $g(c_w) = c_v$. Thus, we have the equivalence between (1) and (2). The equality \eqref{GPS.eq100} is immediate from the previous arguments.
\enddemo

\lemma\label{GPS.lemma20}
Define $O_w = \inte{Q_w}$ for $w \in T$. \\
{\rm (1)}\,\,For any $n \ge 0$ and $w, v \in T_n$, if $w \neq v$, then $O_w \cap O_v = \emptyset$.\\
{\rm (2)}\,\,For any $w, u \in T$, $Q_u \subseteq Q_w$ if and only if $u \in S^m(w)$ for some $m \ge 0$.
\endlemma

\demo
(1)\,\,Note that $O_w = f_w(O)$, where $O = \inte{Q_*}$. Let $w = \word wn$ and let $v = \word vn$. Since $w \neq v$, there exists $k \in \{1, \ldots, n\}$ such that $w_i = v_i$ for any $i = 1, \ldots, k - 1$ and $w_k \neq v_k$. Then $O_w \cap O_v \subseteq  f_{\word w{k - 1}}(O_{w_k} \cap O_{v_k})$. As is shown in the proof of Proposition~\ref{GPS.prop10}, we see that $O_{w_k} \cap O_{v_k} = \emptyset$. This shows $O_w \cap O_v = \emptyset$.\\
(2)\,\,Assume that $Q_u \subseteq Q_w$. Comparing the areas of $Q_u$ and $Q_w$, we see that $|u| \ge |w|$. So there exists $v \in T_{|w|}$ such that $u \in S^{|u| - |w|}(v)$. Suppose $w \neq v$. Then by (1), $O_w \cap O_v = \emptyset$. Since $O_u \subseteq O_v$, we have $O_u \cap O_w = \emptyset$. On the other hand, the assumption shows $O_u \subseteq O_w$. This contradiction yields $w = u$. The converse direction is obvious.
\enddemo

\prop\label{GPS.prop20}
$g \in G_*$ if and only if there exists  a bijective map $g_*: T \to T$ such that, for any $n \ge  0$ and $w \in T_n$, $g_*(w) \in T_n$ and $\vp_{g_*(w)}^{-1}{\circ}g{\circ}\vp_w \in D_J$ and $g(c_w) = c_{g_*(w)}$. In particular, $G_*$ is a subgroup of $D_J$. Moreover, if $g \in G_*$, then $g_*(S^m(w)) = S^m(g_*(w))$ for any $g \in G_*$ and $w \in T$.
\endprop
\demo
Assume that $g \in G_*$. Then for any $n \ge 0$ and $w \in T_n$, there exists $v \in T_n$ such that $g(Q_w) = Q_v$. Define $g_*: T \to T$ by $g_*(w) = v$. If $g_*(w) = g_*(w')$, then we have $g(Q_w) = g(Q_{w'})$, so that $Q_w = Q_{w'}$. Lemma~\ref{GPS.lemma20}-(1) shows that $g_*$ is injective. Since $T_n$ is a finite set, it is surjective as well. By Lemma~\ref{GPS.lemma10}, it follows that $\vp_{g_*(w)}^{-1}{\circ}g{\circ}\vp_w \in D_J$ and $g(c_w) = c_{g_*(w)}$. The converse direction is obvious.\par
It is straightforward to verify that $G_*$ is a group. Let $u \in S^m(w)$, then $g(Q_w) \supseteq g(Q_u)$. Hence $Q_{g_*(w)} \supseteq Q_{g_*(u)}$.
\enddemo
\lemma\label{BAS.lemma10}
For any $g \in G$, there exists $g_*: T \to T$ such that $g_*(T_m) = T_m$ for any $m \ge 0$ and
\[
(\vp_{g_*(w)})^{-1}{\circ}g{\circ}\vp_{w} \in G,\,\,g(c_w) = c_{g_*(w)} \,\,\text{and}\,\,g(Q_w) = Q_{g_*(w)}
\]
for any $w \in T$. In particular, 
\begin{equation}\label{BAS.eq100}
(f_{g_*(w)})^{-1}\circ{g}\circ{f_w} \in G
\end{equation}
for any $w \in T$ and $g \in G$ and $G \subseteq G_*$. Furthermore, for any $g \in G$ and $w \in T$, there exists $h \in G$ such that 
\begin{equation}\label{BAS.eq700}
g(wv) = g(w)h(v)
\end{equation}
for any $v \in T$.
\endlemma
\demo
For simplicity, we denote $g_*$ by $g$ for $g \in G$.
Let $w = \word wm \in T_m$. We define $g_0, g_1, \ldots, g_m \in G$ inductively by
\[
g_0 = g
\]
and
\[
g_{i + 1} = (\vp_{g_i(w_{i + 1})})^{-1}{\circ}g_i{\circ}\vp_{w_{i + 1}}.
\]
Then by \eqref{GPS.eq100} and (A3), we see that $g_i \in G$ for any $i = 0, 1, \ldots, m$. Moreover, set $v_i = g_{i - 1}(w_i)$. Then, for any $k \in \{1, \ldots, m\}$,
\[
(\vp_{\word vk})^{-1}{\circ}g{\circ}\vp_{\word wk} = g_k\quad\text{and}\quad g_{k - 1}(c_{w_k}) = c_{v_k}.
\]
\noindent{\bf Claim}\,\, For any $k \in \{1, \ldots, m\}$, 
\[
g(c_{\word wk}) = c_{\word vk}.
\]
We show this claim by induction on $k$. The case $k = 1$ is included in the assumption (A2). Assume that $g(c_{\word w{k - 1}}) = c_{\word v{k - 1}}$. Then
\begin{align*}
g{\circ}f_{\word w{k - 1}}(f_{w_k}(0)) &= g(r^{k - 1}\vp_{\word w{k - 1}}(f_{w_k}(0)) + c_{\word w{k - 1}})\\
&= r^{k - 1}g{\circ}\vp_{\word w{k - 1}}(f_{w_k}(0)) + g(c_{\word w{k - 1}})\\
&= r^{k - 1}\vp_{\word v{k - 1}}{\circ}g_{k - 1}(f_{w_k}(0)) + c_{\word v{k - 1}}\\
&= f_{\word v{k - 1}}(g_{k - 1}(f_{w_k}(0))) = f_{\word v{k - 1}}(f_{v_k}(0)).
\end{align*}
Thus, we have shown the claim. Now
\begin{multline*}
g(Q_w) = g(f_w(Q_*)) = g(r^m\vp_w(Q_*) + c_w) = r^mg{\circ}\vp_w(Q_*) + g(c_w)\\
 = r^m\vp_v{\circ}g_m(Q_*) + c_v = f_v(g_m(Q_*)) = f_v(Q_*) = Q_v.
\end{multline*}
The rest is to show \eqref{BAS.eq700}. By Lemma~\ref{GPS.lemma10} and what we have shown above, there exists $g' \in G$ such that
\[
g{\circ}f_{wv} = f_{g(wv)}{\circ}g'
\]
On the other hand, there exist $h, h' \in G$ such that
\[
g{\circ}f_{wv} = g{\circ}f_w{\circ}f_v = f_{g(w)}{\circ}h{\circ}f_v = f_{g(w)}{\circ}f_{h(v)}\circ{h'}
\]
Thus we see that $K_{g(wv)} = K_{g(w)h(v)}$. Therefore $g(wv) = g(w)h(v)$.
\enddemo

\lemma\label{BAS.lemma110}
Assume that $Q_s \cap Q_{s'} = b_i(s) = b_j(s')$ for some $s, s' \in S$ and $i, j \in \BbZ_J$. Then
\begin{equation}\label{GPS.eq20}
(f_{s'})^{-1}{\circ}R_{s, s'}{\circ}f_s = (\vp_{s'})^{-1}{\circ}\vp_s{\circ}R_{\rho(i)}
\end{equation}
In particular, $(f_{s'})^{-1}{\circ}R_{s, s'}{\circ}f_s \in G$. Moreover, let $g = (f_{s'})^{-1}{\circ}R_{s, s'}{\circ}f_s$. Then $R_{s, s'}(Q_{sw}) = Q_{s'g(w)}$ for any $w \in T$.
\endlemma

\demo
Let $q_i$ be the midpoint of $b_i$. Then $c_{s'} - c_s = 2(f_s(q_i) - c_s)$. Moreover, an affine map $h : \BbR^2 \to \BbR^2$ equals $R_{s, s'}$ if and only if $h(c_s) = c_{s'}$ and $h(x) = x$ for any $x \in b_i(s)$. Let
\[
h(x) = f_s{\circ}R_{\rho(i)}{\circ}(f_s)^{-1}(x) + c_{s'} - c_s.
\]
Then $h(c_s) = f_s{\circ}R_{\rho(i)}(0) + c_{s'} - c_s = c_s + c_{s'} - c_s = c_{s'}$. If $x \in b_i(s)$, then $(f_s)^{-1}(x) \in b_i$ Hence $R_{\rho(i)}{\circ}(f_s)^{-1}(x) = (f_s)^{-1}(x) - 2q_i$. Thus
\[
h(x) = f_s((f_s)^{-1}(x) - 2q_i) + c_{s'} - c_s = x - 2f_s(q_i) + 2c_s + c_{s'} - c_s = x.
\]
Thus, we have shown 
\[
R_{s, s'}(x) = f_s{\circ}R_{\rho(i)}{\circ}(f_s)^{-1}(x) + c_{s'} - c_s.
\]
Note that $(f_{s'})^{-1}(x) = \frac 1r(\vp_{s'})^{-1}(x - c_{s'})$. This and the above equality imply
\begin{align*}
(f_{s'})^{-1}{\circ}R_{s, s'}{\circ}f_s(x) &= (f_{s'})^{-1}\big((f_s(R_{\rho(i)}(x)) + c_{s'} - c_s\big)\\
&= (f_{s'})^{-1}\big(r\vp_s(R_{\rho(i)}(x)) + c_{s'}\big)\\
&= (\vp_{s'})^{-1}(\vp_s(R_{\rho(i)}(x))).
\end{align*}
So, we have obtained \eqref{GPS.eq20}. Hence by (A4)-(a), $(f_{s'})^{-1}\circ{g}\circ{f_s} \in G$. Let $g = (f_{s'})^{-1}{\circ}R_{s, s'}{\circ}f_s$. Then Lemma~\ref{BAS.lemma10} implies that $g \in G_*$. Hence $g(Q_w) = Q_{g(w)}$. This yields $R_{s, s'}(Q_{sw}) = Q_{s'g(w)}$.
\enddemo

\demo[Proof of Theorem~\ref{BAS.thm00}]
(B3) has already been verified in Lemma~\ref{BAS.lemma10}.
We use induction on $n$. For $n = 1$, (B2) and (B5) are immediate by (A2) and (B5) respectively. (B4) follows from (A4) and Lemma~\ref{BAS.lemma110}. Assume that the desired statements hold up to $n$.  \\
(B2)\,\,By (A2), for any $i \in \BbZ_J$, there exists $(\xi(i), s(i)) \in \BbZ_J \times S$ such that $b_{k(i)}(s(i)) \subseteq b_i$. Define $s_n(i) = s(\xi^n(i))$ and $w^{(n)}(i) = s_1(i)\cdots{s_n(i)}$, where $\xi^n$ is the $n$-th iteration of a self-map $\xi: \BbZ_J \to \BbZ_J$. The inductive argument shows that
\[
b_{\xi^{n + 1}(i)}(w^{(n + 1)}(i)) \subseteq b_{\xi^n(i)}(w^{(n)}(i)) \subseteq b_i
\]
for any $n \ge 1$. So, we have (B2).\\
(B4)\,\,Let $x, y \in T_n$ and $s, t \in S$. Set $w = xs$ and $v = yt$. Assume that $Q_{w} \cap Q_{v} \neq \emptyset$. \\
{\bf Case 1}\,\, Suppose that $x = y$. \\
Then $Q_s \cap Q_t \neq \emptyset$. If $Q_s \cap Q_t = \{f_s(p_i)\} = \{f_t(p_j)\}$, then $Q_w \cap Q_v = \{f_w(p_i)\} = \{f_v(p_j)\}$, so that (b) holds. Otherwise (A4) implies that $Q_s \cap Q_t = b_i(s) = b_j(t)$. By Lemma~\ref{BAS.lemma110}, we see that $(f_t)^{-1}{\circ}R_{s, t}{\circ}f_s \in G$. Let $g = (f_t)^{-1}{\circ}R_{s, t}{\circ}f_s$. Since $f_x{\circ}R_{s, t}{\circ}(f_x)^{-1} = R_{xs, xt} = R_{w, v}$, it follows that
\[
(f_v)^{-1}{\circ}R_{w, v}{\circ}f_w = (f_t)^{-1}{\circ}R_{s, t}{\circ}f_s = g \in G.
\]
Moreover, this shows that $R_{w, v}(Q_{wu}) = Q_{vg(u)}$ for any $u \in T$. Thus we have (b).\\
{\bf Case 2}\,\,Suppose that $x \neq y$.\\
Then $Q_x \cap Q_y \neq \emptyset$. By the induction hypothesis, either $Q_x \cap Q_y = \{f_x(p_j)\} = \{f_y(p_j)\}$ for some $i, j \in \BbZ_J$ or $Q_x \cap Q_y = b_i(x) = b_j(y)$ for some $i, j \in \BbZ_J$. Suppose that the former is the case. Then $Q_x \cap Q_y = Q_w \cap Q_v$. Hence $p_i \in Q_s$ and $p_i = f_s(p_k)$ for some $k \in \BbZ_J$. Similarly, $p_j = f_t(p_l)$ for some $l \in \BbZ_J$. Thus $Q_w \cap Q_s = \{f_w(p_k)\} = \{f_v(p_l)\}$, so that (b) holds. Suppose that the latter is the case. Then the induction hypothesis yields that $g = (f_y)^{-1}\circ{R_{x, y}}\circ{f_x} \in G$ and $R_{x, y}(Q_{xs}) = Q_{yg(s)}$.  Now there are three cases (i), (ii) and (iii). \\
(i)\,\,$Q_w$ and $Q_v$ intersect at their vertices. In this case, we have (b).\\
(ii)\,\,$Q_w \cap Q_v$ is a single point, which is a vertex of one and is not a vertex of the other. Without loss of generality, we may assume that $Q_w \cap Q_v = \{f_w(p_k)\} \subseteq \sd{b_l(v)}{\{f_v(p_{l - 1}), f_v(p_l)\}}$ for some $k, l \in \BbZ$. Since $Q_w \cap Q_v \subseteq Q_x \cap Q_y$ and $b_l(v) \subseteq Q_y$, we see that $b_l(v) \subseteq b_j(y) = Q_x \cap Q_y$.  This shows that $Q_{yg(s)} \cap Q_v$ has a non-empty interior. Hence $Q_{yg(s)} = Q_v$. Since $f_w(p_k)$ is a vertex of $Q_w$, it follows that $R_{x, y}(f_w(p_k))$ is a vertex of $R_{x, y}(Q_w) = Q_{yg(s)} = Q_v$. This contradicts the assumption that $Q_w \cap Q_v$ is not a vertex of $Q_v$. Thus, this case never happens.\\
(iii)\,\,$Q_w \cap Q_v$ is a line segment. Since $Q_w \cap Q_v \subseteq Q_x \cap Q_y$, it follows that $R_{x, y}(Q_w) \cap Q_v$ has a non-empty interior. Therefore, $Q_{yg(s)} = R_{x, y}(Q_w) = Q_{v} = Q_{yt}$. Hence $v = yg(s)$. Moreover, there exist $k, l \in \BbZ_J$ such that $Q_w \cap Q_v = b_k(w) = b_l(v)$. Using (A3), we see that
\[
(f_v)^{-1}{\circ}R_{w, v}{\circ}f_w = (f_{g(s)})^{-1}{\circ}(f_y)^{-1}{\circ}R_{x, y}{\circ}f_x{\circ}f_s = (f_{g(s)})^{-1}{\circ}g{\circ}f_s \in G.
\]
Letting $g' = (f_v)^{-1}{\circ}R_{w, v}{\circ}f_w$, we also obtain that $R_{w, v}(Q_{wu}) = Q_{vg'(u)}$ for any $u \in T$. Thus (a) holds. So, we have shown (B4).\\
(B5)\,\,
For any $s \in S$, Two graphs $(S^n(s), E_{n + 1}^{\ell} \cap S^n(s)^2)$ and $(T_n, E_n^{\ell})$ is isomorphic through the bijection $\s: S^n(s) \to T_n$ defined by $\s(s\word sn) = \word sn$. By the induction hypothesis, we deduce that 
\[
\text{$(S^n(s), E_{n + 1}^{\ell} \cap S^n(s)^2)$ is connected.}\tag{A}
\]
\par
 On the other hand let $(s, s') \in E_1^{\ell}$. Then $Q_s \cap Q_{s'} = b_i(s) = b_j(s)$ for some $i, j \in \BbZ_J$. By (B3), there exists $w \in T_n$ such that $b_k(w) \subseteq b_i$ for some $k \in \BbZ_J$. Since $b_k(sw) \subseteq b_i(s)$, (A4)-(a) implies $R_{s, s'}(Q_{sw}) = Q_{s'g(w)}$, where $g = (f_{s'})^{-1}{\circ}R_{s, s'}{\circ}f_s \in G$. This shows $(sw, s'g(w)) \in E_{n + 1}^{\ell}$. Thus
 \[
 \text{If $(s, s') \in E_1^{\ell}$, then $(w, v) \in E_{n + 1}^{\ell}$ for some $(w, v) \in S^n(s) \times S^n(s')$.}\tag{B}
 \]
Combining (A), (B) and (A5), we have (B5).
\enddemo

Next, we proceed to prove Theorem~\ref{GPS.thm10}. 

\lemma\label{GPS.lemma30}
For any $w \in T$, let $Q_w^{(m)} = \cup_{v \in S^m(w)} Q_v$. Then $Q_w^{(m + 1)} \subseteq Q_w^{(m)}$ for any $m \ge 0$ and
\[
K_w = \bigcap_{m \ge 0} Q_w^{(m)}.
\]
\endlemma

\demo
For a subset $A \subseteq \BbR^2$, define $F(A) = \cup_{s \in S} f_s(A)$. Then by \cite[Theorem~1.1.7]{AOF}, $F_n(A) \to K$ as $n \to \infty$ in the sense of the Hausdorff metric for any non-empty compact set $A$. Since $f_s(Q_*) \subseteq Q_*$ for any $s \in S$, it follows that $F(Q_*) \subseteq Q_*$. This implies that $F^{n + 1}(Q_*) \subseteq F^n(Q_*)$ and hence $K = \lim_{n \to \infty} F^n(Q_*) = \cap_{n \ge 0} F^n(Q_*)$. Since $F^n(Q) = Q_{\phi}^{(n)}$, we have the desired result for $w = \phi$. For a general $w \in T$,
\[
K_w = f_w(K) = f_w\bigg(\bigcap_{n \ge 0} Q_{\phi}^{(m)}\bigg) = \bigcap_{n \ge 0} f_w(Q_{\phi}^{(n)}) = \bigcap_{n \ge 0} Q^{(n)}_w.
\]
\enddemo

\demo[Proof of Theorem~\ref{GPS.thm10}]
(C2)\,\, We use the same notation as in the proof of Lemma~\ref{BAS.thm00}-(B2). For simplicity, we use $w^{(n)}$ in place of $w^{(n)}(i)$. Then 
\[
K_{w^{(n)}} \subseteq Q_{w^{(n)}}, Q_{w^{(n + 1)}} \subseteq Q_{w^{(n)}}\,\,\text{and}\,\,Q_{w^{(n)}} \cap b_i \neq \emptyset
\]
for any $n \ge 0$. Since $\diam{Q_{w^{(n)}},d_*} \le r^n$, we see that $\bigcap_{n \ge 0} Q_{w^{(n)}}$ is a single point and
\[
\bigcap_{n \ge 0} K_{w^{(n)}} = \bigcap_{n \ge 0} Q_{w^{(n)}} \subseteq K \cap b_i.
\]
(C3)\,\,For any $w \in T$, by Proposition~\ref{GPS.prop20},
\[
g(Q_w^{(m)}) = \bigcup_{v \in T_m} g(Q_{wv}) = \bigcup_{v \in T_m} Q_{g(wv)} = \bigcup_{u \in T_m} Q_{g(w)u} = Q_{g(w)}^{(m)}.
\]
Using Lemma~\ref{GPS.lemma30}, we see
\[
g(K_w) = g\bigg(\bigcap_{m \ge 0} Q_w^{(m)}\bigg) = \bigcap_{m \ge 0} g(Q_w^{(m)}) = \bigcap_{m \ge 0}Q^{(m)}_{g(w)} = K_{g(w)}.
\]
(C4)\,\,Let $\{w^{(n)}\}_{n \ge 0}$ be the same sequence as in the proof of (C2) and let $\{p\} = \bigcap_{n \ge 0} Q_{w^{(n)}} \in K \cap b_i$. Note that $f_w(p) \in K_w \cap b_i(w)$. By (B4)-(1), it follows that
\begin{multline*}
f_w(p) = R_{w, v}(f_w(p)) = R_{w, v}\bigg(\bigcap_{n \ge 0} Q_{ww^{(n)}}\bigg) = \bigcap_{n \ge 0} R_{w, v}(Q_{ww^{(n)}}) \\
= \bigcap_{n \ge 0} Q_{vg(w^{(n)})} \in K_v.
\end{multline*}
Thus $f_w(p) \in K_w \cap K_v$.  Next, note that $g(S^m(u)) = S^m(g(u))$ by Proposition~\ref{GPS.prop20}. This along with (B4)-(1) implies that
\begin{multline*}
R_{w, v}(K_{wu}) = R_{w, v}\bigg(\bigcap_{m \ge 0} Q^{(m)}_{wu}\bigg) = \bigcap_{m \ge 0}\bigcup_{x \in T_m} R_{w, v}(Q_{wux}) \\
= \bigcap_{m \ge 0}\bigcup_{x \in T_m} Q_{vg(ux)} = \bigcap_{m \ge 0} Q_{vg(u)}^{(m)} = K_{vg(u)}.
\end{multline*}
(C5)\,\,Define
\begin{equation}\label{GPS.eq200}
E_n^* = \{(w, v)| w, v \in T_n, K_w \cap K_v \neq \emptyset\}.
\end{equation}
Then \cite[Theorem~1.6.2]{AOF}, it follows that $(K, d_*)$ is connected if and only if $(T_1, E_1^*)$ is connected. By (C4), we see that $E_1^{\ell} \subseteq E_1^*$. Hence, (A1) suffices. 
\enddemo

\setcounter{equation}{0}
\section{Regular polygon-based self-similar sets and\\ associated partition}\label{RPB}

As in the last section, we assume that $(S, \{f_s\}_{s \in S}, G)$ is a $(J, G)$-self-similar system, where $J \ge 3$ and $G$ is a subgroup of $D_J$ throughout this section. Recall that by Proposition~\ref{GPS.prop15}, $\{K_w\}_{w \in T}$ is a partition of the self-similar set $K$ associated with $(S, \{f_s\}_{s \in S})$. \par
In this section, we show that the partition $\{K_w\}_{w \in T}$ combined with a particular self-similar measure $\mu_*$ defined in Theorem~\ref{RPB.thm10} fulfills all the prerequisites to start the investigation of the conductive homogeneity of $(K, d_*)$. Namely, we will verify \cite[Assumption~2.15]{Ki22} with $M_* = 1$. \\
To begin with, we introduce necessary definitions and notations. The first one is from graph theory.

\definition\label{RPB.def00}
Let $(V, E)$ be a non-directed graph where $V$ is the vertices and $E$ is the edges. A sequence $(v(1), \ldots, v(k))$ is called a path of $(V, E)$ if $v(1), \ldots, v(k) \in V$ and $(v(i), v(i + 1)) \in E$ for any $i = 1, \ldots, k - 1$. A path $(v(1), \ldots, v(k))$ of $(V, E)$ is called simple if $w(i) \neq w(j)$ whenever $i \neq j$. For a path $\c = (v(1), \ldots, v(k))$ of $(V, E)$, the length of the path $\c$ is defined as $l(\c) = k$. Moreover, a subset $A \subseteq V$ is connected in $(V, E)$ or $(V, E)$-connected if, for any $v_1, v_2 \in A$, there exists a path $(v(1), \ldots, v(k))$ of $(V, E)$ such that $v(1) = v_1$,$v(k) = v_2$ and $v(i) \in A$ for any $i = 1, \ldots, k$. Moreover $A$ is a connected component of $(V, E)$ if and only if $A$ is $(V, E)$-connected and $A = B$ whenever $B$ is $(V, E)$-connected and $A \subseteq B$.
\enddefinition

\definition\label{RPB.def10}
(1)\,\,Define $\pi: T \to T$ by 
\[
\pi(w) = \begin{cases}
\word w{m - 1}&\,\,\text{if $w \in T_m$ for $m \ge 1$ and $w = \word w{m - 1}w_m$,}\\
\phi &\,\,\text{if $w = \phi \in T_0$.}
\end{cases}
\]
(2)\,\,Define $E_n^* = \{(w, v)| w, v \in T_n, w \neq v, K_w \cap K_v \neq \emptyset\}$. For $A \subseteq T_n$, let $E_n^A = \{(w, v)| w, v \in A, (w, v) \in E_n^*\}$. For $w \in A$ and $M \ge 1$, define $\GG_M^A(w)$ by
\begin{multline*}
\GG_M^A(w) = \{v| v \in A, \text{there exists a path $(w(1), \ldots, w(k))$ of $E_n^A$}\\
\text{ such that $w(1) = w, w(k) = v$ and $k \le M + 1$.}\}
\end{multline*}
For simplicity of notations, for $w \in T$, we write $\GG_M(w) = \GG_M^{T_n}(w)$.
For $x \in K$ and $n \ge 0$, define 
\[
U_M(x: n) = \bigcup_{w \in \GG^n(x)}\bigcup_{v \in \GG_M(w)} K_v,
\]
where $\GG^n(x) = \{w| w \in T_n, x \in K_w\}$.
\enddefinition

\remark
By \eqref{GPS.eq15}, 
\[
\#(\GG^n(x)) \le 6
\]
for any $n \ge 0$ and $x \in K$.
\endremark

\remark
For a sequence $\c = (w(1), \ldots, w(k))$, we occasionally abuse a notation by using $\c$ to denote the set $\{v(1), \ldots, v(k)\}$ if no confusion may occur.
\endremark

The next theorem, along with its tags (1), (2A), (2B), (2C), (3), (4) and (5), corresponds to \cite[Assumption~2.15]{Ki22}, which is a prerequisite to apply the theory of conductive homogeneity. 

\thm\label{RPB.thm10}
For any $M \ge 1$, the following statements {\rm (0), (1), (2A), (2B), (2C), (3), (4)}, and {\rm (5)} hold.\\
{\rm (0)}\,\,The partition $\{K_w\}_{w \in T}$ is  minimal and uniformly finite in the sense of {\rm \cite[Definition~2.5]{Ki22}}.\\
{\rm (1)} $K_w$ is connected for any $w \in T$.\\
{\rm (2A)}\,\,
\[
\diam{K_w, d_*} = r^{|w|}.
\]
{\rm (2B)}\,\,There exist $c_1, c_2 > 0$ such that
\[
B_{*}(x, c_1r^n) \subseteq U_M(x : n) \subseteq B_{*}(x, c_2r^n)
\]
for any $x \in K$ and $n \ge 1$, where $B_*(x, r) = \{y | y \in K, d_*(x, y) < r\}$.\\
{\rm (2C)}\,\,
There exists $c > 0$ such that $B_*(x, cr^{|w|}) \subseteq K_w$ for any $w \in T$.\\
{\rm (3)}\,\,
Let $\mu_*$ be the self-similar measure on $K$ with weight $(\frac 1{N_*}, \ldots, \frac 1{N_*})$, where $N_* = \#(S)$. Then, $\mu_*$ is $\a$-Ahlfors regular, i.e. there exist $c_1, c_2 > 0$ such that
\[
c_1t^{\a} \le \mu_*(B_*(x. t)) \le c_2t^{\a}
\]
for any $x \in K$ and $r \in (0, 1]$. Furthermore 
\begin{equation}\label{RPB.eq100}
\mu_*(K_w) = \sum_{v \in S(w)} \mu_*(K_v)
\end{equation}
for any $w \in T$.\\
{\rm (4)}\,\,There exists $M_0 \ge 1$ such that
\[
\GG_M(u) \cap S^k(w) \subseteq \GG_{M_0}^{S^k(w)}(u)
\]
for any $w \in T$, $k \ge 1$ and $u \in S^k(w)$.\\
{\rm (5)}\,\,There exists $k_* \ge 1$ such that 
\begin{equation}\label{RPB.eq200}
\pi^{k_*}(\GG_{M + 1}(w)) \subseteq \GG_M(\pi^{k_*}(w)).
\end{equation}
for any $w \in T$. In particular, if $J \ge 6$, then the above inclusion holds with $k_* = 1$.
\endthm

An immediate corollary of the above theorem is as follows:

\cor\label{RPB.cor10}
Define $T^{(k)} = \bigcup_{n \ge 0} T_{nk}$. Then for any $M \ge 1$, there exists $k \ge 1$ such that {\rm \cite[Assumption~2.15]{Ki22}} holds if $\{K_w\}_{w \in T}$, $M_*$ and $\mu$ are replaced with $\{K_w\}_{w \in T^{(k)}}$, $M$ and $\mu_*$ respectively. In particular, if $J \ge 6$, then $k = 1$.
\endcor

In the rest of this section, we are going to prove the above theorem piece by piece. To begin with, the statements (1) and (2A) follow from (C5) and the remark after Definition~\ref{GPS.def60}, respectively. Next, we prove the uniform finiteness of $\{K_w\}_{w \in T}$.

\prop\label{RPB.prop05}
Define
\[
L_* = \sup_{w \in T} \#(\GG_1(w)).
\]
Then $L_*$ is finite. Namely, the partition $\{K_w\}_{w \in T}$ is uniformly finite in the sense of {\rm \cite[Definition~2.5]{Ki22}}. Furthermore,
\[
\sup_{w \in T} \#(\GG_M(w)) \le L_*^M
\]
for any $M \ge 1$.
\endprop

\demo
Let $w \in T_n$. Then
\[
\bigcup_{v \in \GG_1(w)} Q_v \subseteq B(c_w, 4r^n),
\]
where $B(x, t) = \{y| y \in \BbR^2, |x - y| < t\}$. Hence, letting $\nu_*$ be the Lebesgue measure, we have
\[
16\pi{r^{2n}} = \nu_*(B(c_w, 4r^n)) \ge \nu_*\bigg(\bigcup_{v \in \GG_1(w)} Q_v\bigg) = \#(\GG_1(w))r^{2n}\nu_*(Q_*).
\]
This implies the finiteness of $L_*$. The rest can be deduced by induction on $M$.
\enddemo

\lemma\label{RPB.lemma10}
There exists $c > 0$ such that if $n \ge 1$, $w, v \in T_n$ and $K_w \cap K_v = \emptyset$, then 
\[
d_*(K_w, K_v) \ge cr^{n - 1},
\]
where $d_*(A, B) = \inf_{x \in A, y \in B} d_*(x, y)$.
\endlemma

\demo
First we define $c_0 > 0$ as
\[
c_0 = \min\{d_*(K_{s_1}, K_{s_2})| s_1, s_2 \in S, K_{s_1} \cap K_{s_2} = \emptyset\}.
\]
Second, define
\[
I_1 = \{(i, s_1, s_2)| i \in \BbZ_J, (s_1, s_2) \in S \times S, (K_{s_1} - q_i) \cap R_{\rho(i)}(K_{s_2} - q_i) = \emptyset\},
\]
where $q_i$ is the midpoint of $b_i$, and
\[
c_1 = \min_{(i, s_1, s_2) \in I} d_*(K_{s_1} - q_i, R_{\rho(i)}(K_{s_2} - q_2)).
\]
Third, define
\begin{multline*}
I_2 = \{(i, s_1, s_2, \psi, c)|i \in \BbZ_J, s_1, s_2 \in S, \psi \in D^*_J, c \in \BbR^2, \\Q_J \cap (\psi(Q_J) + c) = \{p_i\}, K_{s_1} \cap (\psi(K_{s_2}) + c) = \emptyset\}
\end{multline*}
and 
\[
c_2 = \min_{(i, s_1, s_2, \psi, c)} d_*(K_{s_1}, \psi(K_{s_2}) + c).
\]
Finally set $c = \min\{c_0, c_1, c_2\}$. Now we are going to show the above lemma by induction on $n$. First in the case $n = 1$, we see that if $s_1, s_2 \in S = T_1$ and $K_{s_1} \cap K_{s_2} = \emptyset$, then $d_*(K_{s_1}, K_{s_2}) \ge c_0 \ge c$. Thus we have the desired statement for $n = 1$.\\
Next, suppose that the desired statement holds for $n$. Assume that $w, v \in T_{n + 1}$ and $K_w \cap K_v = \emptyset$. Then there are four cases (Case 1), (Case 2), (Case 3) and (Case 4) as follows.\\
(Case 1)\,\,$\pi(w) = \pi(v)$.\\
In this case $w = \pi(w)s_1$ and $v = \pi(w)s_2$ for some $s_1, s_2 \in S$ and $K_{s_1} \cap K_{s_2} = \emptyset$. Hence $d_*(K_w, K_v) \ge c_0r^n$.\\
(Case 2)\,\,$\pi(w) \neq \pi(v)$ and $K_{\pi(w)} \cap K_{\pi(v)} = \emptyset$.\\
In this case, $d_*(K_w, K_v) \ge d_*(K_{\pi(w)}, K_{\pi(v)}) \ge cr^m$.\\
(Case 3)\,\,$\pi(w) \neq \pi(v)$, $K_{\pi(w)} \cap K_{\pi(v)} \neq \emptyset$ and $Q_{\pi(w)} \cap Q_{\pi(v)} = b_i(w) = b_i(w)$ for some $i, j \in \BbZ_J$.\\
In this case, by (C4), we see that $R_{\pi(w), \pi(v)}(K_v) = K_{w'}$ for some $w' \in S(w)$. Hence $w = \pi(w)s_1$ and $w' = \pi(w)s_2$ for some $s_1, s_2 \in S$. Define $\vp$ as $\vp(x) = (f_\pi(w))^{-1}(x) - q_i$. Then $\vp(K_w) = K_{s_1} - q_i$ and $\vp(K_v) = R_{\rho(i)}(K_{s_2} - q_i)$. This shows $(i, s_1, s_2) \in I_1$ and hence it follows that $d_*(K_w, K_v) \ge c_1r^n$.\\
(Case 4) \,\,$\pi(w) \neq \pi(v)$, $K_{\pi(w)} \cap K_{\pi(v)} \neq \emptyset$ and $Q_{\pi(w)} \cap Q_{\pi(v)} = f_w(p_i) = f_v(p_j)$ for some $i, j \in \BbZ_J$.\\
In this case, there exist $\psi \in D_J^*$ and $c \in \BbR^2$ such that $(f_{\pi(w)})^{-1}{\circ}f_{\pi(v)}(x) = \psi(x) + c$. Since $(f_{\pi(w)}^{-1}(Q_{\pi(w)}) = Q_J$ and $(f_{\pi(w)})^{-1}(Q_{\pi(v)}) = \vp(Q_J) + c$, we see that$Q_J \cap (\psi(Q_J) + c) = \{p_i\}$. Moreover, choose $s_1$ and $s_2$ such that $w = \pi(w)s_1$ and $v = \pi(v)s_2$. Then $(f_w)^{-1}(K_w) = K_{s_1}$ and $(f_w)^{-1}(K_{\pi(w)s_2}) = \psi(K_{s_1}) + c$. Therefore, $(i, s_1, s_2, \psi, c) \in I_2$ and hence $d_*(K_w, K_v) \ge c_2r^n$.\\
Combining all the cases, we see that the desired statement holds for $n + 1$ as well. 
\enddemo

\demo[Proof of Theorem~\ref{RPB.thm10}-(2B)]
Let $x \in K$ and let $w \in \GG^n(x)$. Suppose that $d_*(x, y) < rc^{n - 1}$, where $c$ is the constant appearing in Lemma~\ref{RPB.lemma10}. If $v \in T_n$, $y \in K_v$ and $K_w \cap K_v = \emptyset$, then the above lemma shows that $d_*(x, y) \ge cr^{n - 1}$. Hence we see that $K_v \cap K_w \neq \emptyset$. This implies that $v \in \GG_1(w)$ and consequently $y \in U_1(x : n)$. Thus we have shown $B_*(x, cr^{n - 1}) \subseteq U_1(x : n) \subseteq U_M(x : n)$ for any $M \ge 1$. \par
Next, let $y \in U_M(x : n)$. Then there exists a path $(w(1), \ldots, w(M + 1))$ of $(T_n, E_n^*)$ such that $x \in K_{w(1)}$ and $y \in K_{w(M + 1)}$. Hence
\[
d_*(x, y) \le \sum_{i = 1}^{M + 1} \diam{K_{w(i)}, d_*} \le (M + 1)r^n.
\]
This yields $U_M(x : n) \subseteq B_*(x, (M + 2)r^n)$. 
\enddemo

The next lemma implies that $\{K_w\}_{w \in T}$ is minimal.

\lemma\label{RPB.lemma20}
For any $w \in T$, $\sd{K_w}{f_w(\partial{Q_*})} \neq \emptyset$, where $\partial{Q_*}$ is the boundary of the regular $J$-gon $Q_*$ with respect to the Euclidean metric. In particular, the partition $\{K_w\}_{w \in T}$ is minimal in the sense of {\rm \cite[Definition~2.2.1]{GAMS}} and {\rm \cite[Definition~2.5]{Ki22}}.
\endlemma

\demo
Note that once there exists $x \in \sd{K}{\partial{Q_*}}$, $f_w(x)$ belongs to $\sd{K_w}{f_w(\partial{Q_*})}$. So, it is enough to show the claim for $w = \emptyset$. Assume $J \ge 4$. Then for any $s \in S$, there exists $j \in \BbZ_J$ such that $b_j(s) \cap Q_* = \emptyset$. By (C2), there exists $x \in b_j \cap K$. Now $f_s(x) \in b_j(s) \cap K$ and $f_s(x) \notin \partial{Q_*}$. So, we have the desired statement for $J \ge 4$. Let $J = 3$. Suppose $K \subseteq \partial{Q_*}$. Then (C2) and (C5) imply that $b_i \subseteq K$ for some $i \in \BbZ_J$. By (A2), there exist $s \in S$ and $j \in \BbZ_J$ such that $b_j(s) \subseteq b_i$. Moreover (A4) and (A5) implies there exists $t \in S$ such that $(s, t) \in E_1^*$. Now by (C4), $R_{s, t}(b_i) \subseteq K$ but $R_{s, t}(b_i)$ is not contained in $\partial{Q_*}$. This contradicts the assumption that $K \subseteq \partial{Q_*}$. \par Now since
\[
B_w = K_w \bigcap \bigg(\bigcup_{v \in \sd{T_{|w|}}{\{w\}}} K_v\bigg) \subseteq f_w(\partial{Q_*}),
\]
$\{K_w\}_{w \in T}$ is minimal.
\enddemo

Another consequence of this lemma is (2C).

\demo[Proof of Theorem~\ref{RPB.thm10}-(2C)]
Let $x \in \sd{K}{\partial{Q_*}}$. Then $B_*(x, c) \subseteq \sd{K}{\partial{Q_*}}$ for some $c > 0$. Let $x_w = f_w(x)$. Then $B_*(x_w, cr^{|w|}) \subseteq K_w$.
\enddemo

The next lemma is used to show (4).

\lemma\label{RPB.lemma30}
Let $n \ge 1$ and let $A \subseteq T_n$. If $v \in T$, $|v| \ge n$ and $K_v \subseteq \cup_{w \in A} K_w$, then there exist $m \ge 0$ amd $w \in A$ such that $v \in S^m(w)$.
\endlemma

\demo
Set $u = \pi^{|v| - n}(v)$. Assume that $u \notin w$. Then for any $w \in A$, (B4) implies that $Q_w \cap Q_v$ is either empty, a point or a line segment. So, it follows that $Q_w \cap (\sd{K_u}{f_u(\partial{Q_*})}) = \emptyset$. Hence $ (\sd{K_u}{f_u(\partial{Q_*})}) \cap (\cup_{w \in A} K_w) = \emptyset$. This contradicts the assumption of the lemma. Therefore, $u  \in A$ and $v \in S^{|v| - n}(u)$.
\enddemo

\demo[Proof of Theorem~\ref{RPB.thm10}-(4)]
Let $c_1$ and $c_2$ be the constants appearing in (2B). Choose $l \in \BbN$ such that $c_2r^l \le c_1$. Let $u = wu'$ and let $v = wv'$. Moreover, let $x \in \sd{K_{u'}}{f_{u'}(\partial{Q_*})}$ and set $y = f_w(x)$. Then $y \in \sd{K_u}{f_w(\partial{Q_*})}$. Suppose $k >  l$. Using (2B), we see that
\begin{multline*}
K_v \subseteq U_M(y : k + |w|) \cap K_w \subseteq B_*(y, c_2r^{k + |w|}) \cap K_w \subseteq B_*(y, c_1r^{k + |w| - l}) \cap K_w\\
f_w(B_*(x, c_1r^{k - l})) \subseteq f_w(U_M(x : k - l)).
\end{multline*}
Hence 
\[
K_{v'} \subseteq U_M(x : k - l).
\]
By Lemma~\ref{RPB.lemma30}, it follows that $v' \in S^l(\GG_M(\pi^l(u')))$. Hence there exists a path $(v(1), \ldots, v(M + 1))$ of $(S^{k - l}(w), E_{|w| + k - l}^{S^{k - l}(w)})$ such that $u \in S^l(v(1))$ and $v \in  S^l(v(M + 1))$. Since $\#(S^l(v(i))) = N^l$ and $(S^l(v(i)), E_{|w| + k}^{S^l(v(i))})$ is connected for any $i = 1, \ldots, M + 1$, we see that $v \in \GG_{(M + 1)N^l}^{S^k(w)}(u)$. In case $k \in \{1, \ldots, l\}$, since $\#(S^k(w)) \le N^l$, it follows that $v \in \GG_{N^l}^{S^k(w)}(u)$.  Combining these two cases, we have the desired statement with $M_0 = (M + 1)N^l$.
\enddemo

The next piece is the statement (3).

\demo[Proof of Theorem~\ref{RPB.thm10}-(3)]
Note that $r^{\a} = \frac 1{N_*}$. By \cite[Theorem~1.4.5]{AOF} and \eqref{GPS.eq10}, we see that $\mu_*(K_w) = (N_*)^{-|w|}$ for any $w \in T$. This immediately yields \eqref{RPB.eq100}. Next, let $t \in (0, 1]$. Then there exist $n \in \BbN$ and $w \in T_n$ such that $r^{n - 1} \le t < r^{n - 2}$ and $x \in K_w$. Since
\[
K_w \subseteq B_*(x, r^{n - 1}) \subseteq B_*(x, t),
\]
it follows that 
\[
r^{2\a}t^{\a} \le r^{n\a} = \mu_*(K_w) \le \mu_*(B_*(x, t)).
\]
Let $c_1$ be the constant appearing in (2B). If $c_1r \le t$, then $B_*(x, t) \le  1 \le \frac {t^{\a}}{(c_1r)^{\a}}$. In case $t < c_1r$, then there exist $m \ge 1$ such that $c_1r^{m + 1} \le t < c_1t^m$. Making use of (2B), we see that
\[
B_*(x, t) \subseteq B_*(x, c_1r^m) \subseteq U_M(x: m).
\]
By the remark after Definition~\ref{RPB.def10} and  Proposition~\ref{RPB.prop05}, we see that
\[
\mu_*(B_*(x, t)) \le \mu_*(U_M(x: m)) \le \#(\GG^m(x))L_*r^{m\a} \le 6L_*r^{m\a} \le \frac{6L_*}{(c_1r)^{\a}}t^{\a}.
\]
\enddemo

Now we come to the final part of the proof.

\demo[Proof of Theorem~\ref{RPB.thm10}-(5)]
By \cite[Proposition~2.16]{Ki22}, the pieces shown above suffice. Suppose that $J \ge 6$ and that $u, v, w \in T_m$ for some $m \ge 1$, $(u, v, w)$ is a path of $(T_m, E_m^*)$ and $\pi(u) \notin \GG_1(\pi(w))$. \\
{\bf Case 1}:\,\,$Q_{\pi(u)} \cap Q_{\pi(v)} = \{f_{\pi(u)}(p_j)\} = \{f_{\pi(v)}(p_k)\}$ for some $j, k \in \BbZ_J$.\\
In this case, letting $q = f_{\pi(u)}(p_j)$, we see that $Q_u \cap Q_v = \{q\}$. Since $Q_{\pi(w)} \cap Q_v \neq \emptyset$, it follows that $q \in Q_{\pi(w)} \cap Q_{\pi(v} \cap Q_{\pi(u)}$. However, since $J \ge 6$, this is impossible.\\
{\bf Case 2};\,\,$Q_{\pi(u)} \cap Q_{\pi(v)} = b_j(u) = b_k(v)$ for some $j, k \in \BbZ_J$. \\
In this case, similar arguments as in Case 1 show that $Q_{\pi(u)} \cap Q_{\pi(v)} \cap Q_{\pi(w)}$ is a single point. This is impossible unless $J = 6$. Assume $J = 6$. Then $Q_{\pi(u)} \cap Q_{\pi(w)} = b_{j'}(u) = b_{k'}(w)$ for some $j', k' \in \BbZ_J$. Then by (C4), $K_{\pi(u)} \cap K_{\pi(w)} \neq \emptyset$ and hence $(\pi(u), \pi(w)) \in E_{m - 1}^*$. This contradicts to the assumption $\pi(u) \notin \GG_1(\pi(w))$.\\
By the above arguments, we have shown
\begin{equation}\label{RPB.eq10}
\pi(\GG_2(w)) \subseteq \GG_1(\pi(w))
\end{equation}
for any $w \in T$. Now, let $M \ge 2$ and let $(w, w(1), \ldots, w(M))$ be a path in $(T_m, E_m^*)$ for some $m \ge 1$. Then 
by \eqref{RPB.eq10}, we see that $\pi(w(2)) \in \GG_1(\pi(w))$ and hence $\pi(w(M)) \in \GG_M(\pi(w))$. Thus we have shown \eqref{RPB.eq200} with $k_* = 1$.
\enddemo

\setcounter{equation}{0}
\section{Isolated contact points and essential boundaries}\label{ICP}

As in the previous section, we assume that $(S, \{f_s\}_{s \in S}, G)$ is a $(J, G)$-self-similar system, where $J \ge 3$ and $G$ is a subgroup of $D_J$ throughout this section. Recall that by Proposition~\ref{GPS.prop15}, $\{K_w\}_{w \in T}$ is a partition of the self-similar set $K$ associated with $(S, \{f_s\}_{s \in S})$. \par
In this section, we study the effect of a local cut point of $K$. Recall that a point $x$ in a topological space $X$ is called a cut point if $\sd{X}{\{x\}}$ is disconnected and if is called a local cut point if $\sd{U}{\{x\}}$ is disconnected for some neighborhood of $x$. In the case of the self-similar set $K$, Lemma~\ref{COP.lemma10}-(1) shows that such a local cut point, which is called an isolated contact point of cells, is a corner of $Q_w$ for some $w \in T$. The final goal of this section is Theorem~\ref{COP.thm10} giving an effective equivalent condition for the nonexistence of such a cut point.

\definition\label{COP.def10}
Let $(S, \{f_s\}_{s \in S}, G)$ be a $(J, G)$-self-similar system with $J \ge 3$ and let $K$ be the associated self-similar set. Define
\[
Q_n(x) = \bigcup_{w \in \GG^n(x)} Q_w
\]
for $n \ge 1$ and $x \in K$. A point $x \in K$ is called an isolated contact point of cells if there exists $n \in \BbN$ such that
\[
Q_n(x) \backslash \{x\}
\]
is disconnected.
\enddefinition

\begin{figure}[ht]
\centering
\includegraphics[width = 300pt]{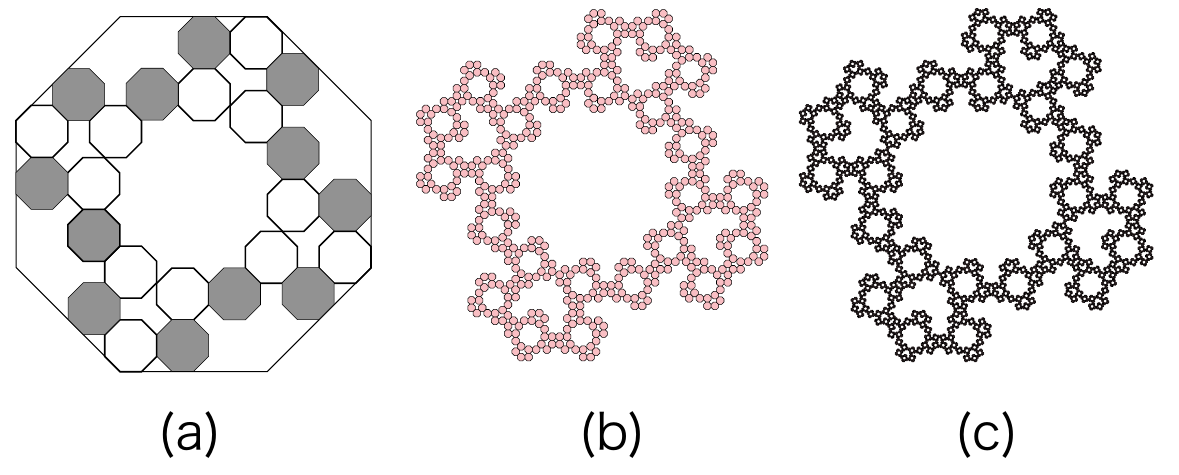}
\caption{$J = 8$, $G = \{D_4^V\}$}\label{Octap}
\end{figure}

\example\label{COP.ex10}
Let $J = 8$ and let $S$ be the collection of small white and grey octagons in Figure~\ref{Octap}-(a). As before, the contraction ratio is determined by the configuration of $s \in S$, and $c_s$ is the centre of the octagon $s$. Define
\[
\vp_s = \begin{cases}
I\quad&\text{if the corresponding octagon $s$ is white,}\\
R_0\quad&\text{if the corresponding octagon $s$ is grey.}
\end{cases}
\]
Then $(S, \{f_s\}_{s \in S}, G)$ is a $(8, G)$-s.s.\,system with $G = D_4^V$, which is generated by $R_{\pi/2}$, and the reflections in lines $p_1p_5, p_2p_6, p_3p_7$ and $p_4p_8$. See Definition~\ref{MTH.def10} for the definition of $D_q^V$ for general $q \ge 3$. The middle figure (b) represents the second stage $\bigcup_{w \in T_2} Q_w$ and the right one (c) represents the corresponding self-similar set $K$. Note that pairs of two white cells intersect at single points in Figure~\ref{Octap}-(a). Let $\{s_1, s_2\}$ be any of such pairs and let $x$ be the single point constituting $Q_{s_1} \cap Q_{s_2}$. Then $x \in K$ and $\sd{Q_1(x)}{\{x\}}$ is disconnected, and hence the point $x$ is an isolated contact point of cells. We will revisit this example as Example~\ref{RNJ.ex20}, where the $p$-conductive homogeneity of $K$ is shown for any $p > \dim_{AR}(K, d_*)$.
\endexample

The next lemma tells us basic properties of an isolated contact point of cells.

\lemma\label{COP.lemma10}
Let $(S, \{f_s\}_{s \in S}, G)$ be a $(J, G)$-self-similar system with $J \ge 3$ and let $K$ be the associated self-similar set. \\
{\rm (1)}\,\,If $x \in K$ and $Q_n(x) \backslash \{x\}$ is disconnected, then, for any $w \in \GG^n(x)$, there exists $i \in \BbZ_J$ such that $x = f_w(p_i)$.\\
{\rm (2)}\,\,A point $x \in K$ is an isolated contact point of cells if and only if the graph $(\GG^n(x), E_n^{\ell}|_{\GG^n(x)})$ is disconnected for some $n \ge 1$.\\
{\rm (3)}\,\,Assume that $J \ge 4$. If $\sd{Q_n(x)}{\{x\}}$ is disconnected, then $\#(\GG^n(x)) = 2$.
\endlemma

\demo
(1)\,\,If $x$ belongs to the interior of $Q_w$, then $\GG^n(x) = \{w\}$ and $\sd{Q_w}{\{x\}}$ is connected. Hence, $x$ belongs to the boundary of $Q_w$. If $x$ does not belong to the vertices of $Q_w$, then $\GG^n(x) = \{w\}$ or $\{w, v\}$ for some $v \in T_n$ and $Q_w \cap Q_v$ is a line segment. In either case, $Q_n(x) \backslash \{x\}$ is connected. This contradiction shows that $x$ is one of the vertices of $Q_w$. This concludes the proof.\\
(2)\,\,Assume that $(\GG^n(x), E_n^{\ell}|_{\GG^n(x)})$ is connected for any $n \ge 1$. Then there exists a path $(w(0), \ldots w(k))$ of $(\GG^n(x), E_n^{\ell}|_{\GG^n(x)})$ such that $\{w(0), \ldots, w(k)\} = \GG^n(x)$. Since $Q_{w(i)} \cap Q_{w(i + 1)}$ is an line segment, we see that $(Q_{w(i)} \cup Q_{w(i + 1)}) \backslash \{x\}$ is connected. Thus $Q_n(x) \backslash \{x\}$ is connected. Therefore, $(\GG^n(x), E_n^{\ell}|_{\GG^n(x)})$ is disconnected for any $n \ge 1$ and hence $x$ is not an isolated contact point of cells. Conversely, assume that $(\GG^n(x), E_n^{\ell}|_{\GG^n(x)})$ is disconnected for some $n \ge 1$. Then, similar arguments as above show that for any $w \in \GG^n(x)$, there exists $i \in \BbZ_J$ such that $x = f_w(p_i)$. Let $C_1, \ldots, C_k$ be the collection of connected components of $(\GG^n(x), E_n^{\ell}|_{\GG^n(x)})$ and define $U_j = \bigcup_{w \in C_j} Q_w$ for $j \in \{1, \ldots, k\}$. If $y \in U_{j_1} \cap U_{j_2}$ for some $j_1 \neq j_2$ and $y \neq x$, then there exist $w_1 \in C_{j_1}$ and $w_2 \in C_{j_2}$ such that $y \in Q_{w_1} \cap Q_{w_2}$. Since $x \in Q_{w_1} \cap Q_{w2}$ as well, we see that $Q_{w_1} \cap Q_{w_2}$ is a line segment and hence $(w_1, w_2) \in E_n^{\ell}$. This contradiction shows that $U_{j_1} \cap U_{j_2} = \{x\}$. Therefore, $\{\sd{U_j}{\{x\}}\}_{j = 1, \ldots, k}$ is the collection of the connected components of $Q_n(x) \backslash \{x\}$. Thus, the point $x$ is an isolated contact point of cells.\\
(4)\,\,Let $\#(Q_n(x)) = k$ and $Q_n(x) = \{w(1), \ldots, w(k)\}$. Then $Q_{w(1)}, \ldots, Q_{w(k - 1)}$ and $Q_{w(k)}$ are distinct $J$-gons sharing the single point $x$, which is one of the vertices of $Q_{w(i)}$ for any $i = 1, \ldots, k$. So, if $J \ge 7$, then $k = 2$. In the case $J = 6$, as we will see in the next proposition, there exists no isolated contact point of cells. Assume that $J = 5$. In this case, we see that $k \le 3$, where $3$ is the maximal number of pentagons sharing a single point. Suppose that $k = 3$. Note that $Q_{w(i)}$ is a parallel translation of either $r^nQ_*^(5)$ or $r^nR_0(Q_*^{(5)}$. So considering all the possible configurations, we conclude that $(w(i_1),w(i_2), w(i_3))$ is a path of $(T_n, E_n^{\ell})$ for some $\{i_1, i_2, i_3\} = \{1, 2, 3\}$. However $\sd{Q_n(x)}{\{x\}}$ is connected in this case. Therefore, $k = 2$ in the case $J = 5$. The remaining case is $J = 4$. In this case, it is immediate to consider all the possible configurations, and we have $k = 2$.
\enddemo

\prop\label{COP.prop10}
Let $(S, \{f_s\}_{s \in S}, G)$ be a $(J, G)$-self-similar system with $J \ge 3$. \\
(1)\,\,If $J = 6$, then there exists no isolated contact point of cells.\\
(2)\,\,If $J \ge 7$, then there exists no isolated contact point of cells if and only if $E_n^* = E_n^{\ell}$ for any $n \ge 1$,.
\endprop

\demo
(1)\,\,Let $J = 6$. Suppose that $x \in K$ is an isolated contact point of cells. By Lemma~\ref{COP.lemma10}, for any $w \in \GG^n(x)$, there exists $i \in \BbZ_J$ such that $x = f_w(p_i)$. By (B5), $\{Q_w\}_{w \in T_n}$ constitute a connected subset of a hexagonal lattice and $x$ is one of its vertex. Hence $\#(\GG^n(x)) \le 3$. In any case, $Q_n(x) \backslash \{x\}$ is connected. This contradiction shows the desired conclusion.\\
(2)\,\,Assume that $J \ge 7$. Then, it follows that $\#(\GG^n(x)) \le 2$ for any $n \ge 1$ and $x \in K$. Suppose that there exists an isolated contact point of cells, say $x \in K$. Then $\GG^n(x) = \{w, v\}$ for some $w \neq v \in T_n$. Moreover, $Q_w \cap Q_v = \{x\}$ and hence $(w, v) \in \sd{E_n^*}{E_n^{\ell}}$. Conversely, suppose that $\sd{E_n^*}{E_n^{\ell}} \neq \emptyset$. Let $(w, v) \in \sd{E_n^*}{E_n^{\ell}}$. Then $Q_w \cap Q_v = K_w \cap K_v$ consists of a single point, and the single point must be an isolated contact point of cells.
\enddemo

To describe an equivalent condition for the nonexistence of isolated contact points of cells, we need the following notion.

\definition\label{COP.def20}
For $g \in D_J$, define $g_*: \BbZ_J \to \BbZ_J$ by $g(b_i) = b_{g_*(i)}$. For simplicity, we sometimes omit the subscript $*$ in $g_*$ and use $g$ to denote $g_*$.
Define
 \begin{multline*}
(\BbZ_J)^e   = \{g_*(i)| g \in G, i \in \BbZ_J, \\
\text{there exist $n \in \BbN$, $w, v \in T_n$ such that $w \neq v$ and $Q_w \cap Q_v = b_i(w)$}\},
\end{multline*}
where the upper case $e$ of $(\BbZ_J)^e$ represents the word ``essential'' and $b_s$ is called an essential boundary segment if $s \in (\BbZ_J)^e$. Moreover, for $B \subseteq \BbZ_J$, define
\[
(B)^e = B \cap (\BbZ_J)^e.
\]
\enddefinition

\definition\label{COP.def30}
For $X \subseteq \BbZ_J$, define
\[
G(X) = \{g(i)| g \in G, i \in X\}.
\]
For $i \in X$, we use $G(i)$ to denote $G(\{i\})$. A subset $X \subseteq \BbZ_J$ is called $G$-invariant if $X = G(X)$ and is called $G$-transitive if $X = G(i)$ for some $i \in X$.
\enddefinition

\remark
A subset $X \subseteq \BbZ_J$ is $G$-transitive if and only if $X = G(i)$ for any $i \in X$. Moreover, if $X$ is $G$-transitive, then it is $G$-invariant.
\endremark

The next proposition gives an effective way to determine $(\BbZ)^e$. 

\prop\label{COP.prop20}
Define
\[
\X = \{X| X \subseteq \BbZ_J, G(X) = X, \text{$X$ satisfies {\rm (E1)} and {\rm (E2)}.}\},
\]
where the conditions {\rm (E1)} and {\rm (E2)} are defined as\\
{\rm (E1)}\,\, For any $(s_1, s_2) \in E_1^{\ell}$, there exist $i_1, i_2 \in X$ such that 
\[
Q_{s_1} \cap Q_{s_2} = b_{i_1}(s_1) = b_{i_2}(s_2).
\]
and\\
{\rm (E2)}\,\,If $i \in \BbZ_J$, $j \in X$, $s \in S$ and $b_i(s) \subseteq b_j$, then $i \in X$.\\
Then
\begin{equation}\label{COP.eq10}
(\BbZ_J)^e = \bigcap_{X \in \X} X
\end{equation}
Moreover, if $X \in \X$ and $X$ is $G$-transitive, then $X = (\BbZ_J)^e$.
\endprop

Now that the conditions (E1) and (E2) concern only $\{b_i(s)\}_{s \in S, i \in X}$, which is a finite set. So, the above proposition is going to be used to identify $(\BbZ_J)^e$ as in the next example, after which we are going to give a proof of the above proposition.\par

\begin{figure}[ht]
\centering
\includegraphics[width = 300pt]{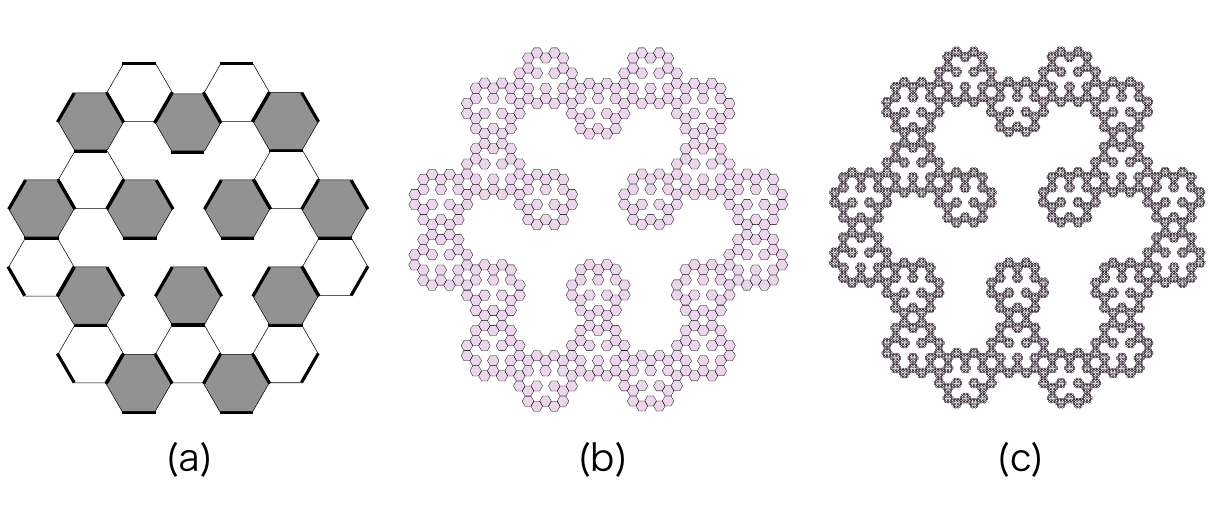}
\caption{$J = 6$, $G = D_3$ and $(\BbZ_6)^e = \{1, 3, 5\}$}\label{HexaZ}
\end{figure}

\example\label{COP.ex20}
Let $J = 6$ and let $S$ be the collection of small white and grey hexagons in Figure~\ref{HexaZ}-(a). As before, the contraction ratio is determined by the configuration of $s \in S$, and $c_s$ is the centre of the octagon $s$. Define
\[
\vp_s = \begin{cases}
I\quad&\text{if the corresponding octagon $s$ is white,}\\
R_0\quad&\text{if the corresponding octagon $s$ is grey.}
\end{cases}
\]
Then $(S, \{f_s\}_{s \in S}, G)$ is a $(6, G)$-s.s.\,system with $G = D_3$. The centre figure, Figure~\ref{HexaZ}-(b), shows the second step and the right figure, Figure~\ref{HexaZ}-(c), illustrate the associated self-similar set $K$. In Figure~\ref{HexaZ}-(a), thick lines are the images of $b_i$ for $i \in \{1, 3, 5\}$. Let $X = \{1, 3, 5\}$. Then it is straightforward to see that $X \in \X$ and $X$ is $D_3$-transitive. By Proposition~\ref{COP.prop20}, it follows that $(\BbZ_6)^e = \{1, 3, 5\}$. We are going to revisit this example as Example~\ref{MIS.ex01}, where $p$-conductive homogeneity of $K$ is shown for any $p > \dim_{AR}(K, d_*)$.
\endexample

\begin{figure}[ht]
\centering
\includegraphics[width = 300pt]{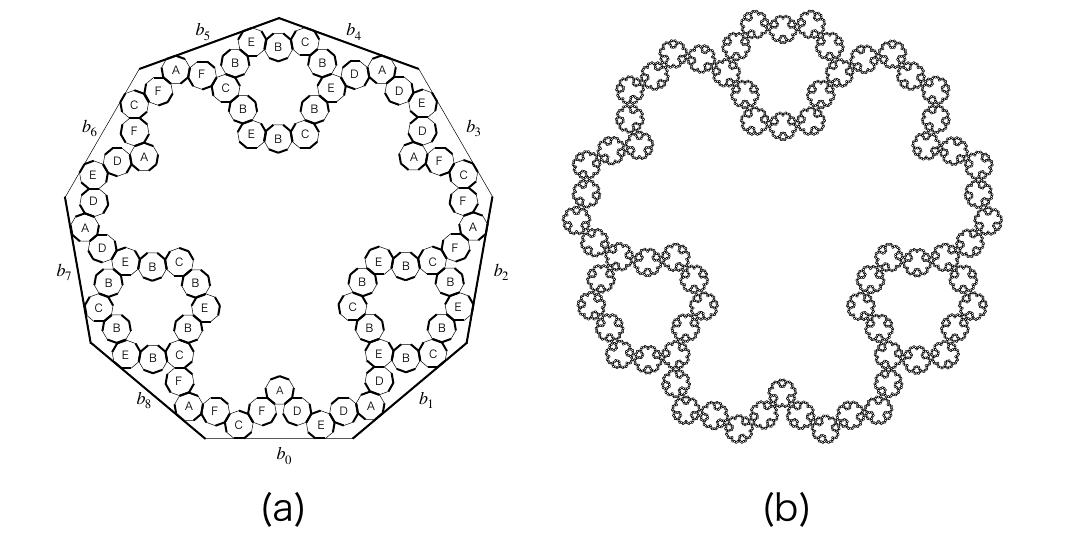}
\caption{$J = 9$, $G = D_3$ and $(\BbZ_9)^e = \{1, 2, 4, 5, 7, 8\}$}\label{nona}
\end{figure}

\example\label{COP.ex30}
Let $J = 9$ and let $S$ be the collection of small nonagons in Figure~\ref{nona}-(a). As before, the contraction ratio is determined by the configuration of $s \in S$, and $c_s$ is the centre of the octagon $s$. The map $\vp_s$ is defined as $I, \Theta_{\pi}, \Theta_{-\frac{2\pi}9}, R_{\frac{2\pi}9}, \Theta_{\frac{2\pi}9}$ and $R_{-\frac{2\pi}9}$ respectively if the corresponding nonagon $s$ in Figure~\ref{nona}-(a) is marked with ``A'', ``B'', ``C'', ``D'', ``E'' and ``F'' respectively. The $(S, \{f_s\}_{s \in S}, G)$ is a $(9, G)$-s.s.\,system with $G = D_3$. The right figure, Figure~\ref{nona}-(b), shows the second step of the generation of the self-similar set $K$. Let $X = \{1, 2, 4, 5, 7, 8\}$. The thick boundary lines of the outer large nonagon correspond to $\{b_i\}_{i \in X}$ and those of small nonagons $s \in S$ correspond to the images of $\{b_i\}_{i \in X}$ by $f_s$. From Figure~\ref{nona}-(a), we see that $X \in \X$ and $X$ is $D_3$-transitive. Thus it follows that $(\BbZ_9)^e = X$. We are going to revisit this example as Example~\ref{MTH.ex100}, where $p$-conductive homogeneity of $K$ is shown for any $p > \dim_{AR}(K, d_*)$.
\endexample

\demo[Proof of Proposition~\ref{COP.prop20}]
First, we are going to show that $(\BbZ_J)^e \in \X$. Suppose that $X = (\BbZ_J)^e$. Then the condition (E1) is immediate from the definition of $(\BbZ_J)^e$. Assume that $i \in \BbZ_J$, $j \in X$, $s \in S$ and $b_i(s) \subseteq b_j$. By the definition of $(\BbZ_J)^e$, there exist $g \in G$, $n \ge 1$ and  $(w, v) \in E_n^{\ell}$ such that $Q_w \cap Q_v = b_{g(j)}(w)$. By (A3), there exists $h \in G$ such that $g{\circ}f_s = f_{g(s)}\circ{h}$. Hence 
\[
b_{h(i)}(g(s)) = f_{g(s)}(h(b_i)) = g(f_s(b_i)) \subseteq g(b_j) = b_{g(j)}.
\]
Since $(\BbZ_J)^e$ is $G$-invariant and $h \in G$, we see that $i \in (\BbZ_J)^e$ if and only if $h(i) \in (\BbZ_J)^e$. So replacing $h(i)$, $g(s)$ and $g(j)$ with $i$, $s$ and $j$ respectively, we may assume that $b_i(s) \subseteq b_j$ and $Q_w \cap Q_v = b_j(w)$ without loss of generality. Let $\xi = (f_v)^{-1}{\circ}R_{w, v}{\circ}f_w$. Then by (B4)-(a), it follows that $\xi \in G$ and $R_{w, v}(Q_{ws}) = Q_{w\xi_*(s)}$. Thus we see that $Q_{ws} \cap Q_{v\xi_*(s)} = b_i(ws)$. Hence $i \in (\BbZ_J)^e$ and we have verified (E2). Thus we have shown that $(\BbZ_J)^e \in \X$.\par
Next we show that if $X \in \X$, then $(\BbZ)^e \subseteq X$. First, using induction, we are going to show that for any $n \ge 1$, the following claims $\rm (E1)_n$ and $\rm (E2)_n$ hold:\,\,
$\rm (E1)_n$\,\,For any $(w, v) \in E_n^{\ell}$, there exists $i_1, i_2 \in X$ such that 
\[
Q_w \cap Q_v = b_{i_1}(w) = b_{i_2}(v).
\]
$\rm (E2)_n$\,\,If $i \in \BbZ_J$, $j \in X$, $w \in T_n$ and $b_i(w) \subseteq b_j$, then $i \in X$.\par
For $n = 1$, they are (E1) and (E2). Assume that ${\rm (E1)_n}$ and ${\rm (E2)_n}$ hold. Let $w = \word w{n + 1} \in T_{n + 1}$ and let $v = \word v{n + 1} \in T_{n + 1}$. Suppose that $(w, v) \in E_{n + 1}^{\ell}$. If $w_1 = v_1$, then $(w', v') \in E_n^{\ell}$, where $w' = w_2{\ldots}w_{n + 1}$ and $v' = v_2{\ldots}v_{n + 1}$. By $\rm (E1)_n$, there exist $i_1, i_2 \in X$ such that $Q_{w'} \cap Q_{v'} = b_{i_1}(w') = b_{i_2}(v')$. Applying $f_{w_1} = f_{v_1}$, we see that $Q_w \cap Q_v = b_{i_1}(w) = b_{i_2}(v)$. In the case $w_1 \neq v_1$, since $(w_1, v_1) \in E_1^{\ell}$, there exist $j_1, j_2 \in X$ such that $Q_{w_1} \cap Q_{w_2} = b_{j_1}(w_1) = b_{j_2}(w_2)$. Suppose that $b_{i_1}(w) = b_{i_2}(v)$ for some $i_1, i_2 \in \BbZ_J$. Since $b_{i_1}(w) \subseteq b_{j_1}(w_1)$, it follows that $b_{i_1}(w') \subseteq b_{j_1}$. The condition $\rm (E2)_n$ shows that $i_1 \in X$. Exactly the same argument shows $i_2 \in X$ as well. Thus we have obtained ${\rm (E1)_{n + 1}}$. Next assume that $i \in \BbZ_J$, $j \in X$, $w \in T_{n + 1}$ and $b_i(w) \subseteq b_j$. Let $w = w_*s$ for some $w* \in T_n$ and $s \in S$. Then $b_i(w) \subseteq b_l(w_*) \subseteq b_j$ for some $l \in \BbZ_J$. Then the condition $\rm (E2)_n$ yields $l \in X$. Since $b_i(s) \subseteq b_l$, we see that $i \in X$ by (E2). Thus we have shown $\rm (E2)_{n + 1}$. \par
So, by induction, we have $\rm (E1)_n$ and $\rm (E2)_n$ for any $n \ge 1$ and the former immediately implies $(\BbZ_J)^e \subseteq X$. Thus, we have shown \eqref{COP.eq10}. \par
Finally, assume that $X \in \X$ and $X$ is $G$-transitive. Above arguments shows $(\BbZ_J)^e \subseteq X$. Therefore, if $i \in (\BbZ_J)^e$, then $G(i) = X$. On the other hand, since $(\BbZ_J)^e$ is $G$-invariant, $G(i) \subseteq (\BbZ_J)^e$. Consequently, $X = (\BbZ_J)^e$.
\enddemo

\lemma\label{PTE.lemma10}
Suppose that there exists no isolated contact point of cells and that $J \ge 4$. Let $i \in (\BbZ_J)^e$. If $u \in T$ and $K_u \cap b_i \neq \emptyset$, then there exists $j \in \BbZ_J$ such that $b_j(u) \subseteq b_i$.
\endlemma

\demo
First assume that there exist $n \ge 1$ and $w, v \in T_n$ such that $Q_w \cap Q_v = b_i(w)$. Let $h = (f_{v})^{-1}{\circ}R_{w, v}{\circ}f_w$. By (B4)-(a), we see that $h \in G$ and $R_{w, v}(Q_{wu}) = Q_{vh(u)}$. Moreover, it follows that
\[
f_w(K_u \cap b_i) = f_v(K_{h(u)} \cap b_{h(i)}) \subseteq Q_{wu} \cap Q_{vh(u)}.
\]
Hence $Q_{wu} \cap Q_{vh(u)} \neq \emptyset$. By (B4), we have either $Q_{wu} \cap Q_{vh(u)} = b_j(wu)$ for some $j \in \BbZ_J$ or $Q_{wu} \cap Q_{vu} = f_{wu}(p_k) = f_{vh(u)}(p_l)$ for some $k, l \in \BbZ_J$. Assume the latter. Then $(wu, vh(u)) \in \sd{E_{n + |u|}^*}{E_{n + |u|}^{\ell}}$. Proposition~\ref{COP.prop10} yields a contradiction if $J \ge 6$ and hence $J = 4$ or $5$. Then it follows that $\#(\GG^{n + |u|}(x)) \le 4$, where $x = f_{wu}(p_k)$. Since $x$ is not a isolated contact point of cells, there exists $u' \in \GG^{n + |w|}(x)$ such that both $(wu, u')$ and $(vg(u), u')$ belong to $E_{n + |w|}^{\ell}$. Let $u' = w'u''$, where $w' \in T_n$ and $u'' \in T_{|u|}$. If $w' \neq w$ and $w' \neq v$, then  $(w, w'), (w', v) \in E_n^{\ell}$, so that $Q_w \cap Q_v$ is a single point. This contradiction shows that $w' = w$ or $w' = v$. In either case, we have $R_{w, v}(Q_{u'}) \neq Q_{u'}$. On the other hand, letting $R_{w, v}(Q_{u'}) = Q_{\tau}$, we see that $(w, \tau), (v, \tau) \in E_{n + |u|}^{\ell}$. Hence $u' = \tau$, so that we end up with a contradiction. Thus, the latter does not happen. Hence the former is the case and $b_j(u) \subseteq b_i$. \par
In general, we can apply the above argument in the case where $i$ and $u$ are replaced by $g_*(i)$ and $g_*(u)$ for some $g_* \in G$. Then we see that there exists $j' \in \BbZ_J$ such that $b_{j'}(g_*(u)) \subseteq b_{g_*(i)}$. Applying $g^{-1}$, we have the desired result.
\enddemo

The next theorem gives a usable criterion for the existence of isolated contact points of cells. In fact, we use it to show the nonexistence of isolated contact points of cells for concrete examples in the later sections.

\thm\label{COP.thm10}
Let $(S, \{f_s\}_{s \in S}, G)$ be a $(J, G)$-self-similar system with $J \ge 4$. There exists no isolated contact point of cells if and only if the following two conditions {\rm (NIC1)} and {\rm (NIC2)} hold;\\
{\rm (NIC1)}\,\, $Q_1(x) \backslash{\{x\}}$ is connected for any $x \in K$\\
{\rm (NIC2)}\,\,There exists $j \in \BbZ_J$ such that $b_j(s) \subseteq b_i$ whenever $K_s \cap b_i \neq \emptyset$ for some $s \in S$ and $i \in (\BbZ_J)^e$.
\endthm

\remark
If $J$ is even, then (NIC2) always holds.
\endremark

\demo
Assume that there exists no isolated contact point of cells. Then the definition of an isolated contact point of cells implies (NIC1). The other property (NIC2) is immediate from Lemma~\ref{PTE.lemma10}. Conversely, assume (NIC1) and (NIC2). Moreover, suppose that there exists an isolated contact point of cells $x \in K$, i.e. $\sd{Q_n(x)}{\{x\}}$ is disconnected for some $n \ge 1$. By Lemma~\ref{COP.lemma10}-(3), we see that $\GG^n(x) = \{w, v\}$ for some $w = \word wn \neq v = \word vn \in T_n$. Let $n_* = \min\{j| j = 1, \ldots, n, w_j \neq v_j\}$ and let $u = \word w{n_* - 1}$, $w' = w_{n_*}\ldots{w_n}$ and $v' = v_{n_*}\ldots{v_n}$. Now replacing $x$, $w$ and $v$ with $f^{-1}(x)$, $w'$ and $v'$ respectively, we may assume that $\sd{Q_n(x)}{\{x\}}$ is disconnected, $\GG^n(x) = \{w, v\}$, $w = \word wn$, $v = \word vn$ and $w_1 \neq v_1$ without loss of generality. Let $l = 1, \ldots, n$. Then $\GG^l(x) = \{\word wl, \word vl\}$ and the intersection of $Q_{\word wl}$ and $Q_{\word vl}$ is the single point $x$ or a line segment. The former is the case when $l = n$ because $\sd{Q_n(x)}{\{x\}}$ is disconnected. In the case $l = 1$, the condition (NIC1) implies that $Q_{w_1} \cap Q_{v_1}$ is a line segment. So, there exists $l_* \in \{1, \ldots, n - 1\}$ such that $Q_{\word w{l_*}} \cap Q_{\word v{l_*}}$ is a line segment and $Q_{\word w{l_* + 1}} \cap Q_{\word v{l_* + 1}}$ is a single point $\{x\}$. Write $w_* = \word w{l_*}$, $v_* = \word v{l_*}$, $s = w_{l_* + 1}$ and $t = v_{l_* + 1}$. Then there exists $i, j \in \BbZ_J$ such that $Q_{w_*} \cap Q_{v_*} = b_i(w_*) = b_j(v_*)$. Suppose that $Q_{w_*s} \cap b_i(w_*)$ is a line segment. Since $R_{w_*, v_*}(Q_{w_*s}) = Q_{v_*t}$, it follows that $Q_{w_*s} \cap b_i(w_*) = Q_{v_*t} \cap b_j(v_*) = Q_{w_*s} \cap Q_{v_*s}$. This contradicts the fact that $Q_{w_*s} \cap Q_{v_*t}$ is a single point $x$. Therefore, we see that $Q_{w_*s} \cap b_i(w_*)$ is a single point $x$. Hence $Q_s \cap b_i = \{(f_{w_*})^{-1}(x)\} \subseteq K$. However, the condition (NIC2) prohibits such a situation. Thus, there exists no isolated contact point of cells.
\enddemo

\setcounter{equation}{0}
\section{Conductive homogeneity and its  consequences}\label{CHC}

In this section, we review the theory of conductive homogeneity in \cite{Ki22} in the case of the self-similar set associated with a $(J, G)$ self-similar system. In particular, we first introduce two constants $\E_{M, p, m}$ and $\s_{p, nm}$, which are called the conductance constant and the neighbour disparity constant, respectively. They are the keys to defining the notion of conductive homogeneity leading to the construction of a counterpart of Sobolev spaces. After giving the definition of the conductive homogeneity, we explain what the consequences of it are, i.e. how a counterpart of Sobolev spaces can be constructed and what their basic properties are. Additionally, in the final part, we present a useful equivalent condition of the conductive homogeneity which will be used in the later sections.\par
Through this section, we continue to assume that $(S, \{f_i\}_{s \in S})$ is a $(J, G)$ self-similar system with $J \ge 3$, $(T, \A)$ is the associated partition, $K$ is the associated self-similar set and $\mu_*$ is the self-similar measure defined in Theorem~\ref{RPB.thm10}. Recall that $(T_n, E_n^A)$ is a non-directed graph for $A \subseteq T_n$.\par
To begin with, we define conductance constants. The following definitions are given in \cite{Ki22}.

\definition\label{CHC.def10}
Let $n\in \BbN$ and let $p \ge 1$. Moreover, let $A \subseteq T_n$. Define $\ell(A) = \{f| f: A \to \BbR\}$. For $f\in \ell (A)$, define $\E_{p, A}^n(f) $
\[
\E_{p, A}^n(f) = \frac{1}{2} \sum_{(u,v)\in E_n^A} |f(u) - f(v)|^p .
\]
Let $A_1, A_2 \subset A$ with $A_1 \cap A_2 = \emptyset$. 
Define
\[
\E_{p,m}(A_1, A_2, A) = \inf\{\E_{p, A}^{n+m}(f) \mid f \colon S^{m}(A) \to \BbR ,\, f\lvert_{S^{m}(A_1)} \equiv 1,\, f\lvert_{S^{m}(A_2)} \equiv 0 \}.
\]
For $w\in T_n$, $M\geq 1$, define
\[\E_{M,p,m}(w) =  \E_{p,m}(w, {\Gamma_M(w)}^c, T_n),\]
which is called the $p$-conductance constant of $w$ with level $m$
\enddefinition

By \cite[Theorem~4.7.6]{GAMS}, it is known that the asymptotic behavior of $p$-conductances $\E_{M, p, m}(w)$ as $m \to \infty$ determines the Ahlfors regular conformal dimension, which is a quasisymmetric invariant of the space defined as follows.

\definition\label{CHC.def100}
Let $(X, d)$ be a metric space.\
(1)\,\,A metric $\rho$ on $X$ is said to be quasisymmetric to $d$ if $(X, \rho)$ gives the same topology as $d$ and there exists a homeomorphism $h$ from $[0, \infty)$ to itself satisfying $h(0) = 0$ and for any $t > 0$, $\rho(x, z) \le h(t)\rho(x, y)$ whenever $d(x, z) < td(x, y)$.\\
(2)\,\,
Let $\a > 0$ and let $\rho$ be a metric on $X$. A Borel regular measure $\mu$ on $(X, \rho)$ is said to be $\a$-Ahlfors regular with respect to $\rho$ if there exist $c_1, c_2 > 0$ such that 
\[
c_1r^{\a} \le \mu(B_{\rho}(x, r)) \le c_2r^{\a}
\]
for any $x \in X$ and $r \in (0, \diam{X, \rho}]$.\\
(3)\,\,
Define the Ahlfors regular conformal dimension $\dim_{AR}(X, d)$ by
\begin{multline*}
\dim_{AR}(X, d) = \inf\{\a| \a > 0, \text{there exist a metric $\rho$ on $X$}\\
\text{and a Borel regular measure $\mu$ on $(X, d)$ such that}\\
\text{$\rho$ is quasisymmetric to $d$ and $\mu$ is $\a$-Ahlfors regular with respect to $\rho$}\}.
\end{multline*}
\enddefinition

The next theorem gives the relation between $\dim_{AR}(K, d_*)$ and the asymptotic behaviours of conductances in the framework of this paper.

\thm\label{CHC.thm10}
For any $M\geq 1$,
\[
\limsup_{m \to \infty}\left( \sup_{w\in T} \E_{M,p,m}(w)^{1/m}\right) <1 \textrm{ if and only if } p> \dim_{AR} (K,d_{*}) .
\] 
\endthm

\demo
By Theorem~\ref{RPB.thm10}, we have \cite[Assumption~2.15]{Ki22}. Hence, by \cite[Proposition~3.3]{Ki22}, we have the desired statement.
\enddemo

Next, we present the definition of neighbour disparity constant, $\s_{p, m, n}$.

\definition\label{CHC.def20}
For $n\in \BbN$, $A\subset T_n$, $f\in \ell(A)$, define
\[
(f)_{A} = \frac{1}{\sum_{v\in A} \mu(K_v)} \sum_{w \in A}f(w)\mu_*(K_w).
\] 
Let $n,m \geq 1$. Define
 \[ 
 \s_{p,m,n} = \sup_{(w,v)\in E^{*}_n}\left( \sup_{f\in \ell(S^m(w,v))} \frac{|(f)_{S^{m}(w)} - (f)_{S^{m}(v)}|^p}{\E_{p, S^{m}(w,v)}^{n+m}(f)}  \right),
 \]
where $S^{m}(w,v)= S^{m}(w) \cup S^{m}(v)$. 
We call $\s_{p,m,n}$ the \textbf{neighbor disparity constant with index $(n, m)$}.
\enddefinition

Now, the $p$-conductive homogeneity is defined as follows.

\definition[Conductive homogeneity]\label{CHC.def30}
Let $p \ge 1$. $(K,d_{*})$ is said to be $p$-conductively homogeneous if and only if  
there exist $c >0$ and $M\geq 1$ such that
\[
\sup_{w \in T}\E_{M,p,m}(w)\sup_{n\geq 1}\s_{p,m,n} \leq  c 
\]
for any $m\geq 1$. 
\enddefinition

The next theorem tells us the reason why we use the word ``homogeneity'' in the last definition.

\thm[{\cite[Theorem~8.1]{Ki22}}]\label{CHC.thm20}
  $(K,d_{*})$ is  $p$-conductively homogeneous if and only if there exist $c_1,c_2 >0$, $\s_p >0$ and $M\geq 1$ such that
\[ 
c_1(\s_p)^{-m} \leq \E_{M,p,m}(w) \leq c_2(\s_p)^{-m}
\]
for any $w \in T\setminus \{\phi\}$, $m \geq 1$ and 
\[
c_1(\s_p)^m \leq \s_{p,m,n} \leq c_2(\s_p)^m
\]
  for any $n,m \geq 1$.
\endthm

Note that Theorem~\ref{CHC.thm10} shows that $\s_p < 1$ if and only if $p > \dim_{AR}(K, d_*)$.\par
The conductive homogeneity enables us to construct a kind of ``Sobolev spaces'' $\W^p$ as a scaling limit of discrete energies $\E_p^m$. Hereafter, for simplicity, we use $\s$ to denote $\s_p$ when no confusion may occur.

\definition\label{CHC.def40}
 Let $C(K)$ be the set of all continuous functions on $K$. 
  For each $f \in L^p(K,\mu)$, $m\geq 1$, define $P_m \colon L^p(K,\mu) \to \ell(T_m)$ by
\[
P_{m}(f)(w) = \frac{1}{\mu(K_w)} \int_{K_w} f(x) \mu(dx)
\]
for  $w\in T_m$ . Moreover, define
\[
\mathcal{N}_p(f) = \left( \sup_{m\geq 0}\s^m \E_{p}^m(P_m f) \right)^{1/p},
\]
and define
\[
\W^p = \{f\mid f\in L^p(K,\mu_*),\,\mathcal{N}_p(f)<\infty \}.
\] 
 \enddefinition
 
The non-negative valued functional $\mathcal{N}_p$ is known to be a semi-norm. The function space $\W^p$ is a counterpart of $(1, p)$-Sobolev space in the case of a smooth space like a domain of $\BbR^n$. The following theorem is extracted from \cite[Theorems~3.21, 3.22, 3.23, 3.35 and 4.6]{Ki22}. It gives basic properties which are analogous to those of the conventional $(1, p)$-Sobolev space of a domain of $\BbR^n$ for $p > n$. Note that
\[
\dim_{AR}(K, d_*) \le \dim_H(K, d_*) < 2.
\]

\thm\label{CHC.thm30}
Suppose that $p > \dim_{AR}(K,d_{*})$ and $(K,d_{*})$ is $p$-conductively homogeneous. Then there exists $\hE_p: \W^p \to \BbR$ having the following properties:\\
{\rm (1)} $\hE_p(f) = 0$ if and only $f$ is a constant on $K$.\\
{\rm (2)} $\hE_p(f)$ is a semi-norm of $\W^p$ and it is equivalent to $\mathcal{N}_p(f)$, i.e. there exist $c_1, c_2 > 0$ such that 
\[
c_1\mathcal{N}_p(f) \le \hE_p(f) \le c_2 \mathcal{N}_p(f)
\]
for any $f \in \W^p$. \\
{\rm (3)} $(\W^p, \norm{\cdot}_{p}+\hE_p)$ is a reflective Banach space, where $\norm{\cdot}_{p}$ is $L^p$ norm with respect to the measure $\mu_*$. \\
{\rm (4)} $\W^p$ is a dense subset of $(C(K), \norm{\cdot}_{\infty})$. Moreover, the inclusion map $\W^p \to C(K)$ is a continuous map from $(\W^p, \norm{\cdot}_{p}+\hE_p)$ to $(C(K), \norm{\cdot}_{\infty})$. \\
{\rm (5)}  $\hE_p$ has Markov property, i.e.
  for any $f \in \W^p$, 
  \[
  \hE_p(\ol{f}) \leq \hE_p(f).
  \] 
where $\ol{f}= (0 \vee f) \wedge 1$. \\
{\rm (6)} There exist $c_3, c_4 > 0$ such that for any $x, y \in K$,
\[
  c_3d_{*}(x,y)^{- \frac{\log \s}{\log r}} \leq \sup_{f\in \W^p, \hE_p(f)\neq 0} \frac{|f(x) - f(y)|^p}{\hE_p(f)} \leq c_4 d_{*}(x,y)^{-\frac{\log \s}{\log r}},
\]
where $r$ is the contraction ratio of the map $f_s$ for $s \in S$.\\
{\rm (7)}{\rm (Self-similarity of energy)} For any $f \in \W^p$ and $s\in S$, $f\circ f_s \in \W^p$ and 
  \[
  \hE_p(f) = \sum_{s \in T_1} \s \hE_p(f\circ f_s).
  \]
  
In particular, in the case $p=2$, $(\hat{\E_2}, \W_2)$ is a strongly local regular Dirichlet form on $L^2(K,\mu)$, and the corresponding diffusion process has the continuous heat kernel $p(t,x,y) :(0,\infty) \times K \times K \to \BbR$. Moreover, there exists $\beta_{*}\geq 2$, $c_1,c_2,c_3, c_4 >0$ such that  
\[
p(t,x,y)\leq c_1t^{- \frac{\alpha_H}{\beta_{*}}}\exp\left( - c_2\left(\frac{|x - y|^{\beta_{*}}}{t} \right)^{\frac{1}{\beta_{*} - 1}} \right)
\]
for any $(t,x,y) \in (0,\infty)\times K \times K$, and 
\[
c_3t^{- \frac{\alpha_H}{\beta_{*}}} \leq p(t,x,y)
\]
for any $y\in B(x,c_4 t^{\frac{\alpha_H}{\beta_{*}}})$, where $\alpha_{H} = \dim_H(K, d_*)$. 
\endthm

Next, we give a useful equivalent condition for $p$-conductive homogeneity for $p > \dim_{AR}(K, d_*)$. 

\thm[{\cite[Theorem~3.33]{Ki22}}]\label{CHC.thm40}
Suppose $p > \dim_{AR}(K,d_{*})$. $(K,d_{*})$ is $p$-conductively homogeneous if and only if 
there exists $M\geq 1$ and $\{c(k)\}_{k \ge 1} \subseteq (0, \infty)$ such that 
\[
 \sup_{z\in T}\E_{M,p,m}(z) \leq c(k)\E_{p,m}(u,v,T_k)
 \]
for any $m\geq 1$ and $u\neq v \in T_k$.
\endthm

This condition is called the ``knight move'' condition. The word ``knight move'' goes back to \cite{BB1}, where Barlow and Bass constructed the Brownian motion of the Sierpinski carpet.

\setcounter{equation}{0}
\section{Conductive homogeneity of $(J, G)$-s.s.\,system I: simple cases}\label{CHP}

In this section, we give three sufficient conditions for the conductive homogeneity of the self-similar set $K$ associated with a $(J, G)$-self-similar system in Theorems~\ref{MTH.thm10}, \ref{MTH.thm20} and \ref{MTH.thm30} and present their examples. Results in these three theorems are simple in the sense that they only concern the number $J$ or the group $G$ and its action on the boundary segments $\{b_i\}_{i \in \BbZ_J}$. Consequently, their statements are simple as well, but proofs require involved arguments that occupy the next three sections. In fact, theorems in this section are essentially consequences of two results, Theorems~\ref{BAS.thm10} and \ref{BAS.thm20}, which will constitute the backbones of all the results on the conductive homogeneity in this paper.\par
As in the previous sections, $(S, \{f_i\}_{s \in S})$ is a $(J, G)$ self-similar system with $J \ge 3$, $(T, \A)$ is the associated partition, $K$ is the associated self-similar set and $\mu_*$ is the self-similar measure defined in Theorem~\ref{RPB.thm10} throughout this section.\par
Now we present the main results of this section, Theorems~\ref{MTH.thm10}, \ref{MTH.thm20} and \ref{MTH.thm30}. Their proofs will be given in Section~\ref{PTM}. The first theorem concerns the case $J = 3$. We don't need any extra conditions in this case at all.

\thm\label{MTH.thm10}
If $J = 3$, then $(K, d_*)$ is $p$-conductively homogeneous for any $p > \dim_{AR}(K, d_*)$.
\endthm

Note that there is no condition on global symmetries, i.e. the group $G$, so that $G$ can be trivial as in the next example.

\example[= Example~\ref{GPS.ex10}]\label{BAS.ex10}
Let $J = 3$ and let $(S, \{f_s\}_{s \in S}, \{I\})$ be the $(3, \{I\})$-s.s.\,system introduced in Example~\ref{GPS.ex10}. Since $J = 3$, the above theorem shows the $p$-conductive homogeneity of $K$ for any $p > \dim_{AR}(K, d_*)$.
\endexample

The second theorem concerns the case where $G$ is relatively large, i.e. the self-similar set has plenty of global symmetry.

\thm\label{MTH.thm20}
If $\BbZ_J$ is $G$-transitive, then $(K, d_*)$ is $p$-conductively homogeneous for any $p > \dim_{AR}(K, d_*)$.
\endthm

If $\BbZ_J$ is $G$-transitive, we have a complete classification of $G$ in the proposition below.

\definition\label{MTH.def10}
Assume that $J$ is even. Define
\[
D_{J/2}^V = Rot_{J/2} \cup \{R_{\theta_i}| i = 0, 1, \ldots, J/2 - 1\},
\]
where $\theta_i = \frac{2\pi}{J}i + \frac{\pi}J - \frac{\pi}2$.
\enddefinition

Recall that $p_i = \frac 1{\cos{\frac{\pi}J}}(\cos{\theta_i}, \sin{\theta_i})$. Hence $R_{\theta_i}$ is a reflection in the line $p_ip_{i + \frac{J}2}$. In fact, the group $D_{J/2}^V$ is isomorphic to $D_{J/2}$ but it does not preserve $Q_*^{(J/2)}$.

\prop\label{MTH.prop10}
$\BbZ_J$ is $G$-transitive if and only if
\[
G = \begin{cases}
D_J,\,\,Rot_J\,\,\text{or}\,\,D_{J/2}^V &\quad\text{if $J$ is even,}\\
D_J\,\,\text{or}\,\,Rot_J&\quad\text{if $J$ is odd.}
\end{cases}
\]
\endprop

To prove the above proposition, we need the following lemma.

\lemma\label{MTH.lemma100}
Assume that $J$ is even. Let $\wG$ be a subgroup of $D_J$. If $\#(\wG) = J$, then $\wG$ is $Rot_J$ or $D_{J/2}$ or $D_{J/2}^V$.
\endlemma

\demo
Let $J = 2q$ for some $q \ge 2$. Define $\wG_{\Theta} = \{g | g \in \wG, \det{g} = 1\}$ and $\wG_R = \{g| g \in \wG, \det{g} = -1\}$. Suppose that $\wG_R \neq \emptyset$. Choose $R \in \wG_{R}$. Define $R^*: \wG \to \wG$ by $R^*(g) = R{\circ}g$. Then the map $R^*$ is one to one and $R^*(\wG_R) = \wG_{\Theta}$. Thus $\#(\wG_R) = \#(\wG_{\Theta}) = q$. This implies that $\wG_{\Theta} = Rot_q$. If $R = R_{\theta_i}$ for some $i \in \{0, \ldots, q - 1\}$, then $\wG = D_q^V$. Otherwise $\wG = D_q$. In the case $\wG_R = \emptyset$, then $\wG = Rot_J$.
\enddemo

\demo[Proof of Proposition~\ref{MTH.prop10}]
If $\BbZ_J$ is $G$-transitive, then $\#(G) = 2J$ or $J$. In the case $\#(G) = 2j$, then $G = D_J$. If $J$ is odd $\#(G) = J$, then the only subgroup of $D_J$ with $\#(G) = J$ is $Rot_J$. If $J$ is even and $\#(G) = J$, then the above lemma gives three possibilities. However,  the case $G = D_{J/2}$ is excluded because $\BbZ_J$ is not $D_{J/2}$-transitive.
\enddemo

\example[= Example~\ref{GPS.ex13}]\label{BAS.ex15}
Let $J = 8$ and let $(S, \{f_s\}_{s \in S}, D_8)$ be the $(8, D_8)$-s.s.\,system introduced in Example~\ref{GPS.ex13}. Since $G = D_8$, the last theorem shows that the $p$-conductive homogeneity of $K$ for any $p > \dim_{AR}(K, d_*)$. In \cite{AndrewsIV}, a diffusion process was constructed on this self-similar set. This result corresponds to the $2$-conductive homogeneity.
\endexample

\example[= Example~\ref{GPS.ex15}]\label{BAS.ex20}
Let $J = 5$ and let $(S, \{f_s\}_{s \in S}, D_5)$ be the $(5, D_5)$-s.s.\,system introduced in Example~\ref{GPS.ex15}. Since $G = D_5$, the last theorem shows the $p$-conductive homogeneity of $K$ for any $p > \dim_{AR}(K, d_*)$.
\endexample

\begin{figure}
\centering
\includegraphics[width = 300pt]{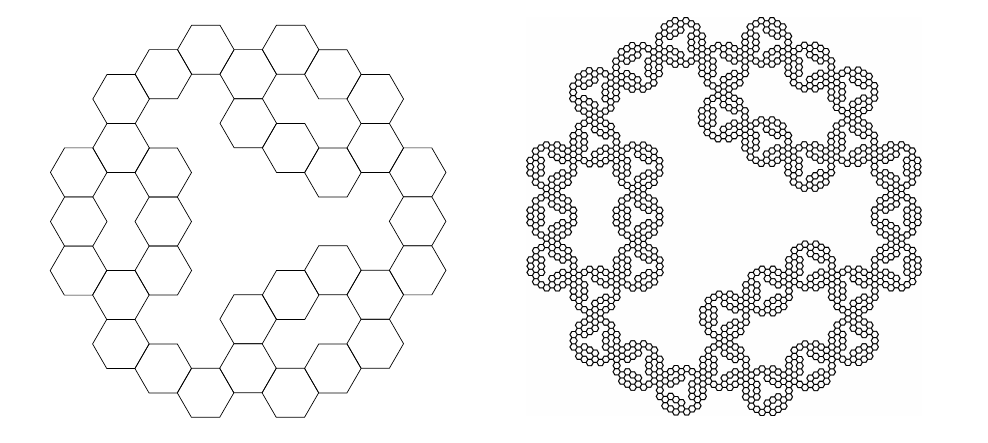}
\caption{$J = 6$ and $G = D_3^V$}\label{HexaDV}
\end{figure}

\example[Figure~\ref{HexaDV}]\label{RNJ.ex10}
Let $J = 6$. As before, the set $S$ is given by the collection of the small hexagons in the left-hand side of Figure~\ref{HexaDV}, and the contraction ratio $r$ and $\{c_s\}_{s \in S}$ are determined accordingly. Moreover, $\vp_s$ is an identity map for every $s \in S$. Then we have a $(6, G)$-s.s.\,system with $G = D_3^V$, so that the last theorem applies.

\endexample

\example[=Example~\ref{COP.ex10}]\label{RNJ.ex20}
Let $J = 8$ and let $(S, \{f_s\}_{s \in S}, D_4^V)$ be the $(8, D_4^V)$-s.s.\,system introduced in Example~\ref{COP.ex10}. By the last theorem, it is $p$-conductively homogeneous for any $p > \dim_{AR}(K, d_*)$ although there are isolated contact points of cells.
\endexample

The third theorem concerns a case where $G$ can be less symmetric, but it is still big enough to generate the essential boundary segments $(\BbZ_J)^e$.

\thm\label{MTH.thm30}
Assume that there exists no contact point of cells and that $(\BbZ_J)^e$ is $G$-transitive. Then $(K, d_*)$ is $p$-conductively homogeneous for any $p > \dim_{AR}(K, d_*)$.
\endthm

Note that if $\BbZ_J$ is $G$-transitive as we assumed in Theorem~\ref{MTH.thm20}, then $(\BbZ_J)^e = \BbZ_J$ and hence the essential boundary segments $(\BbZ_J)^e$ is $G$-transitive as we assume in Theorem~\ref{MTH.thm30}. So, Theorem~\ref{MTH.thm20} and Theorem~\ref{MTH.thm30} are essentially the same, but the latter requires the nonexistence of contact points of cells for a technical reason.

\example[= Example~\ref{COP.ex20}]\label{MIS.ex01}
Let $(S, \{f_s\}_{s \in S}, D_3)$ be the $(6, D_3)$-s.s.\,system given in Example~\ref{COP.ex20}, where we have shown $(\BbZ_6)^e = \BbZ_6^1$. Hence Theorem~\ref{MTH.thm30} implies that the associated self-similar set $K$ is $p$-conductively homogeneous for any $p > \dim{AR}(K, d_*)$.
\endexample

\example[= Example~\ref{COP.ex30}]\label{MTH.ex100}
Let $(S, \{f_s\}_{s \in S}, D_3)$ be the $(9, D_3)$-s.s.\,system given in Example~\ref{COP.ex10}, where we have shown $(\BbZ_9)^e = \{1, 2, 4, 5, 7, 8\}$. Since $(\BbZ_9)^e$ is $D_3$-transitive, Theorem~\ref{MTH.thm30} implies that the associated self-similar set $K$ is $p$-conductively homogeneous for any $p > \dim_{AR}(K, d_*)$.
\endexample

\setcounter{equation}{0}
\section{First backbone theorem}\label{FBT}
Behind Theorems~\ref{MTH.thm10}, \ref{MTH.thm20} and \ref{MTH.thm30}, we have two theorems, Theorems~\ref{BAS.thm10} and \ref{BAS.thm20}, that are the backbones all results on the conductive homogeneity of self-similar sets associated with $(J, G)$-s.s.\,systems. In this section, we present the first one and give its proof.\\
As in the previous sections, $(S, \{f_i\}_{s \in S})$ is a $(J, G)$ self-similar system with $J \ge 3$, $(T, \A)$ is the associated partition, $K$ is the associated self-similar set and $\mu_*$ is the self-similar measure defined in Theorem~\ref{RPB.thm10} throughout this section.

\definition\label{BAS.def10}
(1)\,\,For $A \subseteq T$, define 
\[
K(A) = \bigcup_{w \in A} K_w \quad \text{and}\quad \partial{A} = \{s| s \in \BbZ_J, b_s \cap K(A) \neq \emptyset\}.
\]
(2)\,\, For $n \ge 0$, define $\s_n: T \to  T$ by
\[
\s_n(w_1\ldots{w_k}) = \begin{cases}
w_{n + 1}\ldots{w_k}\quad&\text{if $k \ge n + 1$,}\\
\phi\quad&\text{otherwise.}
\end{cases}
\]
For $u \in T$, $n_1 \in \{0, 1, \ldots, |u|\}$, $n_2, m \ge 0$, let $l = m - n_2 - |u| + n_1$ and define
\[
\H_{n_1, n_2, m}(u) = \begin{cases}
\bigcup_{g \in G. v \in T_{n_2}} S^{l}(vg(\s_{n_1}(u)))\quad&\text{if $l \ge 0$}\\
\{\pi^{-l}(vg(\s_{n_1}(u)))| g \in G, v \in T_{n_2}\}\quad&\text{if $l < 0$}.
\end{cases}
\]
Moreover, for $A \subseteq T$, define
\[
\H_{n_1, n_2, m}(A) = \bigcup_{u \in A} \H_{n_1, n_2, m}(u).
\]
(3)\,\,
For  $w \in T$, $m \ge 0$ and $M \ge 1$, we define $\C_{M, m}(w)$ by
\begin{multline*}
\C_{M, m}(w) = \{(\c(1), \ldots, \c(k))| \c(i) \in S^m(\sd{\GG_M(w)}{\{w\}})\,\,\text{for any $i = 1, \ldots, k$}, \\
\text{there exist $\c(0) \in S^m(w)$ and $\c(k + 1) \in S^m(\GG_M(w)^c)$}\\
\text{such that $(\c(0), \c(1), \ldots, \c(k), \c(k + 1))$ is a path of $(T_{n + m}, E_{n + m}^*)$}\},
\end{multline*}
where $n = |w|$. Moreover, we define $\C_{M, m}^{\ell}(w)$ by the same expression as $\C_{M, m}(w)$ except the graph $(T_{n + m}, E_{n + m}^*)$ is replaced by $(T_{n + m}, E_{n + m}^{\ell})$.\\
\enddefinition

Note that $\H_{n_1, n_2, m}(u) \subseteq T_m$.\par
The next theorem is the first backbone theorem giving a sufficient condition for the conductive homogeneity in terms of families of paths $\C_{M, m}(w)$'s.

\thm\label{BAS.thm10}
Let $M \ge 1$ and let $(j_*, n_*) \in (\BbN \cup \{0\})^2$. For $w \in T$, $m \ge 1$ and $\c \in \C_{M, m}(w)$, define
\begin{multline*}
\P_{\c} = \Bigg\{A \Bigg| A \subseteq \bigcup_{\substack{|j| \le j_*,\\ 0 \le n_2 \le n_*}} \H_{|w| + j, n_2, m}(\c), \\
\,\,\text{$A$ is $(T_m, E_m^*)$-connected and $(\BbZ_J)^e \subseteq \partial{A}$.}\Bigg\}
\end{multline*}
Assume that $\P_{\c} \neq \emptyset$ for any $w \in T$, $m \ge 1$ and $\c \in \C_{M, m}(w)$. Then $(K, d_*)$ is $p$-conductive homogeneous for any $p > \dim_{AR}(K, d_*)$. Moreover, if there exists no isolated contact point of cells, then $\C_{M, m}(w)$ in the above statement is replaced by $\C_{M, m}^{\ell}(w)$.
\endthm

The simplest case of the above theorem is when $j_* = n_* = 0$. In such a case, the set we will deal with is $\H_*(\c)$ defined below.

\definition\label{BAS.def20}
For $\c  \in \C_{M, m}(w)$, define $\H_*(\c)$ by
\[
\H_*(\c) = \H_{|w|, 0, m}(\c)
\]
\enddefinition

The rest of this section is devoted to a proof of Theorem~\ref{BAS.thm10}.

\lemma\label{BAS.lemma500}
For any $u \in T$, $\H_{n_1, n_2, m}(u)$ is $G$-invariant and
\[
\#(\H_{n_1, n_2, m}(u)) \le 
\begin{cases} \#(G)(N_*)^{m - |u| + n_1}\quad&\text{if $l \ge 0$,}\\
\#(G)(N_*)^{n_2}\quad&\text{if $l < 0$}
\end{cases}
\]
\endlemma

\demo
Let $u \in T$ and set $u' = \s_{n_1}(u)$. By Lemma~\ref{BAS.lemma10}, for any $g, h \in G$ and $v \in T_{n_2}$, there exist $h_1, h_2, h_3 \in G$ such that
\[
h{\circ}f_{vg(u')} = f_{h(vg(u'))}{\circ}h_1, \quad h{\circ}f_v = f_{h(v)}{\circ}h_2
\]
and 
\[
h_2{\circ}f_{g(u')} = f_{h_2(g(u'))}\circ{h_3}.
\]
Combining them, we obtain
\[
f_{h(v)}{\circ}f_{h_2(g(u))}\circ{h_3} = f_{h(vg(u'))}\circ{h_1}.
\]
This shows that $h(vg(\s_n(u))) = h(v)(h_2{\circ}g)(\s_n(u))$. Hence $h(\H_{n_1, n_2, m}) \subseteq \H_{n_1, n_2, m}$ for any $h \in G$, so that $\H_{n_1, n_2, m}(u)$ is $G$-invariant. The rest is straightforward by the definition.
\enddemo

\definition\label{BAS.def40}
Let $A_1, A_2 \subseteq T_n$ satisfying $A_1 \cap A_2 = \emptyset$. For $\# \in \{*, \ell\}$, define
\begin{multline*}
\C_m^{\#}(A_1, A_2) = \{(v(1), \ldots, v(k))| v(i) \in S^m(\sd{T_n}{(A_1 \cup A_2)})\,\,\text{for any $i = 1, \ldots, k$}, \\
\text{there exist $v(0) \in S^m(A_1)$ and $v(k + 1) \in S^m(A_2)$}\\
\text{such that $(v(0), v(1), \ldots, v(k), v(k + 1))$ is a path of $(T_{n + m}, E_{n + m}^{\#})$}\}
\end{multline*}
\begin{multline*}
\A_m^{\#}(A_1, A_2) = \{f| f: T_{n + m} \to [0, \infty), \\
\sum_{i = 1}^l f(v(i)) \ge 1\,\,\text{for any $(v(1), \ldots, v(k)) \in \C_m^{\#}(A_1, A_2)$}\}
\end{multline*}
and
\[
\M_{p, m}^{\#}(A_1, A_2) = \inf_{f \in \A_m^{\#}(A_1, A_2)} \sum_{u \in T_{n + m}} |f(u)|^p,
\]
which is called the $p$-modulus of the family of paths $\C_m^{\#}(A_1, A_2)$.
In particular, for $w \in T_n$ and $M \ge 1$, define
\[
\M_{M, p, m}^{\#}(w) = \M_{p, m}^{\#}(\{w\}, \GG_M(w)^c).
\]
If no confusion can occur, we omit $*$ in $\C_m^*(A_1, A_2), \A_m^*(A_1, A_2), \M_{p, m}^{*}(A_1, A_2)$ and $\M_{M, p, m}^*(w)$.
\enddefinition

\remark
The family of paths $\C_{M, m}(w)$ defined in Definition~\ref{BAS.def10} coincides with $\C_m(\{w\}, \GG_M(w)^c)$.
\endremark

Since $\C_m^{\ell}(A_1, A_2) \subseteq \C_m*(A_1, A_2)$, we see that $\A_m^{\ell}(A_1, A-2) \supseteq \A_m^*(A_1, A_2)$, so that
\begin{equation}\label{BAS.eq05}
\M_{p, m}^{\ell}(A_1, A_2) \le \M_{p, m}^*(A_1, A_2).
\end{equation}

\lemma\label{BAS.lemma150}
Assume that $(S, \{f_s\}_{s \in S})$ has no isolated contact point of cells. Then 
\begin{equation}\label{BAS.eq07}
\M_{p, m}^{\ell}(A_1, A_2) \le \M_{p, m}^*(A_1, A_2) \le 12^{p + 1}\M_{p, m}^{\ell}(A_1, A_2)
\end{equation}
for any $n \ge 1$, $A_1, A_2 \subseteq T_n$ with $A_1 \cap A_2 = \emptyset$.
\endlemma

\demo
By Proposition~\ref{COP.prop10}, $\M_{p, m}^{\ell}(A_1, A_2) = \M_{p, m}^*(A_1, A_2)$ if $J \ge 6$. So we may suppose that $J = 3, 4$ or $5$. For $v \in T_{n + m}$, define
\[
H_v = \{u| u \in T_{n + m}, Q_u \cap Q_v \neq \emptyset\}.
\]
Note that $\#(H_v) \le 12$ for any $v \in T_{n + m}$ and $\#(\{v| v \in T_{n + m}, u \in H_v\}) \le 12$ for any $u \in T_{n + m}$.\par
Let $(w(1), \ldots, w(k)) \in \C_m^*(A_1, A_2)$. Choose $w(0) \in S^m(A_1)$ and $w(k + 1) \in S^m(A_2)$ such that $(w(0), w(1)), (w(k), w(k + 1)) \in E_{n + m}^*$. Suppose that $Q_{w(i)} \cap Q_{w(i + 1)}$ is a single point for some $i \in \{0, \ldots, k\}$. Since the single point is not an isolated contact point of cells, there exists a path $(v_1, \ldots, v_l)$ of $(T_{n + m}, E_{n + m}^{\ell})$ such that $v_1 = w(i)$, $v_l = w(i + 1)$ and $\{v_1, \ldots, v_l\} \subseteq H_{w(i)}$. Replace $(w(i), w(i + 1))$ with $(v_1, \ldots, v_l)$ whenever $Q_{w(i)} \cap Q_{w(i + 1)}$ is a single point. Denote the resulting path by $\c_1 = (w'(0), w'(1), \ldots, w'(k'), w'(k' + 1))$, where $w'(0) = w(0)$ and $w'(k' + 1) = w(k + 1)$. Then $\c_1$ is a path of $(T_{n + m}, E_{n + m}^{\ell})$. Let $i_1 = \max\{i| w'(i) \in S^m(A_1)\}$ and let $i_2 = \min\{i| w'(i) \in S^m(A_2)\}$. Then $(w'(i_1 + 1), \ldots, w'(i_2 - 1)) \in \C_m^{\ell}(A_1, A_2)$ and $\{w'(i_1 + 1), \ldots, w'(i_2 - 1)\} \subseteq \bigcup_{i = 1, \ldots, k} H_{w(i)}$. 
So, Using \cite[Lemma~C.4]{Ki22}, we see that
\[
\M_{p, m}^*(A_1, A_2) \le 12^{p + 1}\M_{p, m}^{\ell}(A_1, A_2).
\]
\enddemo

\demo[Proof of Theorem~\ref{BAS.thm10}]
Set $n = |w|$. Let
\[
\H(u) = \bigcup_{j \in \{-j_*, \ldots, 0, \ldots,  j_*\}, n_2 \in \{0 1, \ldots,n_*\}} \H_{|w| + j, n_2, m}(u).
\]
for $u \in T_{n + m}$. Define $H: T_{n + m} \to \P(T_{k + m})$, where $\P(T_{k + m})$ is the collection of subsets of $T_{k + m}$, by
\[
H_u = \bigcup_{x \in T_k} x\H(u)
\]
for $u \in T_{n + m}$. Then 
\begin{multline}\label{BAS.eq10}
\#(H_u) \le \#(T_k)\sum_{j = -j_*}^{j_*} \sum_{n_2 = 0}^{n_*}\#(\H_{n + j, n_2, m}(u)) \\
\le (2j_* + 1)(n_* + 1)(N_*)^{k + j_* + n_*}\#(G)
\end{multline}
for any $u \in T_{n + m}$. On the other hand, let 
\[
\I_{n + j, n_2, m}(v) = \{u| u \in S^m(\GG_M(w)), \s_k(v) \in \H_{n + j, n_2, m}(u)\}
\]
If $j \ge n_2$, then $\pi^{j - n_2}\s_{n_2 + k}(v) = g(\s_{n + j}(u))$ for any $v \in T_{k + m}$ and $u \in \I_{n + j, n_2, m}$. Therefore, 
\[
\I_{n + j, n_2, m}(v) \subseteq \{w_1u_2g(\pi^{j - n_2}\s_{n_2 + k}(v) |  w_1 \in \GG_M(w), u_2 \in T_j\}.
\]
This implies $\#(\I_{n + j, n_2, m}(v)) \le (L_*)^M(N_*)^j$. If $j < n_2$, then for any $v \in T_{k + m}$ and $u \in \I_{n + j, n_2, m}(v)$, there exists $u_2 \in T_{n_2 - j}$ such that $\s_{n_2 + k}(v)u_2 = g(\s_{n + j}(u)$. Therefore, 
\[
\I_{n + j, n_2, m}(v) = \{w_2g(\s_{n_2 + k}(v)u_2)| g \in G, w_2 \in \pi^{n_2 - j}(\GG_M(w)), u_2 \in T_{n_2 - j}\}
\]
So, we have $\#(\I_{n + j, n_2, m}(v)) \le \#(G)(L_*)^M(N_*)^{n_2 - j}$. Combining two cases, we see that
\[
\#(\I_{n + j, n_2, m}(v)) \le \#(G)(L_*)^M(N_*)^{n_* + j_*}.
\]
Thus, for any $w \in T_n$ and $v \in T_{k + m}$,
\begin{multline}\label{BAS.eq20}
\#(\{u| u \in S^m(\GG_M(w)), v \in H_u\}) \le \sum_{j = -j_*}^{j_*}\sum_{n_2 = 0}^{n_*}\#(\I_{n + j, n_2, m})\\
 \le (2j_* + 1)(n_* + 1)\#(G)(L_*)^M(N_*)^{n_* + j_*}
\end{multline}
Now, let $u, v \in T_k$. Then by (B5), there exists a path $(u(0), u(1), \ldots, u(l_*), u(l_* + 1))$ of $(T_k, E_k^{\ell})$ such that $u(0) = u$ and $u(l_* + 1) = v$. Let 
\[
g_i = (f_{u(i + 1)})^{-1}{\circ}R_{u(i), u(i + 1)}{\circ}f_{u(i)}
\]
for $i = 0, 1, \ldots, l_*$. Then by (B4)-(a), it follows that $g_i \in G$ for any $i = 0, 1, \ldots, l_*$. \\
{\bf Claim}\,\,Let $\c \in \C_{M, m}(w)$. Assume that $\P_{\c} \neq \emptyset$. Then there exists $\c' \in \C_m(u, v)$ such that $\c' \subseteq \bigcup_{j = 1, \ldots, l} H_{v(j)}$, where $\c = (v(1), \ldots, v(l))$.\\
Proof of Claim\,\,
Let $A \in \P_{\c}$. For $i = 0, 1, \ldots, l_* + 1$, define $A_i \subseteq T_m$ inductively by
\[
A_0 = A\quad\text{and}\quad A_{i + 1} = g_i(A_i).
\]
Set $h_0 = I \in D_J$ and $h_i = g_{i - 1}{\circ}\cdots{\circ}{g_0}$ for $i = 1, \ldots, l_* + 1$. 
Since $h_i \in G$ and $A_i = h_i(A)$, it follows that
\[
A_i  \subseteq \bigcup_{j = 1, \ldots, l}\H(v(j))
\]
for any $i = 0, 1, \ldots, l_* + 1$.  This shows
\begin{equation}\label{BAS.eq30}
\bigcup_{i = 0, 1, \ldots, l_* + 1} u(i)A_i \subseteq \bigcup_{j = 1, \ldots, l} H_{v(j)}.
\end{equation}
Moreover, using Lemma~\ref{BAS.lemma201}, we see that
\[
R_{u(i), u(i + 1)}^*(u(i)A_i) = u(i + 1)A_{i + 1}
\]
and hence
\begin{equation}\label{BAS.eq40}
R_{u(i), u(i + 1)}(K(u(i)A_i)) = K(u(i + 1)A_{i + 1})
\end{equation}
for any $i = 0, 1, \ldots, l_*$. 
Since $(u(i), u(i + 1)) \in E^{\ell}_k$,  there exist $s_i, t_i \in (\BbZ_J)^e$ such that $Q_{u(i)} \cap Q_{u(i + 1)} = b_{s_i}(u(i)) = b_{t_i}(u(i + 1))$ for any $i = 0, 1, \ldots, l_*$. Since $(h_i)^{-1} \in G$, we see that $(h_i)^{-1}(s_i) \in (\BbZ_J)^e$. Hence
\[
K(A) \cap (h_i)^{-1}(b_{s_i}) = K(A) \cap b_{(h_i)^{-1}(s_i)} \neq \emptyset.
\]
Applying $h_i$ and $f_{u(i)}$, we obtain
\[
K(u(i)A_i) \cap b_{s_i}(u(i)) \neq \emptyset.
\]
Furthermore, by \eqref{BAS.eq40}, 
\begin{multline*}
K(u(i)A_i) \cap b_{s_i}(u(i)) = R_{u(i), u(i + 1)}(K(u(i)A_i) \cap b_{s_i}(u(i)))\\ = K(u(i + 1)A_{i + 1}) \cap b_{t_i}(u(i + 1)).
\end{multline*}
Hence $K(u(i)A_i) \cap K(u(i + 1)A_{i + 1}) \neq \emptyset$ for any $i = 0, \ldots, l_*$. So, there exists $(\omega(1), \ldots, \omega(k_*)) \in \C_m(u, v)$ such that 
\begin{equation}\label{BAS.eq50}
\{\omega(1), \ldots, \omega(k_*)\} \subseteq \bigcup_{i = 1, \ldots, l_*} u(i)A_i\subseteq \bigcup_{j = 1, \ldots, l} H_{v(j)},
\end{equation}
where the second inclusion follows from \eqref{BAS.eq30}. Thus we have shown the claim.\qed\\
Suppose that $\P_{\c} \neq \emptyset$ for any $\c \in \C_{M, m}(w)$. Then since \eqref{BAS.eq10}, \eqref{BAS.eq20} and \eqref{BAS.eq50} constitute the assumptions of \cite[Lemma~C.4]{Ki22}, the above claim yields
\begin{equation}\label{BAS.eq60}
\M_{M, p, m}(w) \le c_k\M_{p, m}(\{u\}, \{v\}),
\end{equation}
where 
\[
c_k = (N_*)^{kp}(L_*)^M\big((2j_* + 1)(n_* + 1)\#(G)(N_*)^{n_* + j_*}\big)^{1 + p}
\]
 Using \cite[Lemma~4.2]{Ki22}, we see that
\begin{equation}\label{BAS.eq70}
\E_{M, p, m}(w) \le 2(L_*)^pc_k\E_{p, m}(\{u\}, \{v\}).
\end{equation}
So, applying Theorem~\ref{CHC.thm40}, we obtain the desired conclusion.\par
Finally in the case where there exists no isolated contact point of cells, if $\P_{\c} \neq \emptyset$ for any $\c \in \C_{M, m}^{\ell}(w)$, then we have
\[
\M_{M, p, m}^{\ell}(w) \le c_k\M_{p, m}(\{u\}, \{v\}).
\]
by \cite[Lemma~C.4]{Ki22}. Using Lemma~\ref{BAS.lemma150}, we obtain modified versions of \eqref{BAS.eq60} and \eqref{BAS.eq70} where the constant $c_k$ is replaced with $12^{p + 1}c_k$. The rest of the argument is the same as above.
\enddemo

\setcounter{equation}{0}
\section{Second backbone theorem}\label{SBT}
In this section, we are going to study detailed structure of paths in $\C_{M, m}(w)$ to apply the first backbone theorem, Theorem~\ref{BAS.thm10}, presented in the last section. \par
The setting of this section is exactly the same as the previous two sections.

\definition\label{BAS.def53}
Let $m \ge n \ge 0$. For $w = \word wm \in T_m$, define 
\[
[w]_n = \word wn.
\]
\enddefinition

In the following lemma and definition, we are going to resolve a path $\c$ of $(T_n, E_n^*)$  into a path of $(T_k, E_k^*)$  and a sequence of paths of  $(T_{n - k}, E_{n - k}^*)$ for $k \in \{0, 1, \ldots, n\}$.

\lemma\label{MIS.lemma15}
Let $\c = (\c(1), \ldots, \c(l))$ be a path of $(T_n, E_n^*)$. For $k \in \{0, \ldots, n\}$, there exist $l_k(\c) \in \BbN$ and an increasing sequence $\{j_1, j_2, \ldots, j_{l_k(\c)}\} \subseteq \BbN$ such that $j_0 = 0$, $j_{l_k(\c)} = l$, 
\[
[\c(j_{j_i - 1} + 1)]_k = [\c(j_{i - 1} + 2)]_k =  \cdots = [\c(j_i)]_k
\]
for every $i = 1, \ldots l_k(\c)$, 
\[
[\c(j_i)]_k \neq [\c(j_i + 1)]_k
\]
for every $i= 1, \ldots, l_k(\c) - 1$.
\endlemma

\definition\label{PTE.def20}
Let $\c = (\c(1), \ldots, \c(l))$ be a path of $(T_n, E_n^*)$. For $k \in \{0, \ldots, n\}$, define
\[
[\c]_k = ([\c(j_1)]_k, \ldots, [\c(j_{l_k(\c)})]_k),
\]
where $l_k(\c)$ and $j_1, \ldots, j_{l_k(\c)}$ are those given in Lemma~\ref{MIS.lemma15}. Moreover, define a family of paths $\{\s_k(\c)_1, \ldots, \s_k(\c)_{l_k(\c)}\}$ of $(T_{n - k}, E_{n - k}^*)$ by
\[
\s_k(\c)_i = (\s_k(\c(j_{i - 1} + 1)), \ldots, \s_k(\c(j_i)))
\]
for $i = 1, \ldots, l_k(\c)$.
\enddefinition

\definition\label{PTE.def23}
Let $\c_1 = (\c_1(1), \ldots, \c_1(l))$ amd $\c_2 = (\c_2(1), \ldots, \c_2(k))$ be paths of $(T_n, E_n^*)$. If $K_{w(\c_1(l))} \cap K_{\c_2(l)} \neq \emptyset$, we define the concatenation of $\c_1$ and $\c_2$, $\c_1 \vee \c_2$, by
\[
\c_1 \vee \c_2 = \begin{cases}
(\c_1(1), \ldots, \c_1(l), \c_2(1), \ldots, \c_2(k))\quad & \text{if $\c_1(l) \neq \c_2(1)$,}\\
(\c_1(1), \ldots, \c_1(l), \c_2(2), \ldots, \c_2(k))\quad & \text{if $\c_1(l) = \c_2(1)$.}
\end{cases}
\]
\enddefinition

Using the above notation, we have
\begin{equation}\label{BAS.eq800}
\c = [\c]_k(1)\s_k(\c)_1 \vee [\c]_k(2)\c_k(\c)_2 \vee \cdots \vee [\c]_k(l(\c))\s_k(\c)_{l_k(\c)}.
\end{equation}
Note that if $\c$ is a path of $(T_n, E_n^*)$, then $[\c]_n = \c$ and $l_n(\c)$ coincides with the length of $\c$ denoted by $l(\c)$. In general, it follows that $l_k(\c) = l([\c]_k)$.

\definition\label{PTE.def21}
For $(w, v) \in E_n^{\ell}$, define
\[
g_{w, v} = (f_w)^{-1}{\circ}R_{w, v}{\circ}f_v.
\]
\enddefinition

The following lemma is deduced from (B4)-(a).

\lemma\label{BAS.lemma201}
For $w \in T$, define $T(w) = \cup_{m \ge 0} S^m(w)$, where $S^0(w) = \{w\}$. Let $(w, v) \in E_n^{\ell}$. Set $g = g_{v, w} \in G$. Define $R^*_{w, v}: T(w) \to T(v)$ by
\[
R^*_{w, v}(u) = vg(\s_n(u))
\]
for $u \in T(w)$. Then $R_{w, v}(K_{u}) = K_{R^*_{w, v}(u)}$ for any $u \in T(w)$.
\endlemma

\lemma\label{BAS.lemma200}
Let $\c$ be a path of $(T_n, E_n^*)$ and let $k = \{0, \ldots, n\}$. Suppose that $[\c]_k$ is a path of $(T_k, E_k^{\ell})$. Define 
\[
g_j = g_{[\c]_k(1), [\c]_k(2)}{\circ}\cdots{\circ}g_{[\c]_k(j - 1), [\c]_k(j)}
\]
for $j \in \{1, \ldots, l_*\}$, where $l_* = l_k(\c)$, and
\[
\s_k^F(\c) = \s_k(\c)_1 \vee g_2(\s_k(\c)_2) \vee \cdots  \vee g_{l_*}(\s_k(\c)_{l_*}).
\]
Then $\s_k^F(\c)$, which is called the $k$-folding of $\c$, is a path of $(T_{n - k}, E_{n - k}^*)$.
\endlemma

The next theorem gives one of basic properties of a path $\c$ of $(T_{n + m}, E_{n + m}^*)$. 

\thm\label{BAS.thm20}
Let $\c$ be a path of $(T_{n + m}, E_{n + m}^*)$. Suppose that $[\c]_n(l_n(\c)) \notin \GG_{M_J}([\c]_n(1))$, where
\[
M_J = \begin{cases}
J - 2\quad&\text{if $J$ is even,}\\
2J - 2\quad&\text{if $J$ is odd,}
\end{cases}
\]
Then there exists 
\[
A \subseteq \bigcup_{\substack{g \in G\\j = 2, \ldots, l_n(\c) - 1}} g(\s_n(\c)_j)
\]
such that $A$ is $(T_m, E_m^*)$-connected and either $\#(\partial{A}) \ge 3$ or $J$ is even and $\partial{A} = \{i, i + J/2\} \subseteq (\BbZ_J)^e$. In particular, suppose that $[\c]_n$ is a path of $(T_n, E_n^{\ell})$. Define 
\[
\widehat{\c} = w(2)\s_n(\c)_2 \vee \ldots \vee w(l_n(\c) - 1)\s_n(\c)_{l_n(\c) - 1}.
\]
Then either
\[
\#((\partial\s_n^F(\widehat{\c}))^e) \ge 3
\]
or\\
$J$ is even and $(\partial\s_n^F(\widehat{\c}))^e = \{i, i + J/2\} \subseteq (\BbZ_J)^e$.
\endthm
 
The rest of this section is devoted to a proof of Theorem~\ref{BAS.thm20}.

\lemma\label{BAS.lemma220a}
Let $\c$ be a path of $(T_{n + m}, E_{n + m}^*)$ and let 
\[
[\c]_n = (w(1), \ldots, w(l_n(\c))).
\]
Assume that $K_{w(1)} \cap K_{w(l_n(\c))} = \emptyset$.
Then there exists a path $\c'$ of $(T_{n + m}, E_{n + m}^*)$ satisfying the following properties (A),(B) and (C):\\
{\rm (A)}\,\,
$\c(1) = \c'(1)$ and $\c(l_{n + m}(\c)) = \c'(l_{n + m}(\c'))$.\\
{\rm (B)}\,\,Set $[\c']_n = (v(1), \ldots, v(l_n(\c'))$. Then $v(j - 1) \neq v(j + 1)$ for any $j \in \{2, \ldots, l_n(\c') - 1\}$,\\
{\rm (C)}\,\,
\[
\bigcup_{\substack{g \in G\\j \in \{2, \ldots, l_n(\c') - 1\}}} g(\s_n(\c')_j) \subseteq \bigcup_{\substack{g \in G\\j \in \{2, \ldots, l_n(\c) - 1\}}} g(\s_n(\c)_j).
\]
In particular, if $[\c]_n$ is a path of $(T_n, E_n^{\ell})$, then so is $[\c']_n$ and $\s_n^F(\c) = \s_n^F(\c')$.
\endlemma

\demo
Write $l_* = l_n(\c)$ for simplicity. Define
\[
D(\c)  = \{i| w(j - 1) = w(j + 1), j = 2, \ldots, l_* - 1\}.
\]
If $D(\c) = \emptyset$, then the original $\c$ itself has the desired property. Otherwise, let $j_* = \min D(\c)$. Set $w = w(j_* - 1) = w(j_* + 1)$ and $v = w(j_*)$. Then there are two cases:\\
Case I: $K_w \cap K_v = Q_w \cap Q_v = \{p\}$ for some $p$.\\
In this case, if $\s_n(\c)_{j_* - 1} = (u(1), \ldots, u(k))$ and $\s_n(\c)_{j_* + 1} = (u'(1), \ldots, u(k'))$, then $u(k) = u'(1)$. Replacing the part 
\[
w\s_n(\c)_{j_* - 1} \vee v\s_n(\c)_{j_*} \vee w\s_n(\c)_{j_* + 1)}
\]
of $\c$ with
\[
(wu(1), \ldots, wu(k), wu'(2), \ldots, wu'(k')),
\]
we obtain a path $\xi(\c)$ of $(T_{n + m}, E_{n + m}^*)$.\\
Case II: $Q_w \cap Q_v = b_i(w) = b_{i'}(v)$ for some $i, i' \in \BbZ_J$. \\
In this case, 
\[
R_{w,v}\big(K(v\s_n(\c)_{j_*})\big) = K(wg_{w, v}(\s_n(\c)_{j_*})).
\]
So, replacing the part
\[
w\s_n(\c)_{j_* - 1} \vee v\s_n(\c)_{j_*} \vee w\s_n(\c)_{j_* + 1)}
\]
of $\c$ with
\[
w\s_m(c)_{j_* - 1} \vee wg(\s_n(\c)_{j_*}) \vee w\s_n(\c)_{j_* + 1}
\]
we obtain a path $\xi(\c)$ of $(T_{n + m}, E_{n + m}^*)$. \par
In both cases, letting $\xi(\c) = (u(1), \ldots, u(k_*)$, where $k_* = l_n(\xi(\c))$, we see that $k_* = l_n(\xi(\c)) = l_n(\c) - 2$, $\#(D(\xi(\c))) \le \#(D(\c)) - 1$ and 
\[
\bigcup_{\substack{g \in G\\j \in \{2, \ldots, k_* - 1\}}} g(\s_n(\xi(\c))_j) \subseteq \bigcup_{\substack{g \in G\\j \in \{2, \ldots, j_* - 1\}}} g(\s_n(\c)_j).
\]
Iterating $\xi$, we eventually see that $D(\xi^N(\c)) = \emptyset$ for some $N \in \BbN$. Letting $\c' = \xi^N(\c)$, we have the desired properties. Now suppose that $[\c]_n$ is a path of $(T_n, E_n^{\ell})$. Then Case I never happens and $[\xi(\c)]_n$ is a path of $(T_n, E_n^{\ell})$. Moreover, since $g_{w(j_* - 1), w(j_*)}{\circ}g_{w(j_*), w(j_* + 1)}$ is an identity, we see that $\s_n^F(\c) = \s_n^F(\xi(\c))$. This immediately implies the desired properties of $\c'$.
\enddemo

\lemma\label{BAS.lemma231}
Let $\c$ be a path of $(T_{n + m}, E_{n + m}^*)$ with $K_{[\c]_n(1)} \cap K_{[\c]_n(l_*)} = \emptyset$, where $l_* = l_n(\c)$. Assume that $[\c]_n(i - 1) \neq [\c]_n(i + 1)$ for any $i = 2, \ldots, l_* - 1$. Moreover, assume that if
\begin{equation}\label{BAS.eq200}
A \subseteq \bigcup_{\substack{g \in G\\j = 2, \ldots, l_* -1}} g(\s_n(\c)_j)
\end{equation}
and $A$ is $(T_m, E_m^*)$-connected, then $\#(A) \le 2$. Then $[\c]_n$ is a path of $(T_n, E_n^{\ell})$.
\endlemma

\demo
Set $Y_j = \partial{\s_n(\c)_j}$ for $j = 1, \ldots, l_*$. Since $\s_n(\c)_j$ is $(T_n, E_n^*)$-connected, we see that $\#(Y_j) \le 2$. Furthermore, for $j = 2, \ldots, l_* -1$, since $w(j - 1) \neq w(j + 1)$, we see that $\#(Y_j)  = 2$. \\
{\bf Claim} Let $j \in \{2, \ldots, l_* - 1\}$. Assume that $Y_j = \{i, i + 1\}$ and $p_i \in K$ for some $i \in \BbZ_J$. Then there exist $i', i'' \in \BbZ_J$ such that $p_{i'}, p_{i''} \in K$ and $f_{w(j)}(p_i) = f_{w(j + 1)}(p_{i'}) = f_{w(j - 1)}(p_{i''})$. Moreover, 
\[
\begin{cases}
Y_{j + 1} &= \{i', i' + 1\}\quad\text{if $j \le l_* - 2$},\\
Y_{j - 1} &= \{i''. i'' + 1\}\quad\text{if $3 \le j$}.
\end{cases}
\]
Proof of Claim:\,\,First we give a proof for the statements on $i'$. Since $\#(Y_j) = 2$, one of the following situations (a) or (b) occurs.\\
(a)\,\,$K_{w(j)} \cap K_{w(j + 1)} = \{f_{w(j)}(p_i)\}$\\
(b)\,\,$Q_{w(j)} \cap Q_{w(j + 1)} = b_{i_1}(w(j)) = b_{i_2}(w(j + 1))$, where $i_1 \in \{i, i + 1\}$.\\
In the case (a), the statements of the claim are obvious. In the case (b), let $(w, v) = (w(j), w(j + 1))$ for simplicity. Then since $R_{w, v}(K_w) = K_v$, we see that $R_{w, v}(f_{w}(p_i)) = f_{w}(p_i) \in K_v$. Therefore there exists $i' \in \BbZ_J$ such that $p_{i'} \in K$ and $f_v(p_{i'}) = f_w(p_i)$. Suppose that $j \le l_* - 1$. Let $g = (f_w)^{-1}{\circ}R_{w, v}{\circ}f_v$ and define $\tau = \s_n(\c)_j \vee g(\s_n(\c))_{j + 1}$. Then $\tau$ is a path in $(T_m, E_m^*)$ contained in the right-hand side of \eqref{BAS.eq200}. Moreover, it follows that $\partial{\tau} = Y_j \cup g(Y_{j + 1})$. Since $\#(\partial{\tau}) \le 2$, we see that $\partial{\tau} = Y_j$ and hence $g(Y_{j + 1}) = Y_j$. Since $R_{w, v}(f_w(p_i)) = f_v(p_{i'})$, it follows that $Y_{j + 1} = \{i', i' + 1\}$. Thus we have show the desired statements on $i'$. A proof of those on $i''$ is entirely similar. \qed\\
Now suppose that $(w(1), \ldots, w(l_*))$ is not a path of $(T_n, E_n^{\ell})$. Then there exists $k \in \{1, \ldots, l_* - 1\}$ and $i(k), i(k + 1) \in \BbZ_J$ such that $Q_{w(k)} \cap Q_{w(k + 1)} = f_{w(k)}(p_{i(k)}) = f_{w(k + 1)}(p_{i(k + 1)})$. Then for $j \in \{k, k + 1\} \in \{1, \ldots, l_*\}$, we see that $Y_j = \{i(j), i(j) + 1\}$ and $p_{i(k)} \in K$. Then inductive use of the above claim shows that there exists $i_j \in \BbZ_J$ for any $j \in \{1, \ldots, l_*\}$ such that $p_{i_j} \in K$ and $f_{w(j)}(p_{i_j})$ is independent of $j$ and belongs of every $K_{w(j)}$. Consequently it follows that $K_{w(1)} \cap K_{w(l_* )} \neq \emptyset$. This contradiction shows the lemma.
\enddemo

\demo[Proof of Theorem~\ref{BAS.thm20}]
First assume that $[\c]_n$ is a path of $(T_n, E_n^{\ell})$. Using Lemma~\ref{BAS.lemma220a}, we replace $\widehat{\c}$ with $(\widehat{\c})'$ satisfying $\s_n^F(\widehat{\c}) = \s_n^F((\widehat{\c})')$. Denoting the modified path again by $\c$, we may assume that 
\begin{equation}\label{BAS.eq250}
[\c]_n(j - 1) \neq [\c]_n(j + 1)
\end{equation}
for any $j = 3, \ldots, l_n(\c) - 2$ without loss of generality. Write $l_n(\c) = l_*$ and let $[\c]_n = (w(1), \ldots, w(l_*))$. Then there exist subsets $\{i_j\}_{1 \le j \le l_* -1}$ and $\{k_j\}_{2 \le j \le l_*}$ of $(\BbZ_J)^e$ such that
\[
Q_{w(i)} \cap Q_{w(i + 1)} = b_{i_j}(w(j)) = b_{k_{j + 1}}(w(j + 1))
\]
for any $j = 1, \ldots, l_* - 1$. Letting $g_j = (f_{w(j)})^{-1}{\circ}R_{w(j), w(j + 1)}{\circ}f_{w(j + 1)}$, we see that $g_j(k_{j + 1}) = i_j$. Moreover, define
\[
h_j = 
\begin{cases}
I\,\,&\text{if $j = 1$,}\\
g_2{\circ}g_3{\circ}\cdots{\circ}g_j&\text{if $j = 2, \ldots, l_* - 1$}
\end{cases}
\]
and $\xi_j = h_j(k_{j + 1})$ for $j = 2, \ldots, l_* - 1$. Then it follows that
\[
\{\xi_1, \xi_2, \ldots, \xi_{l_* - 1}\} \subseteq (\partial\s_n(\widehat{\c}))^e.
\]
Assume that 
\begin{equation}\label{BAS.eq260}
\#((\partial{\s_n^F(\widehat{\c})})^e) \subseteq \{\a_1, \a_2\}
\end{equation}
for some $\a_1, \a_2 \in (\BbZ_J)^e$ with $\a_1 \neq \a_2$. By \eqref{BAS.eq250}, it follows that $k_j \neq i_j$ for $j \in \{3, \ldots, l_* - 2\}$. Since $\xi_{j - 1} = h_{j - 1}(k_j)$ and $\xi_j = h_j(k_{j + 1}) = h_{j - 1}(i_j)$, we see that $\xi_{j - 1} \neq \xi_j$. Due to \eqref{BAS.eq260}, it follows that $\xi_j = \xi_{j + 2}$ for any $j = 2, \ldots, l_* - 4$ and $\xi_j \neq \xi_{j + 1}$ for $j = 2, \ldots, l_* - 3$. This implies that 
\[
R_{w(j), w(j + 1)}(b_{w(j)}(k_j)) = b_{w(j + 2)}(k_{j + 2}) = b_{w(j + 1)}(i_{j + 1})
\]
 for any $j = 2, \ldots, l_* - 3$. Let $\ell_j$ be the straight line containing $b_{w(j)}(i_j)$. Then $\ell_j$ and $\ell_{j + 2}$ are symmetric with respect to $\ell_{j + 1}$ for any $j = 2, \ldots, l_* - 4$. Then unless $J$ is even and $\d(\a_1,\a_2) = J/2$, $\ell_j$'s for $j = 2, \ldots, l_* - 2$ share a common point $p_*$ and the angle between $\ell_j$ and $\ell_{j + 1}$ dose not depends on $j \in \{2, \ldots, l_* - 3\}$. More precisely, the angle $\theta$ is given by
\[
\theta = \Big(1- \frac{2k}J\Big)\pi
\]
where $k = 1, \ldots, (J - 2)/2$ if $J$ is even and $k = 1, \ldots, (J - 1)/2$ if $J$ is odd. So, $Q_{w(2)}, Q_{w(3)}, \ldots, Q_{w(l_* - 3)}$, and $Q_{w(l_* - 1)}$ form a part of a circle around $p_*$. Note that $\ell_1$ and $\ell_2$ also pass through $p_*$ and $Q_{w(1)}$ and $Q_{w(l_*)}$ are also parts of the same circle. Since $Q_{w(1)} \cap Q_{w(l_*)} = \emptyset$, it follows that
\[
(M + 2) \theta \le (k + 2)\theta < 2\pi.
\]
This shows that $M < M_J$. Therefore, if $M \ge M_J$, the remaining possibility is that $|\d(\a_1, \a_2) = J/2$ for any $i = 1, \ldots, k$. This concludes a proof in the case of $[\c]_n$ is a path of $(T_n, E_n^{\ell})$. \par
Now we proceed to the general case. Using Lemma~\ref{BAS.lemma220a} and replacing $\c$ with $\c'$, we may assume that $[\c]_n(j - 1) \neq [\c]_n(j + 1)$ for any $j = 2, \ldots, l_n(\c) - 1$. Suppose that if
\[
A \subseteq \bigcup_{\substack{g \in G\\j \in \{2, \ldots, l_n(\c) - 1\}}} g(\s_n(\c)_j),
\]
and $A$ is $(T_m, E_m^*)$-connected, then $\#(A) \le 2$. By Lemma~\ref{BAS.lemma231}, we see that $[\c]_n$ is a path of $(T_n, E_n^{\ell})$. Applying the above argument of this proof, we obtain the desired statements.
\enddemo

\setcounter{equation}{0}
\section{Proofs of Theorems~\ref{MTH.thm10}, \ref{MTH.thm20} and \ref{MTH.thm30}}\label{PTM}

\definition\label{BAS.def30}
(1)\,\,Let $n \ge 1$. Two subsets $A$ and $B$ of $T_n$ are said to be alternated if there exist $x_1, x_2 \in K(A) \cap \partial{Q_*}$ and $y_1, y_2 \in K(B) \cap \partial{Q_*}$ such that $\{x_1, x_2\} \cap \{y_1, y_2\} \neq \emptyset$ or the connected component of $\sd{\partial{Q_*}}{\{x_1, x_2\}}$ including $y_1$ is different from that including $y_2$.\\
(2)\,\,Let $n \ge 1$ and let $A \subseteq T_n$. For a subgroup $G \subseteq D_J$, define $G(A)$ by
\[
G(A) = \bigcup_{g \in G} g(A).
\]
\enddefinition

\lemma\label{BAS.lemma350}
For any $A \subseteq T$ and any $g \in D_J$, 
\begin{equation}\label{BAS.eq710}
g(\partial{A}) = \partial{g(A)}.
\end{equation}
\endlemma

\demo
Note that $g(K(A)) = K(g(A))$. Let $s \in \partial{g(A)}$. Then $g(K(A)) \cap b_s \neq \emptyset$. Hence $K(A) \cap b_{g^{-1}(s)} \neq \emptyset$ and so $g^{-1}(s) \in \partial{A}$. This implies $\partial{g(A)} \subseteq g(\partial{A})$. The other direction is entirely the same.
\enddemo

\lemma\label{BAS.lemma400}
{\rm (1)}\,\,Let $n \ge 1$. If $A$ and $B$ are subsets of $T_n$, $(T_n, E_n^*)$-connected and alternated, then $A \cup B$ is $(T_n, E_n^*)$-connected. In particular, if $J \ge 5$, then $A \cap B \neq \emptyset$.\\
{\rm (2)}\,\,Suppose that $J = q_1q_2$ for some $q_1, q_2 \in \BbN$. Let $n \ge 1$ and let $A \subseteq T_n$. If $A$ is $(T_n, E_n^*)$-connected and there exist $k, l \in \partial{A}$ such that $\d(k, l) \ge q_2 + 1$. Let $g = \Theta_{2\pi/q_1}$. Then $g^j(A)$ and $g^{j + 1}(A)$ are alternated for any $j \ge 0$ and $Rot_{q_1}(A)$ is $(T_n, E_n^*)$-connected. In particular, if $J \ge 5$ and $A$ is $(T_n, E_n^{\ell})$-connected, then $g^j(A) \cap g^{j + 1}(A) \neq \emptyset$ for any $j \ge 0$ and $Rot_{q_1}(A)$ is $(T_n, E_n^{\ell})$-connected.
\endlemma

\remark
Note that if $A \subseteq T_m$ for some $m \ge 1$, then $K(A)$ is connected if and only if $A$ is $(T_m, E_m^*)$-connected.
\endremark

\demo
(1)\,\,
Choose $x_1, x_2, y_1$ and $y_2$ in the same manner as Definition~\ref{BAS.def30}. Let $(w(1), \ldots, w(k))$ be a path of $(T_n, E_n^*)$ contained in $A$ and satisfying $x_1 \in K_{w(1)}$ and $x_2 \in K_{w(2)}$. Define $\a_0 = x_1$, $\a_k = x_2$ and
\[
\a_i = \begin{cases}
\text{the midpoint of $Q_{w(i)} \cap Q_{w(i + 1)}$}\quad&\text{if $(w(i), w(i + 1)) \in E_n^{\ell}$,}\\
\text{the single point of $Q_{w(i)} \cap Q_{w(i + 1)}$}\quad&\text{otherwise}.
\end{cases}
\]
for $i = 1, \ldots, k - 1$. Define $\c_1$ as the broken line connecting the points $\a_0$, $c_{w(1)}$, $\a_1$, $c_{w(2)}$, $\a_2$, $\ldots$, $c_{w(k - 1)}$, $\a_{k - 1}$, $c_{w(k)}$ and $\a_k$. Next let $(v(1), \ldots, v(l))$ be a path of $(T_n, E_n^*)$ contained in $B$ and satisfying $y_1 \in K_{v(1)}$ and $y_2 \in K_{v(l)}$. Construct a broken line $\c_2$ in the same as $\c_1$ replacing $x_1, x_2$ and $(w(1), \ldots, w(k))$ with $y_1, y_2$ and $(v(1), \ldots, v(l))$ respectively. Let $\b_0, \b_1, \ldots, \b_{l - 1}$ and $\b_l$ be the counterparts of $\a_0, \a_1, \ldots, \a_{k - 1}$ and $\a_k$ respectively. Since $A$ and $B$ are alternated, we see that $\c_1$ and $\c_2$ has an intersection. Let $x_* \in \c_1 \cap \c_2$. \\
{\bf Claim}\,\,If $A \cap B = \emptyset$, then either $x_* \in \{x_1, x_2\} \subseteq K$ or $\{x_*\} = K_{w(i)} \cap K_{w(i + 1)} = Q_{w(i)} \cap Q_{w(i + 1)}$ for some $i = 1, \ldots, k - 1$.\\
{Proof of Claim}\,\,
If $x_*$ belongs to an interior of $Q_{w(i)}$ for some $i = 1, \ldots, k$, then $w(i) \in A \cap B$. Therefore, if $A \cap B = \emptyset$, then $x_* \in \{\a_0, \a_1, \ldots, \a_k\} \cap \{\b_0, \b_1, \ldots, \b_l\}$. Suppose that $x_* = \a_i = \b_j$ for some $i \in \{0, 1, \ldots, k\}$ and $j \in \{0, 1, \ldots, l\}$. Suppose that $i \in \{1, \ldots, k - 1\}$ and that $(w(i), w(i + 1)) \in E_n^{\ell}$. Since $x_i$ belongs to the interior of $Q_{w(i)} \cap Q_{w(i + 1)}$, we see that $\{w(i), w(i + 1)\} \cap B \neq \emptyset$. So this contradiction shows that $\{x_*\} = Q_{w(i)} \cap Q_{w(i + 1)} = K_{w(i)} \cap K_{w(i + 1)}$. In the case $i \in \{0, k\}$, we immediately see that $x_* \in \{x_1, x_2\}$. Thus we have shown the claim.\par
Now, if $A \cap B \neq \emptyset$, then $A \cup B$ is $(T_n, E_n^*)$-connected. In the case when $A \cap B = \emptyset$, the above claim implies $x_* \in K$. Since $x_* \in K_{w(i)}$ and $x_* \in K_{v(j)}$ for some $i, j$, we see that $(w(i), v(j)) \in E_n^*$. Thus $A \cup B$ is $(T_n, E_n^*)$-connected. \par
Finally assume that $J \ge 5$ and that $A \cap B = \emptyset$. Let $\theta_*$ be the angle between $b_0$ and $b_1$. Then $\theta_* = (1 - \frac 2J)\pi$ and $2\theta_* > \pi$ if $J \ge 5$. Hence
\begin{equation}\label{BAS.eq500}
\#(\{w| w \in T_n, x_* \in K_w\}) \le \begin{cases}
1\quad&\text{if $x_* \in \partial{Q_*}$,}\\
3\quad&\text{if $x_*$ belongs to the interior of $Q_*$.}
\end{cases}
\end{equation}
By the above argument, it follows that $x_* \in K_{w(i)} \cap K_{v(j)}$ and $w(i) \neq v(j)$. By \eqref{BAS.eq500}, the point $x_*$ belongs to the interior of $Q_*$. Then due to the above claim, we see that there exist $i \in \{1, \ldots, k - 1\}$ and $j \in \{1, \ldots, l - 1\}$ such that $x_* \in K_{w(i)} \cap K_{w(i + 1)} \cap K_{v(j)} \cap K_{v(j + 1)}$. By \eqref{BAS.eq500}, we have contradiction. Thus we see that $A \cap B \neq \emptyset$.\\
(2)\,\,For simplicity of the arguments, we let $k = 0$. Since $\d(k, l) \ge q_2 + 1$, it follows that $q_2 < l < J - q_2$. Choose $x_1$ and $x_2$ such that $x_1 \in K(A) \cap b_0$ and $x_2 \in K(A) \cap b_l$. Let $y_i = g(x_i)$ for $i = 1, 2$. Then $y_1 \in b_{q_2} \cap K(g(A))$ and $y_1 \in \b_{l + q_2} \cap K(g(A))$. Since $0 < q_2 < l < l + q_2 < J$, we see that $A$ and $g(A)$ are alternated. Thus using (1), we see that $A \cup g(A)$ is $(T_n, E_n^*)$-connected. Similar arguments show that $g^j(A) \cup g^{j + 1}(A)$ is $(T_n, E_n^*)$-connected for any $j \ge 0$. Since $Rot_{q_1}(A) = \cup_{j = 0, 1, \ldots, q_2 -1} g^j(A)$, we see that $Rot_{q_1}(A)$ is $(T_n, E_n^*)$-connected. If $J \ge 5$, the conclusion of (1) implies that $g^j(A) \cap g^{j + 1}(A) \neq \emptyset$ for any $j \ge 0$. Hence $Rot_{q_1}(A)$ is $(T_n, E_n^{\ell})$-connected if $A$ is $(T_n, E_n^{\ell})$-connected.
\enddemo

\demo[Proof of Theorem~\ref{MTH.thm10}]
Suppose $J = 3$. Choose $M \ge M_3 = 4$. Let $w \in T$ and let $\c \in \C_{M, m}(w)$. Then there exist $u \in S^m(w)$ and $v \in S^m(\GG_M(w))^c)$ such that $(u, \c(1)), (\c(l(\c)), v) \in (T_{n + m}, E_{n + m}^*)$, where $n = |w|$.  Applying Theorem~\ref{BAS.thm20} to $\widetilde{\c} = (u, \c(1), \ldots, \c(l(\c)), v)$, we obtain $A \subseteq \H_*(\c)$ which is $(T_m, E_m^*)$-connected and $\#(\partial{A}) \ge 3$. Since $\#(\BbZ_3) = 3$, it follows that $\partial{A} = \BbZ_3$ and hence $(\BbZ_3)^e \subseteq \partial{A}$. So, Theorem~\ref{BAS.thm10} immediately shows the desired statement.\par
\enddemo

\demo[Proof of Theorem~\ref{MTH.thm20}]
Let $M \ge M_J$, $w \in T$ and $m \ge 1$. Let $\c \in \C_{M, m}(w)$ and define $\widetilde{\c}$ in the same way as in the proof of Theorem~\ref{MTH.thm10}. By Theorem~\ref{BAS.thm20}, there exist $A \subseteq \H_*(\c)$, and $i_1, i_2 \in \partial{A}$ such that $A$ is $(T_m, E_m^*)$-connected and $b_{i_1} \cap b_{i_2} = \emptyset$. Note that $\d(i_1,i_2) \ge 2$. Now we have two cases.\\
{\bf Case 1}: $G = D_J$ or $Rot_J$.\\
Due to Theorem~\ref{MTH.thm10}, it is enough to consider the case when $J \ge 4$. Applying Lemma~\ref{BAS.lemma400}-(2) with $q_1 = J$ and $q_2 = 1$, we see that $Rot_J(A)$ is $(T_m, E_m^*)$-connected. Since $Rot_J(A) \subseteq \H_*(\c)$ and $\partial(Rot_J(A)) = \cup_{g \in Rot_J}\partial{g(A)} = \BbZ_J \subseteq (\BbZ_J)^e$, Theorem~\ref{BAS.thm10} yields the statement of the theorem.\\
{\bf Case 2}: $J$ is even and $G = D_{J/2}^V$.\\
In the case $J = 4$, our assertion coincides with \cite[Theorem~4.13]{Ki22} with $L = 2$. So, we assume $J \ge 6$. First we construct $B \subseteq \H_*(\c)$ such that $B$ is $(T_m, E_m^*)$-connected and there exist $j_1, j_2 \in \partial{B}$ such that $\d(j_1, j_2) \ge 3$ and $j_1 - j_2$ is odd. Suppose that $i_1 - i_2$ is odd. Then $\d(i_1, i_2) \ge 3$. Hence, it is enough to set $B = A$. Suppose that $i_1 - i_2$ is even. Then, without loss of generality, we have $i_2 = i_1 + 2k$ for some $k \in \BbN$ with $k \le J/4$. Let $R = R_{\theta_{i_1}}$. Note that $R$ is a reflection in the line $p_{i_1}p_{i_1 + J/2}$. Since $R(A)$ contains $R(i_1) = i_1 + 1$ and $R(i_2) = i_1 - 2k + 1$, $A$ and $R(A)$ are alternated. By Lemma~\ref{BAS.lemma400}-(1), letting $B = A \cup R(A)$, we see that $B$ is $(T_m, E_m^*)$-connected and $B \subseteq \H_*(\c)$. Since $\partial{B}$ contains four distinct points $i_1, i_2, R(i_1)$ and $R(i_2)\}$ and $J \ge 6$, there exist $j_1, j_2 \in \partial{B}$ such that $\d(j_1, j_2) \ge 3$. Thus, we have obtained $B$ having all the required properties. Letting $q_1 = J/2$ and $q_2 = 2$ and applying Lemma~\ref{BAS.lemma400}-(2), we see that $Rot_{q_1}(B)$ is $(T_m, E_m^*)$-connected. Since $Rot_{q_1} \subseteq D_{J/2}^V$, it follows that $Rot_{q_1}(B) \subseteq \H_*(\c)$. Moreover, since $\BbZ_J = Rot_{q_1}(\{i_1, i_1 + 1\}) \subseteq \partial{Rot_{q_1}(B)}$, Theorem~\ref{BAS.thm10} yields the statement of the theorem.
\enddemo

Finally, we proceed to a proof of Theorem~\ref{MTH.thm30}.

\lemma\label{MTH.lemma10}
Let $G$ be a subgroup of $D_J$. Define $q_1 = \#(G \cap SO(2))$. Set $q_2 = J/q_1$ and $\theta = \frac{2\pi}{q_1}$. Then $q_2 \in \BbN$. Furthermore, let $X \subseteq \BbZ_J$ is $G$-transitive. Then either (A) or (B) holds:\\
{\rm (A)}\,\,$G = Rot_{q_1}$ and $X = Rot_{q_1}(i)$ for some $i \in \BbZ_J$.\\
{\rm (B)}\,\,$G = Rot_{q_1} \cup R_{\tau}(Rot_{q_1})$ for some $\tau \in \BbR$. There exists $i \in \BbZ_J$ such that $X = Rot_{q_1}(i) \cup Rot_{q_1}(\{R_{\tau}(i)\})$.
\endlemma

\demo
Let $G_+ = G \cap SO(2)$ and let $G_- = \sd{G}{G_+}$. Then $G_- = \{g| g \in G, \det{g} = -1\}$ and $G_+ = Rot_{q_1}$. Since $Rot_{q_1}$ is a subgroup of $Rot_J$, the number $q_1$ divides into $J$. First, we consider the case when $G_+ = G$. In this case, if $X$ is $G$-transitive, then $X$ is an orbit of any/some $i \in X$ by $G = Rot_{q_1}$. Hence $X = Rot_{q_1}(i)$ for some $i \in \BbZ_J$. Second, we assume that $G_- \neq \emptyset$. In this case, choose $h \in G_-$ and define $h^*: G_+ \to G_-$ by $h^*(g) = h{\circ}g$. Then $h^*$ is bijective. Moreover, there exists $\tau \in \BbR$ such that $R_{\tau}$. Thus we see that $G = Rot_{q_1} \cup R_{\tau}(Rot_{q_1})$. As in the other case, if $X$ is $G$-transitive, then $X$ is an orbit of any/some $i \in X$ by $G$. Hence $X = Rot_{q_1}(i) \cup Rot_{q_1}(R_{\tau}(i))$.
\enddemo

\lemma\label{MTH.lemma20}
Let $G$ be a subgroup of $D_J$. Define $q_1 = \#(G \cap SO(2))$ and $q_2 = J/q_1$. Assume that $X \subseteq \BbZ_J$ and $X$ is $G$-transitive. Suppose that $A_0 \subseteq T_m$, $A_0$ is $(T_m, E_m^*)$-connected and there exist $i_1, i_2 \in X \cap (\partial{A_0})^e$ such that $\d(i_1, i_2) \ge q_2 + 1$. Then there exists $A \subseteq G(A_0)$ such that $A$ is $(T_m, E_m^*)$-connected and $X \subseteq (\partial{A})^e$.
\endlemma

\demo
By Lemma~\ref{BAS.lemma400}-(2), we see that $Rot_{q_1}(A_0)$ is $(T_m, E_m^*)$-connected. In the case (A) of Lemma~\ref{MTH.lemma10}, letting $A = Rot_{q_1}(i_1)$, we see that $A$ satisfies the desired properties. In the case (B) of Lemma~\ref{MTH.lemma10},  $Rot_{q_1}(A_0)$ and $Rot_{q_1}(R_{\tau}(A_0))$ are alternated. Define $A = Rot_{q_1}(A_0) \cup Rot_{q_1}(R_{\tau}(A_0))$. Then by Lemma~\ref{BAS.lemma400}-(1), $A$ is $(T_m, E_m^*)$-connected. Moreover, since $X = Rot_{q_1}(i_1) \cup Rot_{q_1}(R_{\tau}(i_1))$, we see that $X \subseteq (\partial{A})^e$. 
\enddemo

\demo[Proof of Theorem~\ref{MTH.thm30}]
Let $w \in T$ and let $\c \in \C^{\ell}_{M_J, m}(w)$. Set $B = \s_n^F(\c)$. Then Theorem~\ref{BAS.thm20} implies that $\#(({\partial}B)^e) \ge 3$ or $J$ is even and $({\partial}B)^e = \{k, k + J/2\}$ for some $k \in \BbZ_J$. Keeping this in mind, we are going to show the following claim.\\
{\bf Claim}\,\,There exists $A \subseteq G(B)$ such that $A$ is $(T_m, E_m^*)$-connected and $(\BbZ_J)^e \subseteq \partial{A}$.\par
By Theorem~\ref{BAS.thm10}, the above claim implies the theorem. The rest of the proof is devoted to showing the claim. There are two cases.\\
{\bf Case 1}: $(\BbZ_J)^e = Rot_{q_1}(i)$.\\
First assume that $\#(({\partial}B)^e) \ge 3$. Suppose that there exists $i_1, i_2 \in ({\partial}(B))^e$ with $\d(i_1, i_2) \ge q_2 + 1$. Letting $X = (\BbZ_J)^e$ and applying Lemma~\ref{MTH.lemma20} with $A_0 = B$, we have the claim. Otherwise it follows that $(\partial{B})^e = \{i_1, i_2, i_3\}$, $\d(i_1, i_2) = \d(i_2, i_3) = d(i_3, i_1) = q_2$.  Then $q_1 = 3$ and $(\BbZ_J)^e = \{i_1,, i_2, i_3\} = (\partial{B})^e$. Hence, letting $A = B$, we have the claim. \par
Next assume that $\#(({\partial}B)^e) = 2$. Then $({\partial}B)^e = \{i, i + J/2\}$. If $q_1 > 2$, then $d(i, i + J/2) \ge q_2 + 1$. Set $i_1 = i$ and $i_2 = i + J/2$. By the same arguments as above, we have the claim. If $q_1 = 2$, then $(\BbZ_J)^e = \{i, i + J/2\} = (\partial{B})^e$, so that it is enough to let $A = B$.\\
{\bf Case 2}: $(\BbZ_J)^e = Rot_{q_1}(i) \cup Rot_{q_1}(R_{\tau}(i)\})$ and $R_{\tau}(i) \notin Rot_{q_1}(i)$.\\
First, assume that $q_1 = 2$. Then $(\BbZ_J)^e = \{i, i + J/2\} \cup \{i', i' + J/2\}$, where $i' = R_{\tau}(i)$. By the property of $(\partial{B})^e$ given at the beginning of the proof, it follows that $\{k, k + J/2\} \subseteq (\partial{B})^e$ for some $k \in \{i, i'\}$. Therefore $B$ and $R_{\tau}(B)$ is alternated. Define $A = B \cup R_{\tau}(B)$. Using Lemma~\ref{BAS.lemma400}-(1), we verify the claim.\par
Next assume $q_1 > 2$. If there exists $i_1, i_2 \in (\partial{B})^2$, then Lemma~\ref{MTH.lemma20} suffices. Hence assume that $\d(i_1, i_2) \le q_2$ for any $i_1, i_2 \in (\partial{B})^e$. Since $J/2 > q_2$, this assumption implies that $\#(\partial{B})^e \ge 3$. Let $X_0 = Rot_{q_1}(i)$ and let $X_1 = Rot_{q_1}(R_{\tau}(i)\})$. Then $(\partial{B})^e \cap X_j$ contains more than two point for some $j \in \{0, 1\}$. Let $j = 0$ for simplicity. Thus there exists $k \in (\partial{B})^e$ such that $\{k, k + q_2\} \subseteq (\partial{B})^e \cap X_0$. Now there exists $l \in X_1$ such that $l = k + k'$ for some $k' \in \{1, 2, \ldots, p_2 - 1\}$. Note that both $l - q_2$ and $l + q_2$ belong to $X_1$. Since $q_1 > 2$, it follows that $l + q_2 \neq l - q_2$ and hence $\d(l + q_2, l - q_2) \ge q_2$. This implies $\d(k, l + q_2) \ge q_2 + 1$. Note that $\Theta_{\frac{2\pi}J}(\{k, k + q_2\}) = \{l, l + q_2\}$. Since $B$ and $\Theta_{\frac{2\pi}J}(B)$ are alternated, if $A_0 = B \cup \Theta_{\frac{2\pi}J}(B)$, then $k, l + q_2 \in (\partial{A_0})^e$ and $A_0$ is $(T_m, E_m^*)$-connected by Lemma~\ref{BAS.lemma400}-(1). Now Lemma~\ref{MTH.lemma20} suffices.
\enddemo

\setcounter{equation}{0}
\section{Conductive homogeneity of $(J, G)$-s.s.\,system II: less symmetric with even $J$}\label{ENJ}
As in the previous sections, $(S, \{f_i\}_{s \in S})$ is a $(J, G)$ self-similar system with $J \ge 3$, $(T, \A)$ is the associated partition, $K$ is the associated self-similar set and $\mu_*$ is the self-similar measure defined in Theorem~\ref{RPB.thm10} throughout this section. Furthermore, we will assume that $J$ is even in this section. Then as we have observed in Lemma~\ref{MTH.lemma100}, subgroups of $D_V$ with $J$ elements are $Rot_J$, $D_{J/2}$ and $D_{J/2}^V$. In the case where $G = Rot_J$ or $D_{J/2}^V$, Theorem~\ref{MTH.thm20} yields the conductive homogeneity of the self-similar set $K$. Consequently, a natural question is what happens when $G = D_{J/2}$ where $\BbZ_J$ is not $G$-transitive. In fact, we will study the case when $G = Rot_{J/2}$ as well in this section. In either case, the orbits of $G$-action on $\BbZ_J$ are 
\begin{equation}\label{ENJ.eq10}
\BbZ_J^0 = \{i| i \in \BbZ_J, i \equiv 0 \mod 2\}\,\,\text{and}\,\, \BbZ_J^1 = \{i| i\in \BbZ_J, i \equiv 1 \mod 2\}.
\end{equation}
Moreover $G$-invariant subset of $\BbZ_J$ are $\BbZ_J$, $\BbZ_J^0$ and $\BbZ_J^1$ and hence $(\BbZ_J)^e$ must be one of them. In particular if $(\BbZ_J)^e = \BbZ_J^k$ for some $k \in \{0, 1\}$, then $(\BbZ_J)^e$ is $G$-transitive and Theorem~\ref{MTH.thm30} ensures the conductive homogeneity. In conclusion, if $G$ is $D_{J/2}$ or $Rot_{J/2}$, the remaining case is when $(\BbZ_J)^e = \BbZ_J = \BbZ_J^0 \cup \BbZ_J^1$. To determine the conductive homogeneity in such a case, we have to deal with not only the group $G$  and its action of $\BbZ_J$ but also how the contractions $\{f_s\}_{s \in S}$ map the boundaries $\{b_s\}_{s \in S}$. More precisely, our sufficient condition of the conductive homogeneity is given in terms of the notions defined below. Note that we do not need to assume that $J$ is even in the following definitions.

\definition\label{MIS.def200}
Let $Y \subseteq \BbZ_J$ and let $A \subseteq T_n$. Define
\begin{multline*}
E_n^{\ell}(Y) = \{(u, v) | (u, v) \in E_n^{\ell}, \text{there exist $i_1, i_2 \in Y$}\\
\text{such that $Q_w \cap Q_v = b_{i_1}(u) = b_{i_2}(v)$.}\}.
\end{multline*}
and
\[
\partial_YA = \{i| i \in \BbZ_J, \text{there exist $w \in A$ and $i' \in Y$ such that $b_{i'}(w) \subseteq b_i$}\}.
\]
In particular, we write $\partial_{\ell}A = \partial_{\BbZ_J}A$. Moreover, for $X \subseteq \BbZ_J$, define
\[
\Con_T(n, X) = \big\{A\big| A \subseteq T_n, \text{$A$ is a connected component of $(T_n, E_n^{\ell}(X^c))$}\big\}
\]
and
\[
F_{\partial}(n, X) = \{(\partial_{X^c}A)^c| A \in \Con_T(n, X)\}.
\]
In particular, in the case $n = 1$, we omit $n$ in $\Con_T(n, X)$ and $F_{\partial}(n, X)$ and write $\Con_T(X)$ and $F_{\partial}(X)$ respectively.
\enddefinition

\remark
Assume that $J$ is even.  If $b_{i'}(w) \cap b_i \neq \emptyset$ for $i, i' \in \BbZ_J$ and $w \in T$, then $b_{i'}(w) \subseteq b_i$. Hence $\partial{A} = \partial_{\ell}A$ for any $n \ge 1$ and $A \subseteq T_n$.
\endremark

The next theorem is the main result of this section. It, or more precisely its contrapositive, gives a sufficient condition of the conductive homogeneity.

\thm\label{MIS.thm20}
Assume that $J = 2q$ for some $q \ge 3$, $G = D_q$ or $Rot_q$, and that there exists no isolated contact point of cells. If there exists $p > \dim_{AR}(K, d_*)$ such that $(K, d_*)$ is not $p$-conductive homogeneous, then either
\begin{equation}\label{MIS.eq500}
\BbZ_J^i \in F_{\partial}(\BbZ_J^i) 
\end{equation}
for some $i \in \{0, 1\}$
or
\begin{equation}\label{MIS.eq510}
\BbZ_J^1 \in F_{\partial}(\BbZ_J^0)\quad\text{and}\quad \BbZ_J^0 \in F_{\partial}(\BbZ_J^1).
\end{equation}
Moreover, if $J = 6$ or $8$, then \eqref{MIS.eq510} does not occur.
\endthm

\remark
Under the assumption of Theorem~\ref{MIS.thm20}, it follows that $E_n^* = E_n^{\ell}$ for any $n \ge 1$ by Lemma~\ref{MIS.lemma05}.
\endremark

The statement of the above theorem is not so complicated but a proof, which is given in the next two sections, requires quite involved arguments and occupies a considerable number of pages.\par

\begin{figure}[ht]
\centering
\includegraphics[width = 300pt]{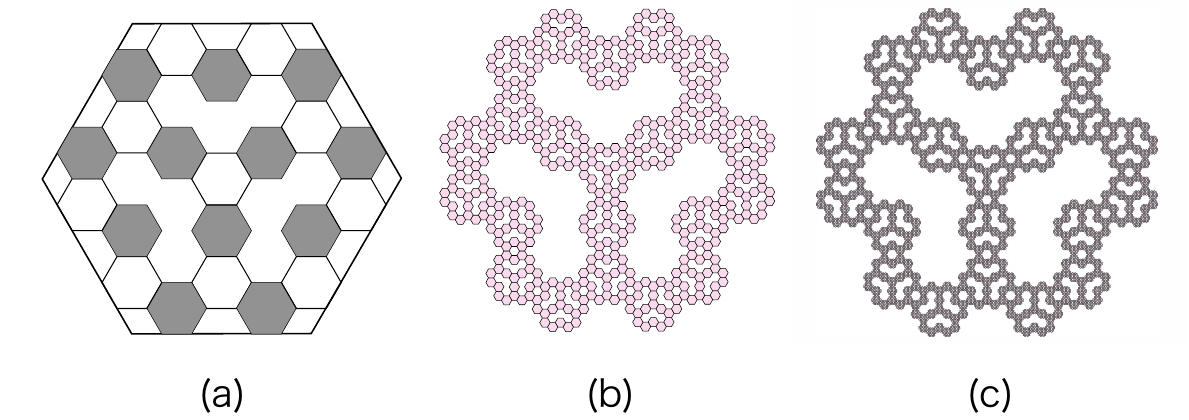}
\caption{$J = 6$, $G = \D_3$}\label{HexaD3a}
\end{figure}

\begin{figure}[ht]
\centering
\includegraphics[width = 200pt]{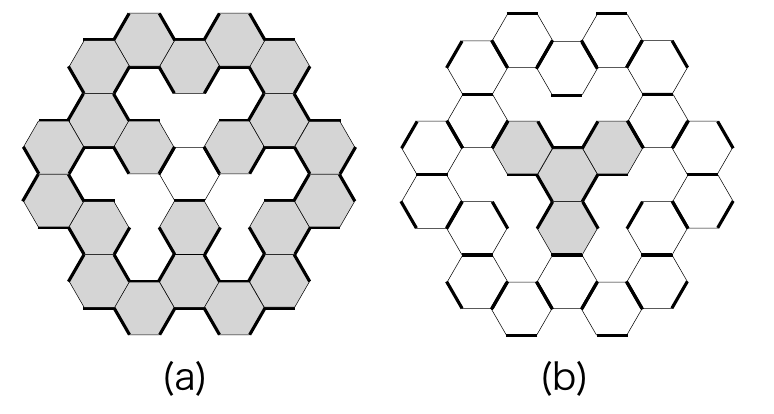}
\caption{$\Con_T(\BbZ_6^0)$ and $\Con_T(\BbZ_6^1)$}\label{HexaD3c}
\end{figure}

\example[Figure~\ref{HexaD3a}-(c)]\label{MIS.ex00}
Let $J = 6$. Define $S$ as the collection of small white and grey hexagons in Figure~\ref{HexaD3a}-(a) and define
\[
\vp_s = \begin{cases}
I\quad&\text{if the corresponding hexagon $s$ is white,}\\
R_0\quad&\text{if the corresponding hexagon $s$ is grey.}
\end{cases}
\]
Then $(S, \{f_s\}_{s \in S}, G)$ is a $(6, G)$-s.s.\,system with $G = D_3$. The middle figure (b) represents the second stage $\bigcup_{w \in T_2} Q_w$ and the right one (c) represents the corresponding self-similar set $K$. In this case, $(\BbZ_6)^e = \BbZ_6$. The connected components $\Con_T(\BbZ_6^0)$ is represented by Figure~\ref{HexaD3c}-(a), where thick lines are removed edges corresponding to the images of $\BbZ_6^0$. Consequently, there are only two connected components, say $A_1$ and $A_2$, which are colored light grey and white, respectively. Then $\partial_{\BbZ_6^1}A_1 = \BbZ_6$ and $\partial_{\BbZ_6^1}A_2 = \emptyset$. Hence it follows that $F_{\partial}(\BbZ_6^0) = \{\BbZ_6, \emptyset\}$. In the case of $\Con_T(\BbZ_6^1)$ illustrated in Figure~\ref{HexaD3c}-(b), each connected component consists of a single cell except the light grey region. Anyway, we see that $\partial_{\BbZ_6^0}A = \emptyset$ for any connected component $A \in \Con_T(\BbZ_6^1)$ and hence $F_{\partial}(\BbZ_6^1) = \{\BbZ_6\}$. Thus Theorem~\ref{MIS.thm20}-(b) shows that $(K, d_*)$ is $p$-conductively homogeneous for any $p > \dim_{AR}(K, d_*)$.

\endexample

\begin{figure}[ht]
\centering
\includegraphics[width = 300pt]{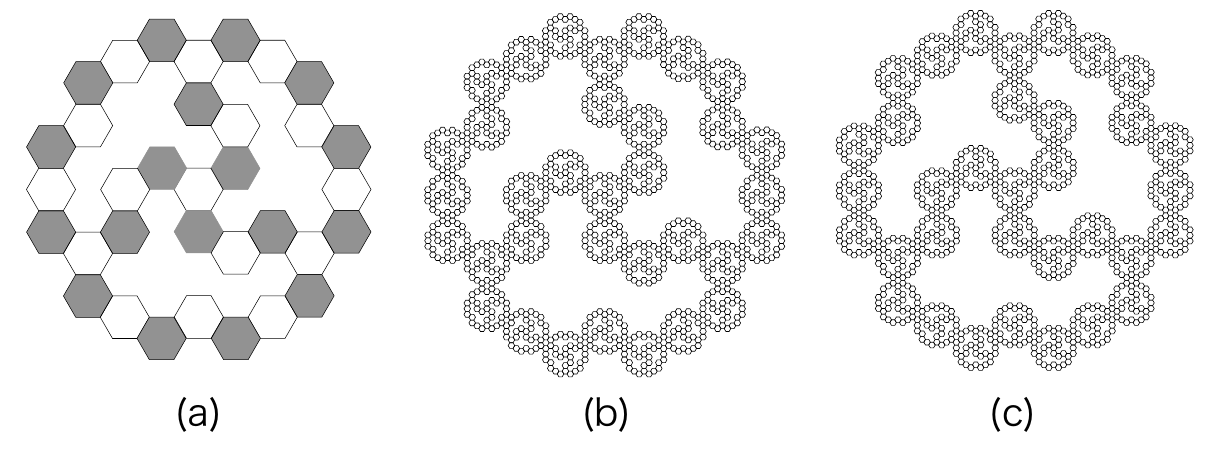}
\caption{Hexagon-based isomers: $J = 6$, $G = Rot_3$}\label{HexaJH0}
\end{figure}

\begin{figure}
\begin{minipage}[b]{6cm}
\centering
\includegraphics[width = 120pt]{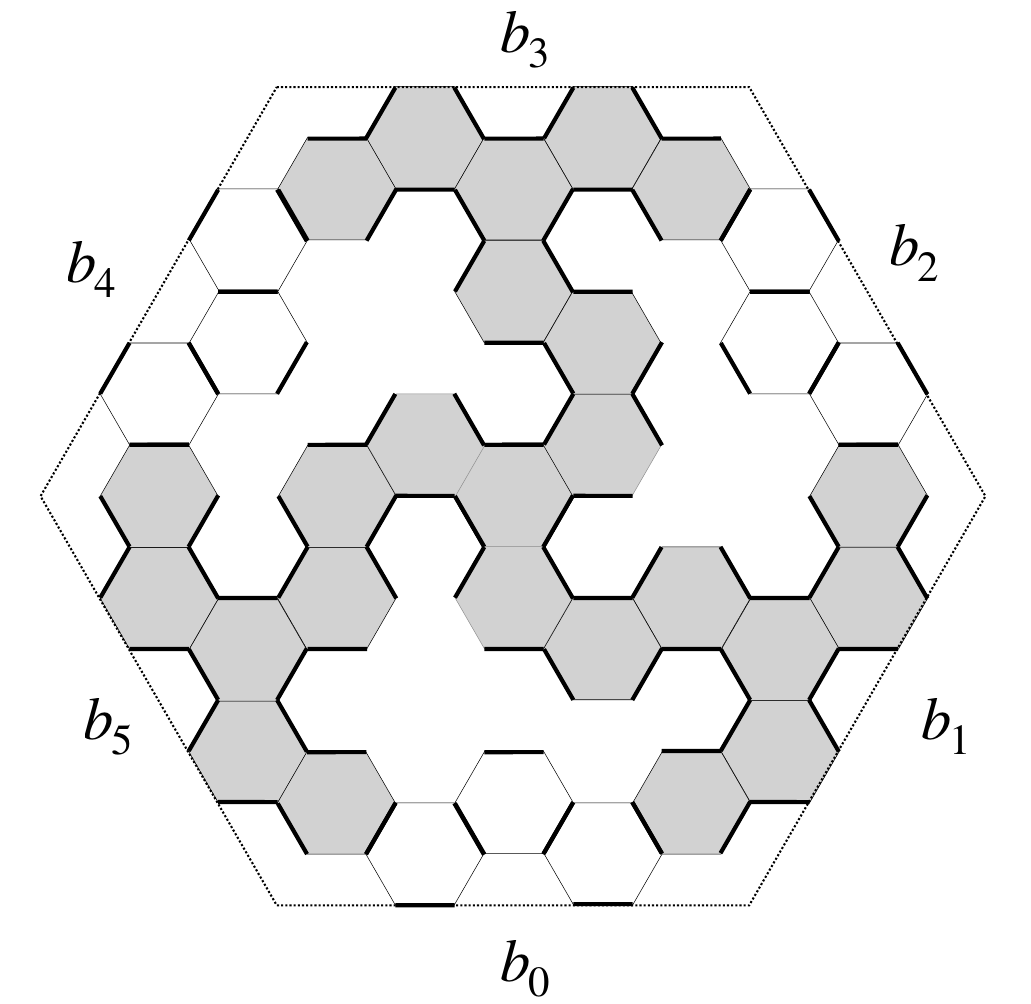}
\end{minipage}
\begin{minipage}[b]{6cm}
\hspace{-20pt}
\centering
\includegraphics[width=114pt]{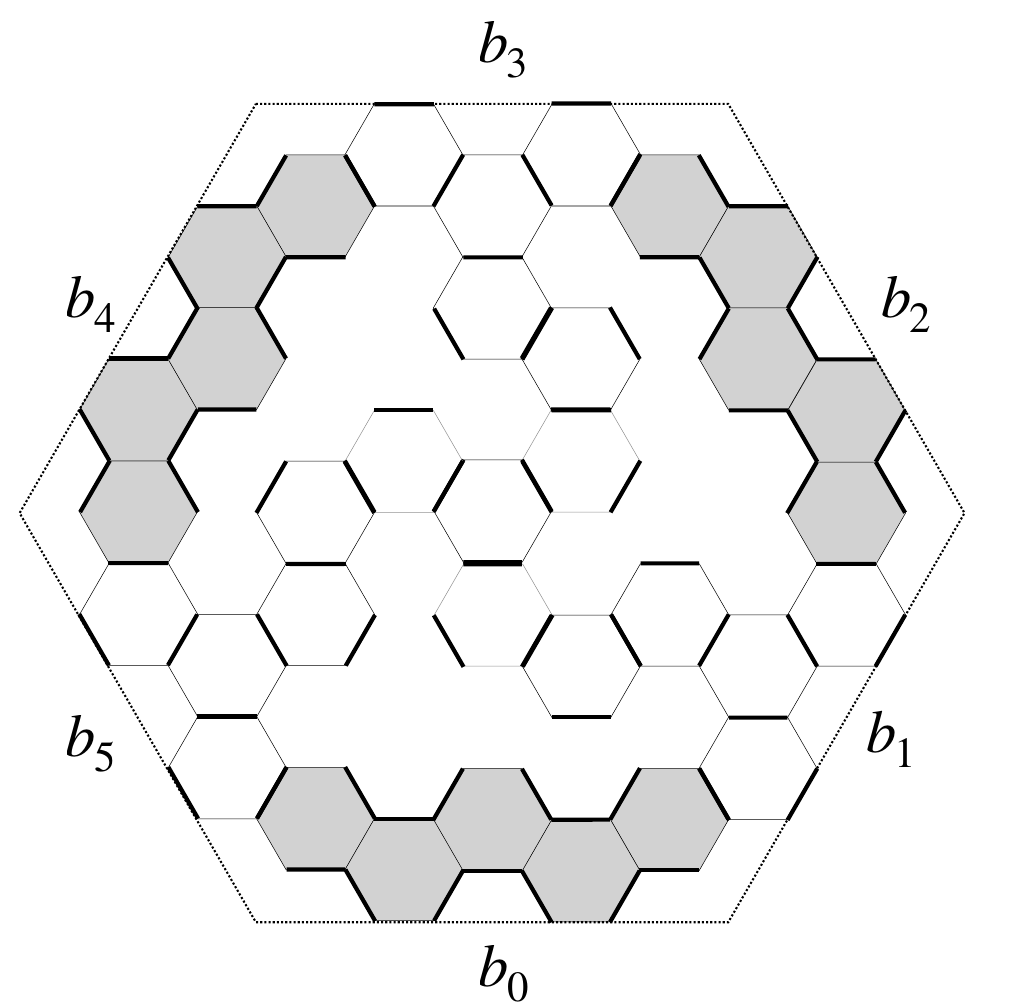}
\vspace{3pt}
\end{minipage}
\caption{$\Con_T(\BbZ_J^0)$ and $\Con_T(\BbZ_J^1)$}\label{HexaJH1}
\end{figure}

\example[Isomers]\label{MIS.ex10}
Let $J = 6$ and let $G = Rot_3$. Define $S$ as the collection of small hexagons in Figure~\ref{HexaJH0}-(a). The contraction ratio $r$ is determined by the configuration of all the hexagons in $S$. To be exact, this example is two examples, Jekyll and Hyde, sharing the same $S$, but having quite different natures on $F_{\partial}(\BbZ_6^i)$. See Figure~\ref{HexaJH0}. They are a kind of isomers. We present $\Con_T(\BbZ_6^0)$ and $\Con_T(\BbZ_6^1)$ for these two isomers in Figure~\ref{HexaJH1}, where thick lines are disconnected and light grey regions correspond to some of the components of $\Con_T(\BbZ_J^i)$. Among two figures in Figure~\ref{HexaJH1}, one of them represents $\Con_T(\BbZ_6^0)$ and the other represents $\Con_T(\BbZ_6^1)$. ``Which is which'' depends on an isomer.\\
{\bf Case Jekyll}(Figure~\ref{HexaJH0}-(b))\,\,\,In this case,
\[
\vp_s = \begin{cases}
I\quad&\text{if the hexagon $s$ is colored white in Figure~\ref{HexaJH0}-(a),}\\
R_0\quad&\text{if the hexagon $s$ is colored grey in Figure~\ref{HexaJH0}-(a).}
\end{cases}
\]
for $s \in S$. The right figure of Figure~\ref{HexaJH1} represents $\Con_T(\BbZ_6^0)$. Let $A_j$ be the element of $\Con_T(\BbZ_6^0)$ corresponding to the light grey region intersecting with $b_j$ for $j = 0, 1, 2$. Then $\partial_{\BbZ_6^1}A_j = \{j\}$ and hence $\{j\}^c = \sd{\BbZ_J}{\{j\}} \in F_{\partial}(\BbZ_J^0)$. For any other $A \in \Con_T(\BbZ_6^0)$, we see that $\partial_{\BbZ_6^1}A = \emptyset$. Consequently, 
\[
F_{\partial}(\BbZ_6^0) = \{\{0\}^c, \{2\}^c, \{4\}^c, \BbZ_6\}.
\]
On the other hand, the left figure represents $\Con_T(\BbZ_6^1)$. Let $A_*$ be the element of $\Con_T(\BbZ_6^1)$ corresponding to the large light grey region. Then $\partial_{\BbZ_J^0}A_* = \BbZ_6^1$ and hence $\BbZ_6^0 \in F_{\partial}(\BbZ_J^1)$. For any other $A \in \Con_T(\BbZ_6^1)$, we have $\partial_{\BbZ_6^0}A = \emptyset$. Thus, it follows that
\[
F_{\partial}(\BbZ_6^1) = \{\BbZ_6^0, \BbZ_6\}.
\]
Now Theorem~\ref{MIS.thm20} shows that $(K, d_*)$ is $p$-conductively homogeneous for any $p > \dim_{AR}(K, d_*)$.\\
{\bf Case Hyde}(Figure~\ref{HexaJH0}-(c))\,\,\,In this case,
\[
\vp_s = \begin{cases}
R_0\quad&\text{if the hexagon $s$ is colored grey in Figure~\ref{HexaJH0}-(a),}\\
I\quad&\text{if the hexagon $s$ is colored white in Figure~\ref{HexaJH0}-(a).}
\end{cases}
\]
for $s \in S$. The correspondence between left and right figures of Figure~\ref{HexaJH1} and $\Con_T(\BbZ_6^i)$ for $i \in \BbZ_2$ is now reversed. By a similar argument to Case Jekyll, we see that
\[
F_{\partial}(\BbZ_6^0) = \{\BbZ_6^0, \BbZ_6\}
\]
and
\[
F_{\partial}(\BbZ_6^1) = \{\{0\}^c, \{2\}^c, \{4\}^c, \BbZ_6\}.
\]
Not that $\BbZ_6^0 \in F_{\partial}(\BbZ_6^0)$, we can not tell anything about the conductive homogeneity by Theorem~\ref{MIS.thm20}.
\endexample

\example[= Examle~\ref{GPS.ex20}]\label{MIS.ex20}
Let $J = 8$ and let $(S, \{f_s\}_{s \in S}, Rot_4)$ be the $(8, Rot_4)$-s.s.\,system introduced in Example~\ref{GPS.ex20}. In this case, $(\BbZ_8)^e = \BbZ_8$. Moreover, it follows that
\[
F_{\partial}(\BbZ_8^0) = \{\BbZ_8, \{1\}^c, \{3\}^c, \{5\}^c, \{7\}^c\}
\]
and
\[
F_{\partial}(\BbZ_8^1) = \{\BbZ_8, \{0\}^c, \{2\}^c, \{4\}^c, \{6\}^c\}.
\]
Hence Theorem~\ref{MIS.thm20}-(2) yields $p$-conductive homogeneity of $(K, d_*)$ for any $p > \dim_{A}(K, d_*)$.
\endexample

\section{Bondary correspondence}\label{BCM}

In this section, we study properties of paths of $(T_n, E_n^{\ell}(X))$ for $X \subseteq \BbZ_J$ to prepare proofs of Theorem~\ref{MIS.thm20} in the previous section. Also, the results in this section will be used in Sections~\ref{WNG} and \ref{POT}, where $J$ is not necessarily even. Throughout this section, $J \ge 3$ and $(S, \{f_s\}_{s \in S}, G)$ is assumed to be a $(J, G)$-s.s.\,system, where $G$ is a subgroup of $D_J$. Note that we don't assume that $J$ is even.                                                                                                                

\lemma\label{BCM.lemma10}
Let $X$ be a $G$-invariant subset of $\BbZ_J$. \\
{\rm (1)}\,\,Let $(w, v) \in E_n^{\ell}$. Suppose that $Q_w \cap Q_v = b_{i_1}(w) = b_{i_2}(v)$ for some $i_1, i_2 \in \BbZ$. Then $i_1 \in X$ if and only if $i_2 \in X$.\\
{\rm (2)}\,\,$E_n^{\ell} = E_n^{\ell}(X) \cup E_n^{\ell}(X^c)$ and $E_n^{\ell}(X) \cap E_n^{\ell}(X^c) = \emptyset$.\\
{\rm (3)}\,\,For any $w \in T$, $(u, v) \in E_n^{\ell}(X)$ if and only if $(wu, wv) \in E_{n + |w|}^{\ell}(X)$.\\
{\rm (4)}\,\,For any $w \in T$ and $g \in G$, there exists $h \in G$ such that
\[
g(b_i(w)) = b_{h(i)}(g(w))
\]
for any $i \in \BbZ_J$. In particular if $(w, v) \in E_n^{\ell}(X)$, then $(g(w), g(v)) \in E_n^{\ell}(X)$.
\endlemma

\demo
(1)\,\,Since $g_{w, v} \in G$, $g_{w, v}(i_2) = i_1$ and $X$ is $G$-invariant, we have the desired result.\\
(2)\,\,This is immediate from (1).\\
(3)\,\,If $Q_u \cap Q_v = b_u(i_1) = b_u(i_2)$, then $Q_{wu} \cap Q_{wv} = b_{wu}(i_1) = b_{wv}(i_2)$. Now it is obvious.\\
(4)\,\,By \eqref{BAS.eq100}, letting $h = (f_{g(w)})^{-1}{\circ}g{\circ}f_w$, we see that $h \in G$. So, $g(b_{i_1}(w)) = f_{g(w)}(h(b_{i_1})) = b_{h(i_1)}(g(w))$. In particular, since $X$ is $G$-invariant, it follows that $h(i) \in X$ whenever $i \in X$. This immediately implies the final part of the statement.
\enddemo

\definition\label{BCM.def20}
Let $X \subseteq \BbZ_J$. Define
\begin{multline*}
\C^{\ell}_{m, \partial}(X) = \{\c|\text{$\c$ is a path of $(T_m, E_m^{\ell}(X))$}, \\\text{and there exist $i_1, i_2 \in X$ and $k_1, k_2 \in (\BbZ_J)^e$ such that}\\
\text{$b_{i_1}(\c(1)) \subseteq b_{k_1}$ and $b_{i_2}(\c(l(\c))) \subseteq b_{k_2}$.}\}
\end{multline*}
In particular, we use $\C^{\ell}_{m, \partial}$ to denote $\C^{\ell}_{m, \partial}(\BbZ_J)$.
\enddefinition

\lemma\label{BCM.lemma20}
Let $X$ be a $G$-invariant subset of $\BbZ_J$. Let $\c$ be a path of $(T_m, E_m^{\ell}(X))$ and let $m = m_1 + m_2$ for some $m_1, m_2 \in \BbZ$. Assume that $l_{m_1}(\c) \ge 3$. Then $\s_{m_1}(\c)_j \in \C_{m_2, \partial}^{\ell}(X)$ for any $j \in 2, \ldots, l_{m_1}(\c) - 1$.
\endlemma

\demo
Let $[\c]_{m_1} = (w(1), \ldots, w(l_*))$ and let $\s_{m_1}(\c)_j = (v_j(1), \ldots, v_j(l_j))$ for any $j = 1, \ldots, l_*$. Note that
\[
\c = (w(1)v_1(1), \ldots, w(1)v_1(l_1), \ldots, w(l_*)v_{l_*}(1), \ldots, w(l_*)v_{l_*}(l_{l_*})).
\]
By Lemma~\ref{BCM.lemma10}-(3), since $(w(j)v_j(k), w(j)v_j(k + 1)) \in E_m^{\ell}(X)$, we see that $(v_j(k), v_j(k + 1)) \in E_{m_2}^{\ell}(X)$. Moreover, Since $(w(j - 1)v_{j - 1}(l_{j - 1}), w(j)v_j(1)) \in E_m^{\ell}(X)$ for any $j = 2, \ldots, l_* - 1$, there exists $k \in (\BbZ_J)^e$ and $i \in X$ such that $b_i(w(j)v_j(1)) \subseteq b_k(w(j))$. This implies $b_i(v_j(1)) \subseteq b_k$. Similarly, for any $j = 2, \ldots, l_* - 1$, there exist $i \in X$ and $k \in (\BbZ_J)^e$ such that $b_i(v_j(l_j)) \subseteq b_k$. Thus we have shown $\s_{m_1}(\c)_j =(v_j(1), \ldots, v_j(l_j)) \in \C^{\ell}_{m_2, \partial}(X)$.
\enddemo

\lemma\label{BCM.lemma30}
Let $X$ be a $G$-invariant subset of $\BbZ_J$. Let 
\[
\widetilde{\c} = (\c(0), \c(1), \ldots, \c(l_*), \c(l_* + 1))
\]
be a path of $(T_{n + m}, E_{n + m}^{\ell})$ with $[\c(0)]_n \neq [\c(1)]_n$ and $[\c(l_*)]_n \neq [\c(l_* + 1)]_n$. Set $\c = (\c(1), \ldots, \c(l_*))$. Assume that $(\partial_{\ell}\s_n^F(\c))^e \subseteq X$. Then $[\widetilde{\c}]_n$ is a path of $(T_n, E_n^{\ell}(X))$. 
\endlemma

\demo
Let $[\c]_n = (w(1), \ldots, w(l))$. Then $[\widetilde{\c}]_n = (w(0), w(1), \ldots, w(l), w(l + 1))$, where $w(0) = [\c(0)]_n$ and $w(l + 1) = [\c(l_* + 1)]_n$. Note that 
\[
\s_n^F(\c) = \s_n^F(\c)_1 \vee g_2(\s_n(\c)_2) \vee \ldots \vee g_{l_*}(\s_n(\c)_{l})
\]
for some $g_2, \ldots, g_{l} \in G$. Since $(\partial_{\ell}\s_n^F(\c))^e \subseteq X$, Lemma~\ref{BCM.lemma10}-(4) shows that $(\partial_{\ell}\s_n(\c)_j)^e \subseteq X$ for any $j = 1, \ldots, l$. Let $v_1(j)$ (resp. $v_L(j)$) be the first (resp. the last) piece of $\s_n(\c)_j$. Then there exist $i_1(j), i_L(j) \in \BbZ_J$ and $k_1(j), k_L(j) \in (\BbZ)^e$ such that $b_{i_1(j)}(v_1(j)) \subseteq b_{k_1(j)}$ and $b_{i_L(j)}(v_L(j)) \subseteq b_{k_L(j)}$. This implies that $k_1(j), k_L(j) \in (\partial_{\ell}\s_n(\c)_j)^e \subseteq X$. Since 
\[
Q_{w(j)} \cap Q_{w(j + 1)} = \begin{cases}
b_{k_1(j + 1)}(w(j + 1))\quad &\text{for $j = 0, \ldots, l - 1$,}\\
b_{k_2(j)}(w(j))\quad &\text{for $j = 1, \ldots, l$,}
\end{cases}
\]
we see that $(w(j), w(j + 1)) \in E_n^{\ell}(X)$ for $j = 0, 1, \ldots, l$.
\enddemo

\lemma\label{BCM.lemma40}
Let $X$ be a $G$-invariant subset of $\BbZ_J$ and let $q_0, q, q_1, q_2 \in \BbN$ satisfying $q < q_0$ and $q = q_1 + q_2$. Let $\widetilde{\c} = (\c(0), \c(1),\ldots, \c(l_*)), \c(l_* + 1))$ be a path of $(T_{q_0}, E_{q_0}^{\ell})$ with $[\c(0)]_q \neq [\c(1)]_q$ and $[\c(l_*)]_q \neq [\c(l_* + 1)]_q$. Define $\c = (\c(1), \ldots, \c(l_*))$. Assume that $(\partial_{\ell}\s_q^F(\c))^e \subseteq X$. Then $g(\s_{q_1}([\c]_q)_j) \in \C_{q_2, \partial}^{\ell}(X)$ for any $j = 1, \ldots, l_{q_1}([\c]_q)$ and $g \in G$. Moreover, $\s_{q_1}^F([\c]_q) = [\s_{q_1}^F(\c)]_{q_2}$ is a path of $(T_{q_2}, E_{q_2}^{\ell}(X))$ and $(\partial_{\ell}\s_{q_1}^F(\c))^e \subseteq \partial_X([\s_{q_1}^F(\c)]_{q_2})$.
\endlemma

\demo
Applying Lemma~\ref{BCM.lemma30} with $n = q$ and $m = q_0 - q$, we see that $[\widetilde{\c}]_q$ is a path of $(T_q, E_q^{\ell}(X))$. Note that $[\widetilde{\c}]_q = ([\c(0)]_q, [\c(1)]_q, \ldots, [\c(l_*)]_q, [\c(l_* + 1)]_q$. Set $m_1 = q_1$ and $m_2 = q_2$. Then Lemma~\ref{BCM.lemma20} yields that $\s_{q_1}([\c]_q)_j \in \C_{q_2, \partial}^{\ell}(X)$ for any $j = 1, \ldots,l_{q_1}([\c]_q)$. Let $[\c]_{q_1} = (w(1), \ldots, w(l))$ and let $[\s_{q_1}(\c)_j]_{q_2} = (v_j(1), \ldots, v_j(k_j))$ for $j = 1, \ldots, l$. Then 
\[
[\c]_q = (w(1)v_1(1), \ldots, w(1)v_1(k_1), w(2)v_2(1), \ldots, w(l)v_l(1), \ldots, w(l)v_l(k_l)).
\]
Hence $l_{q_1}([\c]_q) = l$ and $\s_{q_1}([\c]_q)_j = (v_j(1), \ldots, v_j(k_j)) = [\s_{q_1}(\c)_j]_{q_2}$ for any $j = 1, \ldots, l$. Since 
\[
\s_{q_1}^F([\c]_q) = \s_{q_1}([\c]_q)_1 \vee g_2(\s_{q_1}([\c]_q)_2) \vee \ldots \vee g_l(\s_{q_1}([\c]_q)_l)
\]
for some $q_2, \ldots, q_l \in G$, we see that $\s_{q_1}([\c]_q)$ is a path of $(T_{q_2}, E^{\ell}_{q_2}(X))$. Note that
\[
\s_{q_1}^F(\c) = \s_{q_1}(\c)_1 \vee g_2(\s_{q_1}(\c)_2) \vee \ldots \vee g_l(\s_{q_1}(\c)_l).
\]
For any $\b \in (\partial_{\ell}\s_{q_1}^F(\c))^e$, there exist $j \in \{1, \ldots, l\}$, $\tau \in g_j(\s_{q_1}(\c)_j)$ and $\a \in \BbZ_J$ such that $b_{\a}(\tau) \subseteq b_{\b}$. From now on, we write $g = g_j$ for simplicity. Since $[\tau]_{q_2} = g(v_j(k))$ for some $k$, there exists $\a' \in \BbZ_J$ such that 
\[
b_{\a}(\tau) \subseteq b_{\a'}(g(v_j(k))) \subseteq b_{\b}.
\]
This implies that $\a' \in (\partial_{\ell}\s_q^F(\c))^e \subseteq X$. Thus, $\b \in \partial_X([\s_{q_1}^F(\c)]_{q_2})$.
\enddemo

\section{Proof of Theorem~\ref{MIS.thm20}}\label{PTE}

This section is devoted to a proof of Theorems~\ref{MIS.thm20}. In this section we assume that $J = 2q$ with $q \ge 3$ and that $(S, \{f_s\}_{s \in S}, G)$ is a $(J, G)$-s.s.\,system, where $G$ is a subgroup of $D_J$.\par

\definition\label{PTE.def10}
Let $A \subseteq T_n$ for some $n \ge 0$. $A$ is said to connect $0$ and $1$ if  there exists $i \in \partial{A} \cap \BbZ_J^0$ and $j \in \partial{A} \cap \BbZ_J^1$ such that $b_i \cap b_j = \emptyset$. 
\enddefinition

The following theorem is a corollary of Theorem~\ref{BAS.thm10}.

\thm\label{PTE.thm10}
Let $j_*, n_* \ge 1$. Suppose that $Rot_q \subseteq G$ and that for any $w \in T$ and $\c \in \C_{M, m}(w)$, there exists 
\[
A \subseteq \bigcup_{j \in \{-j_*, \ldots, 0, \ldots,  j_*\}, n_2 \in \{0, 1, \ldots,n_*\}} \H_{|w| + j, n_2, m}(\c)
\]
 such that $K(A)$ is connected and $A$ connects $0$ and $1$. Then $(K, d_*)$ is $p$-conductively homogeneous for any $p > \dim_{AR}(K, d_*)$. Moreover, if there exists no isolated contact point of cells, then $\C_{M, m}(w)$ in the above statement is replaced by $\C_{M, m}^{\ell}(w)$.
\endthm

\demo
Since $A$ connects $0$ and $1$, there exist $s \in \BbZ_J^0 \cap \partial{A}$ and $t \in \BbZ_J^1 \cap \partial{A}$ such that $b_k \cap b_l = \emptyset$. This implies that $\d(k, l) \ge 3$. Applying Lemma~\ref{BAS.lemma400}-(2) with $q_1 = q$ and $q_2 = 2$, we see that $Rot_q(A)$ is $(T_m, E_m^*)$-connected. Moreover, since $\BbZ_J = \{g(k), g(l)| g \in Rot_q\} \subseteq \partial(Rot_q(A))$, Theorem~\ref{BAS.thm10} yields the desired conclusion.
\enddemo

Now we proceed to a proof of Theorem~\ref{MIS.thm20}. In the rest of this section, we assume that $J = 2q$ for some $q \ge 3$, $G \subseteq D_q$ and there exists no isolated contact point of cells.\par

\lemma\label{MIS.lemma10}
Let $J = 6$ and let $n \ge 1$. If there exists $w(1), w(2), w(3) \in T_n$ such that $Q_{w(1)} \cap Q_{w(2)} \cap Q_{w(3)}$ is a single point, then there exists $i \in \BbZ_6$ such that $R_{\theta_i} \in G$, where $\theta_i = \frac{\pi}3(i + 2)$. In particular, $G$ is not a subgroup of $D_3$.
\endlemma

\demo
Let $Q_{w(1)} \cap Q_{w(2)} \cap Q_{w(3)} = \{f_{w(1)}(p_i)\}$. Set $g_j = g_{w(j + 1), w(j)}$ for $j = 1, 2, 3$, where we read $j + 1$ as $1$ for $j = 3$. By (B4)-(a), it follows that $g_j \in G$ for any $i \in \{1, 2, 3\}$. Moreover
\begin{align*}
g_3g_2g_1(p_i) &= f_{w(1)}^{-1}{\circ}R_{w(3), w(1)}{\circ}R_{w(2), w(3)}{\circ}R_{w(1), w(2)}{\circ}f_{w(1)}(p_i)\\
&= f_{w(1)}^{-1}{\circ}f_{w(1)}(p_i) = p_i
\end{align*}
Hence Letting $g = g_3g_2g_1$, we see that $g(p_i) = g(p_i)$ and $\det(g_i) = -1$, where $\det(A)$ is the determinant of a matrix $A$. This shows $g = R_{\theta_i}$. Since $g \in G$, we have shown the lemma.
\enddemo

\lemma\label{MIS.lemma05}
If $J \ge 6$, $G \subseteq D_q$ and there exists no isolated contact point of cells, then $E_n^* = E_n^{\ell}$ for any $n \ge 1$ and $Q_{w(1)} \cap Q_{w(2)} \cap Q_{w(3)} = \emptyset$ for any $n \ge 1$ and any three distinct elements $w(1), w(2), w(3) \in T_n$.
\endlemma

\demo
Proposition~\ref{COP.prop10} shows that $E_n^* = E_n^{\ell}$. The rest of the statement is obvious in the case $J \ge 8$. If $J = 6$, it follows from Lemma~\ref{MIS.lemma10}.
\enddemo

\lemma\label{MIS.lemma20}
Let $\c = (\c(1), \ldots, \c(l)) \in \C_{M, m}^{\ell}(w)$ for some $M, m \ge 1$ and $w \in T$. Set $n = |w|$ and let $l_* = l_{n + j}(\c)$. For any $j \in \{1, \ldots, m\}$, $[\c]_{n + j} \in \C_{M, j}^{\ell}(w)$, $\s_{n + j}(\c)_i \in \C_{m - j, \partial}^{\ell}$ for any $i = 1, \ldots, l_*$ and
\[
\H_{n + j, 0, m}(\c) =  S^j\bigg(\bigcup_{g \in G, i = 1, \ldots, l_*}g(\s_{n + j}(\c)_i)\bigg).
\]
Moreover, let $[\c]_{n + j} = (w(1), \ldots, w(l_*))$ and define
\[
g_i = \begin{cases}
\text{the identity}&\quad\text{if $i = 1$,}\\
g_{w(i - 1), w(i)}&\quad\text{if $i = 2, \ldots, l_*$},
\end{cases}
\]  
\[
h_i = g_1{\circ}\cdots{\circ}g_i
\]
for $i = 1, \ldots, l_*$ and
\[
A_j(\c) = \underset{i = 1, \ldots, l_*}{\vee} h_i(\s_{n + j}(\c)_i). 
\]
Then 
\begin{equation}\label{MIS.eq10}
h_i\big(\s_{n + j}(\c)_i(l_i)\big) = h_{i + 1}\big(\s_{n + j}(\c)_{i + 1}(1)\big)
\end{equation}
for $i = 1, \ldots, l_* - 1$, where $l_i = l(\s_{n + j}(\c)_i)$. Moreover  $A_j(\c)$ is $(T_{m - j}, E_{m - j}^{\ell})$-connected and 
\[
\H_{n + j, 0, m}(\c) = S^j\bigg(\bigcup_{g \in G} g(A_j(\c))\bigg).
\]
\endlemma

\remark
Recalling Lemma~\ref{BAS.lemma200}, we see that $A_j(\c)$ is the $(n + j)$-folding of $\c$, $\s_{n + j}^F(\c)$.
\endremark
\remark
By definition, the element $\s_{n + j}(\c)_i(l_i)$ and $\s_{n + j}(\c)_{i + 1}(1)$  are the last piece of the path $\s_{n + j}(\c)_i$ and the first one of the path $\s_{n + j}(\c)_{i + 1}$.
\endremark

\demo
Since $\c \in \C_{M, m}^{\ell}(w)$,  we see that  $\s_{n + j}(\c)_i \in \C_{m - j, \partial}^{\ell}$ and $h_i(\s_{n + j}(\c)_i) \in \C_{m - j, \partial}^{\ell}$ for any $i = 1, \ldots, l_*$.  On the other hand, Lemma~\ref{BAS.lemma201} shows
\[
w(1)h_i(\s_{n + j}(\c)_i)) = R_{w(1), w(2)}^*{\circ}\cdots{\circ}R_{w(i - 1), w(i)}^*(w(i)\s_{n + j}(\c)_{i})
\]
for any $i = 1, \ldots, l_*$. Moreover, since $\c(i') = w(i)\s_{n + j}(\c)_i(l_i)$ and $\c(i' + 1) = w(i + 1)\s_{n + j}(\c)_{i + 1}(1)$ for some $i'$,  it follows that
\[
R_{w(i), w(i + 1)}^*\big(w(i + 1)\s_{n + j}(\c)_{i + 1}(1)\big) =  w(i)\s_{n + j}(\c)_i(l_i), 
\]
so that we have \eqref{MIS.eq10}. The rest of statements is straightforward by the fact that $A_j(\c) = \s_{n + j}^F(\c)$.
\enddemo

\lemma\label{MIS.lemma30}
Assume that $J = 2q$ for some $q \ge 3$, $G \subseteq D_q$ and there exists no isolated contact point of cells. Let $M \ge M_J$, let $\c \in \C_{M, m}^{\ell}(w)$ and let $n = |w|$. Then for any $j \in \{1, \ldots, m\}$, one of the following three cases {\rm (a), (b)} or {\rm (c)} happens:\\
{\rm (a)}\,\,$(\partial{A_j(\c)})^e = \{s, s + J/2\}$ for some $s \in \BbZ_J$.\\
{\rm (b)}\,\,$\#((\partial{A_j(\c)})^e) = 3$ and $\bigcup_{s \in (\partial{A_j(\c)})^e} b_s$ is disconnected\\
{\rm (c)}\,\,$\#((\partial{A_j(\c)})^e) \ge 4$.
\endlemma

\demo
Let $\gamma = (\c(1), \ldots, \c(l))$ and choose $\c(0)\in S^m(w)$ and $\c(l + 1)\in S^{m}(\Gamma_M(w)^c)$ satisfying $(\c(0),\c(1)),(\c(l),\c(1+1))\in E_{n+m}^l$.
We use the same notations as in Lemma~\ref{MIS.lemma20} and its proof. For example, $l_* = l_{n + j}(\c)$ and $[\c]_{n + j} = (w(1), \ldots, w(l_*))$. First we modify $\c$ to avoid the occurrence of $i$ such that $w(i - 1) = w(i + 1)$. Let
\[
D_j(\c)  = \{i| w(i - 1) = w(i + 1), i = 1, \ldots, l_*\}
\]
and
\[
i(j, \c) = \min D_j(\c)\quad\,\,\text{if $D_j(\c) \neq \emptyset$}.
\]
Note that $i(j, \c) > 1$ because $w(0) \in S^m(w)$ and $w(2) \notin S^m(w)$ by the definition of $\C_{M, m}(w)$. If $D_j(\c) \neq \emptyset$, define $\xi(\c) \in \C_{M, m}(w)$ by replacing $w(i')\s_{n + j}(\c)_{i'}$ with $w(i' - 1)(g_{i'})^*(\s_{n + j}(\c)_{i'})$, where $i'= i(j, \c)$. Then $l_{n + j}(\xi(\c)) = l_{n + j}(\c) - 2$ and $[\xi(\c)]_{n + j} = (w(1), \ldots, w(i'- 1), w(i' + 2), \ldots, w(l_*))$. Moreover, we see that $\H_{n + j, 0, m}(\c) = \H_{n + j, 0, m}(\xi(\c))$ and $A_j(\c) = A_j(\xi(\c))$. Now, iterating $\xi$, we eventually see that $D_j(\xi^N(\c)) = \emptyset$. Replacing $\c$ with $\xi^N(\c)$, we may assume that $D_j(\c) = \emptyset$ without loss of generality.

Set $w(0) = [\c(0)]_{n + j}$ and $w(l_* + 1) = [\c(l + 1)]_{n + j}$. Moreover, define $R_i$ as $R_{w(i), w(i + 1)}$. There exist $\{s_i\}_{i=1}^{\ell_* +1},\{t_i\}_{i=0}^{\ell_*} \subset \BbZ_J$ such that
$Q_{w(i)} \cap Q_{w(i + 1)} = b_{t_i}(w(i)) = b_{s_{i + 1}}(w(i + 1))$ for any $i = 0, \ldots, l_*$. Then
\[
\quad R_i(b_{t_i}(w(i))) = b_{s_{i + 1}}(w(i + 1))
\]
for any $i = 0,1, \ldots, l_*$.
Let $Y=\bigcup_{i=1}^{\ell_*}h_i(\{s_{i},t_{i}\})$, then $Y \subseteq (\partial{A_j(\c)})^e$. 

\noindent Claim: $\d(s_i, t_i) \ge 2$ for any $i = 1, \ldots, l_*$.

\noindent Proof of Claim: If $\d(s_i, t_i) = 1$, then $Q_{w(i - 1)} \cap Q_{w(i)} \cap Q_{w(i + 1)} \neq \emptyset$. This contradicts Lemma~\ref{MIS.lemma05}. Hence, we have the desired statement.

By the above claim, it follows that $\#((\partial{A_j(\c)})^e) \geq 2$. In the case of $\#((\partial{A_j(\c)})^e)  \ge 4$, the statement(c) is satisfied. So, we consider the cases $\#((\partial{A_j(\c)})^e)  = 2$ and $\#((\partial{A_j(\c)})^e)  = 3$. Assume that $\#((\partial{A_j(\c)})^e)  = 3$. If $\bigcup_{s \in Y} b_s$ is connected, then there exists $k \in \BbZ_J$ such that $Y=\{k-1,k,k+1\}$. On the other hand, the above claim implies that there exists $t \in Y$ such that $\d(s, t) \ge 2$. By this contradiction, we see that $\cup_{s \in Y} b_s$ is disconnected. Next assume that $\#((\partial{A_j(\c)})^e)  = 2$. The same arguments as in the proof of Theorem~\ref{BAS.thm20} show that $J$ is even and $(\partial{A_j(\c)})^e= \{s, s + J/2\}$ for some $s \in \BbZ_J$.
\enddemo

If $G \subseteq D_q$, then each edge of $Q_w$ can be labelled as $0$ or $1$. More precisely, since $\BbZ_J^i$ is $G$-invariant for any $i \in \BbZ_2$, Lemma~\ref{BCM.lemma10} implies the following lemma.

\lemma\label{MIS.lemma40}
For any $n \ge 1$, 
\[
E_n^{\ell} = E_n^{\ell}(\BbZ_J^0) \cup E_n^{\ell}(\BbZ_J^1)\quad\text{and}\quad E_n^{\ell}(\BbZ_J^0) \cap E_n^{\ell}(\BbZ_J^1) = \emptyset
\]
Moreover, define $\C_{m, \partial}^i = \C_{m, \partial}^{\ell}(\BbZ_J^i)$ for $i \in \BbZ_2$. Then $g({\c}) \in \C_{m, \partial}^i$ for any $g \in G$, $\c \in \C_{m, \partial}^i$ and $i \in \BbZ_2$.
\endlemma

\demo
The first part of the claims is immediate from Lemma~\ref{BCM.lemma10}-(1). The second part is due to Lemma~\ref{BCM.lemma10}-(4).
\enddemo

\lemma\label{MIS.lemma44}
Let $J \ge 6$ and let $M \ge M_J$. Assume that $Rot_q \subseteq G$. Let $\c \in \C_{M, m}(w)$ for some $w \in T$ and $m \ge j$. If $(\partial{A_j(\c)})^e \subseteq \BbZ_J^i$ for some $i \in \BbZ_2$, then $g^l(A_j(\c)) \cap g^{l + 1}(A_j(\c)) \neq \emptyset$ for any $l \ge 0$, where $g = \Theta_{2\pi/q}$, 
\begin{equation}\label{PTE.eq100}
\big(\partial(Rot_q(A_j(\c)))\big)^e = \BbZ_J^i,
\end{equation}
and $Rot_q(A_j(\c))$ is $(T_{m - j}, E_{m - j}^{\ell})$-connected. In particular, if $J = 6$, then $(\partial{A_j(\c)})^e = \BbZ_J^i$.
\endlemma

\demo
Let $\rho=\Theta_{2\pi/q}$. First, assume that $J=6$. If $\#((\partial{A_j(\c)})^e)=2$, then Lemma~\ref{MIS.lemma30} shows that $(\partial{A_j(\c)})^e = \{s, s + 3\}$ for some $s \in \BbZ_J$. Note that $s$ and $s + 3$ do not belong to the same $\BbZ_J^i$. This contradicts the fact that $(\partial{A_j(\c)})^e \subseteq \BbZ_J^i$. Thus it follows that $\#((\partial{A_j(\c)})^e)\geq 3$. Since $(\partial{A_j(\c)})^e \subseteq \BbZ_6^i$ and $\#(\BbZ_6^i) = 3$, we have $(\partial{A_j(\c)})^e = \BbZ^i_J$. \par
\noindent{\bf Claim 1:}\,\,For any $g \in Rot_3$, $A_j(\c)$ and $g(A_j(\c))$ are alternated.\\
Proof of Claim 1:\,\,Choose $x_s \in K(A_j(\c)) \cap b_s$ and $y_s \in K(g(A_j(\c)))\cap b_s$ for each $s \in \BbZ_J^i$. If $\{x_s\}_{s \in \BbZ_6^i} \cap \{y_s\}_{s \in \BbZ_6^i} \neq \emptyset$, then $A_j(\c) \cap g(A_j(\c))$ is alternated. Assume otherwise. Then, the set $\partial{Q_*}\backslash\{x_s| s \in \BbZ_6^i\}$ is decomposed into three connected components. Let $X_s$ be the connected component whose boundaries are $x_{s - 2}$ and $x_{s + 2}$ for each $s \in \BbZ_J^i$. Note that for any $s \in \BbZ_J$, $X_s \cap b_s = \emptyset$ and hence  $y_s \notin X_s$. Fix $t \in \BbZ_J^i$. Suppose that $y_t \in X_s$  for some $s \in \BbZ_J^i$. Then Two points $y_t$ and $y_{s}$ belong to different connected components of $\partial{Q_*}\backslash\{x_{s - 2}, x_{s + 2}\}$. Hence $A_j(\c)$ and $g(A_j(\c))$ are alternated. Thus, we have shown Claim 1.\par
Now Lemma~\ref{MIS.lemma20} implies that both $A_j(\c)$ and $g(A_j(\c))$ are $(T_{m - j}, E_{m - j}^{\ell})$-connected. Combining this fact with Claim 1 and using Lemma~\ref{BAS.lemma400}-(1), we see that $A_j(\c)) \cap g(A_j(\c)) \neq \emptyset$. Since $g$ is an arbitrary element of $Rot_3$, it follows that $Rot_3(A_j(\c))$ is $(T_{m - j}, E_{m - j}^{\ell})$-connected.\par 
Second, assume that $J \ge 8$.\\
\noindent {\bf Claim 2}:\,\,There exist $k, l \in (\partial{A_j(\c)})^e$ such that $\d(k, l) \ge 4$.\\
\noindent Proof of Claim 2: Let $k_0 \in (\partial{A_j(\c)})^e$. If $\{k_0 - 2, k_0, k_0 + 2\}\supseteq (\partial{A_j(\c)})^e$, then Lemma~\ref{MIS.lemma30} shows that $(\partial(A_j(\c))^e = \{k_0 - 2, k_0, k_0 + 2\}$. Since $\d(k_0 + 2, k_0 - 2) = 4$, we have verified the claim in this case. Otherwise, there exists $k_1 \in (\partial{A_j(\c)})^e$ such that $\d(k_0, k_1) \ge 3$. Since both $k_0$ and $k_1$ belongs to $\BbZ_J^i$, we see that $\d(k_0, k_1) \ge 4$. Thus, we completed a proof of Claim 2.\par
Again Lemma \ref{MIS.lemma20} implies $A_j(\c)$ is $(T_{m - j}, E_{m - j}^{\ell})$-connected. Hence by Lemma~\ref{BAS.lemma400}-(2) with $q_1 = q$ and $q_2 = 2$, we see that $g^l(A_j(\c)) \cap g^{l + 1}(A_j(\c) \neq \emptyset$ for any $l \ge 0$ and $Rot_q(A_j(\c))$ is $(T_{m - j}, E_{m - j}^{\ell})$-connected. Since $\BbZ_J^i$ is $Rot_q$-transitive, we have \eqref{PTE.eq100}.
\enddemo

\demo[Proof of Theorem~\ref{MIS.thm20}-(1)]
Let $J \ge 6$ and let $M \ge M_J$. Since $(\BbZ_J)^e = \BbZ_J^i$, it follows that $(\partial{A_0(\c)})^e \subseteq \BbZ_J^i$ for any $w \in T$, $m \ge 1$ and $\c \in \C_{M, m}(w)$. Let
\[
A = \bigcup_{g \in Rot_q} g(A_0(\c)).
\]
Then Lemma~\ref{MIS.lemma44} shows that $(\partial{A})^e = \BbZ_J^i = (\BbZ_J)^e$ and $K(A)$ is connected. Since $A \subseteq \H_*(\c)$, Theorem~\ref{BAS.thm10} suffices.
\enddemo

Now we proceed to a proof of Theorem~\ref{MIS.thm20}-(2). By Theorem~\ref{MIS.thm20}-(1), we see that $(\BbZ_J)^e = \BbZ_J$ under the assumption of Theorem~\ref{MIS.thm20}-(2).

\definition\label{PTE.def40}
For $A \subseteq T$, define $Q(A) = \bigcup_{w \in A} Q_w$ and define $Q_0(A)$ as the interior of $Q(A)$ as a subset of $(Q_*, d_{Q_*})$, where $d_{Q_*}$ is the restriction of the Euclidean metric onto $Q_*$. Also define $Q^{(m)} = \bigcup_{w \in T_m} Q_w$.
\enddefinition

\lemma\label{MIS.lemma45}
Assume that $J = 2q$ for some $q \ge 3$ and that $G \subseteq D_q$. Let $w \in T$, $m \ge 1$ and $k\in \{0, 1, \ldots, m - 1\}$. Moreover let $\c \in \C_{M, m}^{\ell}(w)$. Set $n = |w|$. Assume that $(\partial{A_k(\c)})^e \subseteq \BbZ_J^i$ for some $i \in \BbZ_2$. Let $k = j_1 + j_2$. Then 
 \begin{equation}\label{MIS.eq520}
[g(\s_{n + j_1}(\c)_l)]_{j_2} \in \C_{j_2, \partial}^i
\end{equation}
for any $g \in G$ and $l = 1, \ldots, l_{n + j_1}(\c)$. Moreover, $[A_{j_1}(\c)]_{j_2}$ is a path of $(T_{j_2}, E_{j_2}^{\ell}(\BbZ_J^i))$ and $(\partial_{\ell}A_{j_1}(\c))^e \subseteq \partial_{\BbZ_J^i}[A_{j_1}(\c)]_{j_2}$.
\endlemma

\demo
Define $q_0 = n + m, q = n + k, q_1 = n + j_1$ and $q_2 = j_2$. Let $\c = (\c(1), \ldots, \c(l_*))$. Since $\c \in \C_{M, m}(w)$, there exist $\c(0) \in S^m(w)$ and $\c(l_* + 1) \in S^m((\GG_M(w))^c)$ such that $(\c(0), \c(1))$ and $(\c(l_*), \c(l_* + 1))$ belong to $E_{n + m}^{\ell}$. Define $\widetilde{\c} = (\c(0), \c(1), \ldots, \c(l_*), \c(l_* + 1))$. Note that since $G \subseteq D_q$, $\BbZ_J^i$ is $G$-invariant. Hence, letting $X = \BbZ_J^i$, we have the desired statements immediately by Lemma~\ref{BCM.lemma40}.
\enddemo

\lemma\label{MIS.lemma50}
Assume that $J = 2q$ for some $q \ge 3$ and $M \ge M_J$ and that $G = Rot_q$ or $D_q$. Let $\c \in \C_{M, m}(w)$ for some $w \in T$ and $m \ge 1$ and let $j \in \{0, 1, \ldots, m\}$.\\
{\rm (1)}\,\,If $A_j(\c)$ dose not connect $0$ and $1$, then $(\partial{A_j(\c)})^e \subseteq \BbZ_J^i$ for some $i\in \BbZ_2$.\\
{\rm (2)}\,\,Suppose that both $A_j(\c)$ and $A_{j + 1}(\c)$ do not connect $0$ and $1$ and that $(\partial{A_j(\c)})^e \subseteq \BbZ_J^{i_1}$ and $(\partial{A_{j + 1}(\c)})^e \subseteq \BbZ_J^{i_2}$ for $i_1, i_2 \in \BbZ_2$. Then there exists $X \in F_{\partial}(\BbZ^{i_2 + 1})$ such that $X \subseteq \BbZ_J^{i_1 + 1}$. Furthermore, suppose that $A$ does not connect $0$ and $1$ whenever $A \subseteq \H_{|w| + j + 1, 1, m}(\c)$ and $A$ is $(T_m, E_m^{\ell})$-connected. Then $X = \BbZ^{i_1 + 1}$.
\endlemma

\demo
\noindent(1)\,\,Suppose that for each $i\in \BbZ_2$, there exists $s_i \in (\partial A_j(\c))^e \cap \BbZ_J^i$. Since $A_j(\c)$ does not connect $0$ and $1$, it follows that $d(s_0, s_1) = 1$. First assume that $\#(\partial{(A_j(\c)})^e) = 2$. Lemma~\ref{MIS.lemma30} implies $\d(s_0, s_1) = J/2$, so that we have a contradiction. Next assume that $\#((\partial A_j(\c))^e) \geq 3$. From Lemma~\ref{MIS.lemma30}, there exists $s \in (\partial{A_j(\c)})^e$ such that $b_s$ and $b_{s_1} \cup b_{s_2}$ are disjoint. Suppose that $s \in \BbZ_J^{i_*}$. Then $\d(s, s_{i_* + 1}) \ge 3$ and hence $A_j(\c)$ connects $0$ and $1$. Thus, we have obtained a contradiction.\\
\noindent(2)\,\,Since $\partial{A_{j + 1}(\c)} \subseteq \BbZ_J^{i_2}$, applying Lemma~\ref{MIS.lemma45} with $k = j + 1$, $j_1 = j$ and $j_2 = 1$, we see that $[A_j(\c)]_1$ is $(T_1, E_1^{\ell}(\BbZ_J^{i_2}))$-connected and 
\begin{equation}\label{MIS.eq30}
(\partial_{\ell}A_j(\c))^e \subseteq \partial_{\BbZ_J^{i_2}}[A_j(\c)]_1.
\end{equation}
Moreover, since $(\partial(A_j(\c)))^e \subseteq \BbZ_J^{i_1}$, Lemma~\ref{MIS.lemma44} implies that 
\[
g^l(A_j(\c)) \cap g^{l + 1}(A_j(\c)) \neq \emptyset
\]
for any $l \ge 0$, where $g = \Theta_{2\pi/q}$, that   
\begin{equation}\label{MIS.eq40}
\big(\partial(Rot_q(A_j(\c))\big)^e = \BbZ_J^{i_1},
\end{equation}
and that $Rot_q(A_j(\c))$ is $(T_{m - j}, E_{m - j}^{\ell})$-connected. By Lemma~\ref{BCM.lemma10}-(4), we see that $g^l([A_j(\c)_1])$ is $(T_1, E_1^{\ell}(\BbZ_J^{i_2})$-connected. Moreover, since $g^l([A_j(\c)]_1]) = [g^l(A_j(\c))]_1$, we see that $g^l([A_j(\c)]_1) \cap g^{l + 1}([A_j(\c)]_1) \neq \emptyset$. Thus let $A_* = Rot_q([A_j(\c)]_1)$. Then $A_*$ is $(T_1, E_1^{\ell}(\BbZ_J^{i_2}))$-connected. By \eqref{MIS.eq30}, it follows that
\[
(\partial_{\ell}Rot_q(A_j(\c)))^e \subseteq \partial_{\BbZ_J^{i_2}}(A_*)
\]
Combining this with \eqref{MIS.eq40}, we see that $\BbZ_J^{i_1} \subseteq \partial_{\BbZ_J^{i_2}}(A_*)$. Let $\wA \in \Con_T(\BbZ_J^{i_2 + 1})$ including $A_*$ and let $X = (\partial_{\BbZ_J^{i_2}}\wA)^c$. Then $X \in F_{\partial}(\BbZ_J^{i_2 + 1})$ and $X \subseteq \BbZ_J^{i_1 + 1}$.\par
Finally, suppose that $A$ does not connect $0$ and $1$ whenever $A \subseteq \H_{|w| + j + 1, 1, m}$ and $A$ is $(T_m, E_m^{\ell})$ connected and that $X \neq \BbZ_J^{i_1 + 1}$. Then there exist $i_1, i_2 \in \partial_{\BbZ_J^{i_2}}(\wA)$ such that $i_1 \in \BbZ_J^{i_1}$, $i_2 \in \BbZ_J^{i_1 + 1}$ and $\d(i_1, i_2) \ge 3$. Since $\wA \in \Con_T(\BbZ^{i_2 + 1})$, there exist a path $(w(1), \ldots, w(i_*))$ of $(T_1, E_1^{\ell}(\BbZ_J^{i_2}))$ and $k_1, k_2 \in \BbZ^{i_2}$ such that $w(1), \ldots, w(i_*) \in \wA$, $b_{w(1)}(k_1) \subseteq b_{i_1}$ and $b_{w(i_*)}(k_2) \subseteq b_{i_2}$. On the other hand, replacing $j$ and $i$ with $j + 1$ and $i_2$ respectively in Lemma~\ref{MIS.lemma44}, we see that 
\begin{equation}\label{MIS.eq60}
\partial{Rot_q(A_{j + 1}(\c))} = \BbZ_J^{i_2}
\end{equation}
and $Rot_q(A_{j + 1}(\c))$ is $(T_{m - j - 1}, E_{m - j - 1}^{\ell})$-connected. Define $A^{(k)}$ inductively as
\[
A^{(1)} = Rot_q(A_{j + 1}(\c))\quad\text{and}\quad A^{(k + 1)} = g_{w(k), w(k + 1)}(A^{(k)})
\]
for $k = 1, \ldots, i_* - 1$. Since $(w(k), w(k + 1)) \in E_1^{\ell}(\BbZ^{i_2})$ and \eqref{MIS.eq60}, using Lemma~\ref{BAS.lemma201}, we see that $w(k)A^{(k)} \cup w(k + 1)A^{(k + 1)}$ is $(T_{m - j}, E_{m - j}^{\ell})$-connected. Moreover, $i_1 \in \partial{w(1)A^{(1)}}$ and $i_2 \in \partial{w(i_*)A^{(i_*)}}$. Let $A = S^j(\cup_{l = 1}^{i_*} w(k)A^{(k)})$. Then $A \subseteq \H_{|w| + j + 1, 1, m}(\c)$, $A$ is $(T_m, E_m^{\ell})$-connected and $A$ connects $0$ and $1$. Thus we have a contradiction, so that $X = \BbZ_J^{i_1 + 1}$.
\enddemo

\demo[Proof of Theorem~\ref{MIS.thm20}-(2)]
Suppose that $(K,d)$ does not have $p$-conductively homogenous for some $p > \textrm{dim}_{AR}(K,d)$. Then, by Theorem~\ref{PTE.thm10}, there exist $w \in T$ and $\c \in \C_{M, m}(w)$ such that if 
\[
A \subseteq \bigcup_{j = 0, 1, 2, n_2 = 0, 1}\H_{|w| + j, n_2, m}(\c)
\] 
and $A$ is $(T_m, E_m^{\ell})$-connected, then $A$ does not connect $0$ and $1$. Since $S^k(A_j(\c))$ belongs to $\H_{|w| + j, 0, m}(\c)$ and is $(T_m, E_m^{\ell})$-connected for $j = 0, 1, 2$, it follows that $S^j(A_j(\c))$ does not connect $0$ and $1$. Hence $A_j(\c)$ does not connect $0$ and $1$ for $j = 0, 1, 2$. By Lemma~\ref{MIS.lemma50}-(1), for $j = 0, 1, 2$, there exist $i_j \in \BbZ_2$ such that $\partial{A_j(\c)} \subseteq \BbZ_J^{i_j}$. Then using Lemma~\ref{MIS.lemma50}-(2), we see that
\[
\BbZ_J^{i_0 + 1} \in F_{\partial}(\BbZ^{i_1 + 1})\quad\text{and}\quad \BbZ_J^{i_1 + 1} \in F_{\partial}(\BbZ^{i_2 + 1}).
\]
If $i_0=i_1$ or $i_1=i_2$, then \eqref{MIS.eq500} follows. Otherwise, we have \eqref{MIS.eq510}.\par
Finally assume that $J = 6$ or $8$ and that \eqref{MIS.eq510} holds. Then for each $i \in \BbZ_2$, there exists $A_i \in \Con_T(\BbZ_J^{i})$ such that $\partial_{\BbZ_J^{i + 1}}{A_i} = \BbZ_J^i$. Then $A_0$ and $A_1$ are alternated. So, for each $i \in \BbZ_2$, there exists a simple path $\c_i = (\c_i(1), \ldots, \c_i(l_i))$ of $(T_1, E_1(\BbZ_J^{i + 1}))$ such that $\{c_i(1), \ldots, \c_i(l_i)\} \subseteq A_i$ for each $i \in \BbZ_2$ and $\c_0$ and $\c_1$ are alternated. By Lemma~\ref{BAS.lemma400}-(1), 
$\c_0 \cap \c_1 \neq \emptyset$. Let $s \in \c_0 \cap \c_1$. Define
\begin{multline*}
B_i = \{k| k \in \BbZ_J^{i + 1}, \text{there exists $t \in T_1$ and $l \in \BbZ_J^{i + 1}$ such that $b_k(s) = b_l(t)$}\\
\text{ or there exists $l \in \BbZ_J^{i}$ such that $b_k(s) \subseteq b_l$.}\}
\end{multline*}
for any $i \in \BbZ_2$. Then $\#(B_i) \ge 2$ for any $i \in \BbZ_2$. On the other hand, Lemma~\ref{MIS.lemma05} shows that if $k_1, k_2 \in B_0 \cup B_1$, then $\d(k_1, k_2) \ge 2$. In the case $J = 6$, let $k_1 \neq k_2 \in B_0$. Then $\d(k_1, k_2) = 2$ and there is only one possible location of $B_1$. So we have a contradiction. In the case $J = 8$, let $k_1 \neq k_2 \in B_0$. Then $k_2 - k_2 = 2$ or $4$. In either case, it is impossible to accommodate $B_1$. Thus we never have \eqref{MIS.eq510} if $J = 6$ or $8$.
\enddemo

\setcounter{equation}{0}
\section{Conductive homogeneity of $(J, G)$-s.s.\,system III: no global symmetry }\label{WNG}

In this section, we deal with the cases where $G = \{I\}$, where $I$ is the identity map.  As in the previous sections, $(S, \{f_s\}_{s \in S}, G)$ is a $(J, G)$-self-similar system, where $J \ge 3$ and $G$ is a subgroup of $D_J$.\par
 Note that if the set of global symmetry $G_*$ is trivial, i.e. $G^* = \{I\}$, then $G$ is trivial but the other direction is not true. In fact, if $J = 3$ (resp. $4$), since $Q^{(1)}$ is a part of the triangle tiling of $Q_*$ (reps. a square tiling of $Q_*$), there exists a folding map $\vp : Q_* \to Q_*$ that is continuous and a similitude on every $Q_s$ with $\vp(Q_s) = Q_*$. See Definition~\ref{WNG.def10} for a generalized version of its definition.  Consequently using the folding map, one can modify the original system of self-similar sets $\{f_s\}_{s \in S}$ to make $G = \{I\}$. See \cite[Theorem~4.22]{Ki22} for details for the case $J = 4$. Moreover, by \cite{Sasaya1}, there always exists a folding map for $J = 5$. So, we may assume that $G = \{I\}$ without loss of generality if $J \in \{3, 4, 5\}$.
 
 \definition[Folding map]\label{WNG.def10}
A continuous map $\psi: K \to K$ is called a folding map for a $(J, G)$-self-similar system $(S, \{f_s\}_{s \in S}, G)$ if, for each $s \in S$, there exists an affine bijective map $\psi_s: Q_s \to Q_*$ such that $\psi|_{K_s} = \psi_s$ and $\psi(K_s) = K$.
\enddefinition

The next proposition tells that the existence of a folding map is virtually equivalent to the condition $G = \{I\}$.

\prop\label{WNG.prop10}
Let $(S, \{f_s\}_{s \in S}, G)$ is a $(J, G)$-self-similar system and let $K$ be the associated self-similar set. Assume that there exists no isolated contact point of cells.\\
{\rm (1)}\,\,$G = \{I\}$ if and only if there exists a folding map $\psi : K \to K$ such that 
\[
(\psi|_{K_s})^{-1} = f_s|_K
\]
for any $s \in S$. This folding map $\psi$ is called the folding map associated with the $(J, \{I\})$-self-similar system $(S, \{f_s\}_{s \in S}, \{I\})$.\\
{\rm (2)}\,\, Suppose that there exists a folding map $\psi: K \to K$.  For $s \in S$, define an affine similitude
\[
F_s = (\psi_s)^{-1},
\]
where $\psi_s:Q_s \to Q_*$ is the affine map appearing in Definition~\ref{WNG.def10}. Then $(S, \{F_s\}_{s \in S}, \{I\})$ is a $(J, \{I\})$-self-similar system and its associated self-similar set is $K$.
\endprop

To prove this proposition, we need the following lemmas.

\lemma\label{WNG.lemma10}
Suppose that $G  = \{I\}$.\\
{\rm (1)}\,\,
Let $w, v \in T_n$. If there exist $i, j \in \BbZ_J$ such that $Q_w \cap Q_v =b_i(w) = b_j(v)$, then $i = j$ and $f_w|_{b_i} = f_v|_{b_i}$.\\
{\rm (2)}\,\,Let $n \ge 1$ and let $w, v \in T_n$. Assume that there exists $i \in \BbZ_J$ such that $Q_w \cap Q_v = b_i(w) = b_i(v)$. Then $R_{w, v}^*(wu) = vu$ for any $u \in T$. In particular,  if $b_j(wu) \subseteq b_i(w)$ for some $j \in \BbZ_J$ and $u \in T_m$, then $Q_{wu} \cap Q_{vu} = b_j(wu) = b_j(vu)$.\\
{\rm(3)}\,\,Let $w, v \in T_{n + m}$. If $[w]_n \neq [v]_n$ and there exists $i, j \in \BbZ_J$ such that $Q_w \cap Q_v = f_w(b_i) = f_v(b_j)$, then $\s_n(w) = \s_n(v)$. \par
\endlemma

\demo
(1)\,\,Since $g_{w, v} \in G = \{I\}$, we see that
\[
R_{w, v}{\circ}f_w  = f_v.
\]
So, $f_w(x) = R_{w, v}(f_w(x)) = f_v(x)$ for any $x \in b_i$. Hence $f_v(b_i) = f_w(b_i) = f_v(b_j)$. This yields $b_i = b_j$ and hence $i = j$.\\
(2)\,\,If $b_j(wu) \subseteq b_i(w)$, then
\[
b_j(wu) = R_{w, v}(b_j(wu)) = R_{w, v}{\circ}f_{w}{\circ}f_u(b_j) = f_v{\circ}f_u(b_j) = b_j(vu).
\]
Thus we have $Q_{wu} \cap Q_{vu} = b_j(wu) = b_j(vu)$.\\
(3)\,\,This immediately follows from (3).
\enddemo

\lemma\label{WNG.lemma02}
Let $(S, \{f_s\}_{s \in S}, G)$ be a $(J, G)$-self-similar system and let $K$ be the associated self-similar set. Then there exist $x, y \in K$ such that $x$ and $y$ are independent.
\endlemma

\demo
Suppose otherwise. Then since $K$ is connected and $K \cap b_i \neq \emptyset$ for any $i \in \BbZ_J$, $K$ must be a line segment $\overline{ab}$ where $K \cap \partial{Q_*} = \{a, b\}$ and $a + b = 0$. Again by the fact that $K \cap b_i \neq \emptyset$ for any $i \in S$, this is possible only when $J = 3$ or $4$. In the case $J = 4$, $K$ must be $\overline{p_0p_2}$ or $\overline{p_1p_3}$. Since $K_s$ is a part of the line segment, it follows that $Q_s \cap Q_t$ is empty or a single point for any $s \neq t \in S$. This contradicts (A5). In the case $J = 3$, $K$ must be a line segment between one of the vertex $p_i$ and the midpoint of $p_{i + 1}$ and $p_{i + 2}$. Let $q$ be the mid point. Choose $s \in S$ such that $q \in K_s$ and choose $t \in S$ such that $Q_s \cap Q_t \neq \emptyset$. Then $Q_s \cap Q_t$ is a single point and hence we have a contradiction to (A5). Thus we have verified the desired statement.
\enddemo

\lemma\label{WNG.lemma03}
For $j = 1, 2$, let $f_j:\BbR^2 \to \BbR^2$ be defined as
\[
f_j(x) = r\vp_j(x) + c_j
\]
for any $x \in \BbR^2$, where $r > 0$, $\vp_j \in O(2)$ and $c_j \in \BbR^2$. Assume that there exists $i \in \BbZ_J$ such that $f_1|_{b_i} = f_2|_{b_i}$ and $f_1(Q_*) \cap f_2(Q_*) = f_1(b_i) = f_2(b_i)$. Then $(\vp_2)^{-1}{\circ}\vp_1 = R_{\rho(i)}$.
\endlemma

\demo
Set $q_0 = (\cos{\rho(i)}, \sin{\rho(i)})$ and let $q_1$ be the midpoint of $b_i$. Since $f_1|_{b_i} = f_2|_{b_i}$, it follows that
\begin{equation}\label{WNG.eq00}
r\vp_1(q_1 + {\a}q_0) = r\vp_2(q_1 + {\a}q_0)
\end{equation}
for any $\a \in \BbR$. Letting $\a = 0$, we obtain $r\vp_1(q_1) + c_1 = r\vp_2(q_1) + c_2$. Since $c_1 - c_2 = -2r\vp_1(q_1)$, it follows that $\vp_1(q_1) = -\vp_2(q_1)$ and hence $(\vp_2)^{-1}{\circ}\vp_1(q_1) = -q_1$. Moreover, using \eqref{WNG.eq00}, we see that $\vp_1(q_0) = \vp_2(q_0)$. Therefore, it follows that $(\vp_2)^{-1}{\circ}\vp_1(q_0) = q_0$. Consequently, $(\vp_2)^{-1}{\circ}\vp_1 = R_{\rho(i)}$.
\enddemo

\demo[Proof of Proposition~\ref{WNG.prop10}]
(1)\,\,Assume that $G = \{I\}$. Define $\psi(x) = (f_s)^{-1}(x)$ for $s \in S$ and $x \in K_s$. What should be verified is the consistency of this definition, i.e. to show $(f_s)^{-1}(x) = (f_t)^{-1}(x)$ if $x \in K_s \cap K_t$ for some $s, t \in S$. 
Suppose that $x \in K_s \cap K_t$ for some $s, t \in S$. If (A4)-(a) holds, then Lemma~\ref{WNG.lemma10}-(1) shows the consistency.
Assume (A4)-(b) holds, i.e. $x = f_s(p_i) = f_t(p_j)$ for some $i, j \in \BbZ_J$. Since $x$ is not an isolated contact point of cells, Lemma~\ref{COP.lemma10}-(2) implies that $(\GG^1(x), E^{\ell}_1|_{\GG^1(x)})$ is connected. Hence there exist $\{s_1, \ldots, s_k\} \subseteq \GG^1(x)$ such that $s = s_1, t = s_k$ and $(s_i, s_{i + 1}) \in E^{\ell}_1$ for any $i = 1, \ldots, k - 1$. For each $i \in \{1, \ldots, k -1\}$, the above argument assuming (A4)-(a) yields that $(f_{s_i})^{-1}(x) = (f_{s_{i + 1}})^{-1}(x)$. Consequently we obtain $(f_s)^{-1}(x) = (f_t)^{-1}(x)$. Thus, the map $\psi: K \to K$ is well-defined and obviously it is a holding map by definition.\par
Conversely, assume that there exists a folding map $\psi$ satisfying $(\psi|_{K_s})^{-1} = f_s|_K$ for any $s \in S$. \\
{\bf Claim 1}\,\,$(\psi_s)^{-1} = f_s$ on $Q_*$ for any $s \in S$.\\
Proof of Claim 1:\,\,By the definition of the folding map, we see that $(\psi_s)^{-1} = f_s$ on $K$ and $(\psi_s)^{-1}(0) = f_s(0) = c_s$. By Lemma~\ref{WNG.lemma02}, there exist $x, y \in K$ such that $x$ and $y$ are independent. Since $(\psi_s)^{-1}(x) = f_s(x)$, we have the claim. \qed\\
Now the condition (A3) is obvious by letting $e_*: S \to S$ be the identity map. To verify (A4)-(a), suppose that $Q_s \cap Q_t = f_s(b_i) = f_t(b_j)$ for some $s, t \in S$ and $i, j \in \BbZ_J$. The definition of the holding map immediately implies that $i = j$ and $f_s|_{b_i} = f_t|_{b_i}$. Using Lemma~\ref{WNG.lemma03}, we see that $(\vp_t)^{-1}{\circ}\vp_s = R_{\rho(i)}$. This shows (A4)-(a) with $G = \{I\}$.\\
(2)\,\,
First we are going to show that $(S, \{F_s\}_{s \in S})$ is a $(J, \{I\})$-self-similar system. Since $\psi_s(Q_s) = Q_*$, we see that $F_s(Q_*) = Q_s$. Hence there exists $\xi_s \in D_J$ such that
\[
F_s(x) = r\xi_s(x) + c_s
\]
for any $x \in Q_*$. This shows (A1). Assume that $f_s(b_j) \subseteq b_i$. Since $f_s(b_j) = F_s(b_{j'})$ for some $j' \in \BbZ_J$, the condition (A2) is verified. Letting $e_*$ be the identity map of $S$, we have (A3). Next suppose that $Q_s \cap Q_t = f_s(b_i) = f_t(b_j)$. Since $\psi_s|_{Q_s \cap Q_t} = \psi_t|_{Q_s \cap Q_t}$, it follows that $Q_s \cap Q_t = F_s(b_k) = F_t(b_k)$ for some $k \in \BbZ_J$. Applying Lemma~\ref{WNG.lemma03}, we have $(\xi_t)^{-1}{\circ}\xi_s = R_{\rho(i)}$. This immediately yields (A4). Finally, since $E_1^l$ stays the same even if we replace $\{f_s\}_{s \in S}$ with $\{F_s\}_{s \in S}$, (A5) follows. Thus we have shown that $(S, \{F_s\}_{s \in S})$ is a $(J, \{I\})$-self-similar system. Since $F_s(K) = K_s$, $K$ is the associated self-similar set.
\enddemo

\lemma\label{WNG.lemma15}
Suppose that $G = \{I\}$ and let $\psi$ be the associated folding map. For $m \ge 2$ and $w = \word wm \in T_m$, define $\psi_{w}: Q_{w} \to Q_*$ by $\psi_{w} = \psi_{w_m}{\circ}\cdots{\circ}\psi_{w_1}$. Then $\psi^m|_{K_{\word wm}} = \psi_{\word wm}|_{K_{\word wm}}$ and, for any $w \in T$ and $n \in \BbN$ with $n \le |w|$, 
\[
\psi_{\word wn}{\circ}f_w = f_{\s_n(w)}.
\]
\endlemma

\demo
The equality $\psi^m|_{K_{\word wm}} = \psi_{\word wm}|_{K_{\word wm}}$ is straightforward from the definition of $\psi_s$. Since $\psi_{w} = (f_w)^{-1}$, we have $\psi_{\word wn}{\circ}f_w = f_{\s_n(w)}$.
\enddemo

Hereafter in this section, we always assume that $G = \{I\}$ and $\psi: K \to K$ is a folding map associated with $(S, \{f_s\}_{s \in S})$. Since $G = \{I\}$, every subset of $\BbZ_J$ is $G$-invariant. This fact enable us to show that $F_{\partial}(n, \cdot)$ defined in Definition~\ref{MIS.def200} is given by the $n$-th iteration of $F_{\partial}(1, \cdot)$ as follows.

\thm\label{WNG.thm10}
Assume that $G = \{I\}$. For any $n, m \ge 1$ and $X \subseteq \BbZ_J$, 
\[
F_{\partial}(n + m, X) = \bigcup_{Y \in F_{\partial}(m, X)} F_{\partial}(n, Y).
\]
\endthm

A proof of this theorem is given in the next section.\par
Let $\P(\BbZ_J)$ be the collection of subsets of $\BbZ_J$. Define $F_{\partial}: \P(\BbZ_J) \to \P(\BbZ_J)$ by 
\[
F_{\partial}(\A) = \bigcup_{X \in \A} F_{\partial}(1, X).
\]
for $\A \in \P(\BbZ_J)$. Then the above theorem shows that if $F^n_{\partial}$ is the $n$-th iteration of $F_{\partial}$, then
\[
F^n_{\partial}(\A) = \bigcup_{X \in \A} F_{\partial}(n, X)
\]
for any $\A \in \P(\BbZ_J)$.

\definition\label{WNG.def40}
Define $\B_l = \{X| X \subseteq \BbZ_J, \#(X) = l\}$ for $l = 0, \ldots, J$ and 
\[
\B_{J - 2}^{os} = \begin{cases}\{\sd{\BbZ_J}{\{s, s + J/2\}}| s = 0, \ldots, J/2 - 1\} &\text{if $J$ is even,}\\
                                                     \emptyset & \text{if $J$ is odd.}
                            \end{cases}
\]
The symbol ``$os$''  in the definition of $\B_{J - 2}^{os}$ represents the word  ``opposite sides''. Moreover define
\[
\B^H = \bigcup_{l = 1, \ldots, J - 3} \B_l \cup \B_{J - 2}^{os}
\]
and
\[
B^L = (\sd{\B_{J - 2}}{\B_{J - 2}^{os}}) \cup \B_{J - 1} \cup \B_J.
\]
\enddefinition

Note that $\B_J = \{\BbZ_J\}$ and $\B_0 = \{\emptyset\}$ and that $F_{\partial}(\B_l) = \B_l$ for $l = 0$ and $J$. Furthermore, we show
\[
F_{\partial}(\B^L) \subseteq \B^L
\]
in the next section.\par
In terms of the iterations of $F_{\partial}$, we give a sufficient condition for the conductive homogeneity in the next theorem, which is the main result of this section.

\thm\label{WNG.thm20}
Assume that $G = \{I\}$ and that there exists no isolated contact point of cells. The following three conditions $(F_{\partial}1)$, $(F_{\partial}2)$ and $(F_{\partial}3)$ are equivalent:\\
$(F_{\partial}1)$\,\,For any $X \in \B^H$ and $n \ge 1$, $X \notin F_{\partial}^n(X)$,\\
$(F_{\partial}2)$\,\,For any $X \in \B^H$, there exists $n \ge 1$ such that either
\[
\emptyset \in F_{\partial}^n(X)\quad\text{or}\quad F_{\partial}^n(X) \subseteq \B^L.
\]
$(F_{\partial}3)$\,\,For any $X \in \B^H$, there exists $n \ge 1$ such that
\[
F_{\partial}^n(X) \subseteq \{\emptyset\} \cup \B^L.
\]
Moreover, if any of the above conditions holds, then $(K, d_*)$ is $p$-conductively homogeneous for any $p > \dim_{AR}(K, d_*)$.
\endthm

This theorem is proven in the next section.\par
To verify the condition of the above theorem, the next proposition is useful.

\prop\label{WNG.prop100}
Assume that $G = \{I\}$ and that $X \subseteq Y \subseteq \BbZ_J$.\\
{\rm (1)}\,\,If $\emptyset \in F_{\partial}^n(Y)$, then $\emptyset \in F_{\partial}^n(X)$.\\
{\rm (2)}\,\,If $F^n_{\partial}(X)$ is contained in $\B^L$, then so does $F^n_{\partial}(Y)$.
\endprop

This proposition is immediate from the next lemma.

\lemma\label{WNG.lemma500}
Assume that $G = \{I\}$ and that $X \subseteq Y \subseteq \BbZ_J$. Then, for any $Z \in F_{\partial}(n, X)$, there exists $Z' \in \F_{\partial}(n, Y)$ such that $Z \subseteq Z'$. Conversely, for any $Z' \in F_{\partial}(n, Y)$, there exists $Z \in \F_{\partial}(X)$ such that $Z \subseteq Z'$.
\endlemma

\demo
Let $A \in \Con_T(n, X)$. Then since $X^c \supseteq Y^c$, there exist $B_1, \ldots, B_k \in \Con_T(n, Y)$ such that
\[
A = \bigcup_{i = 1}^k B_i.
\]
Then $\partial_{X^c}A \supseteq \partial_{Y^c}B_i$ for any $i = 1, \ldots, k$. This implies the desired statement.
\enddemo

In the rest of this section, we present two examples to illustrate the usefulness of Theorem~\ref{WNG.thm20}.

\begin{figure}[ht]
\hspace*{0pt}
\begin{minipage}[b]{3cm}
\centering
\includegraphics[width = 100pt]{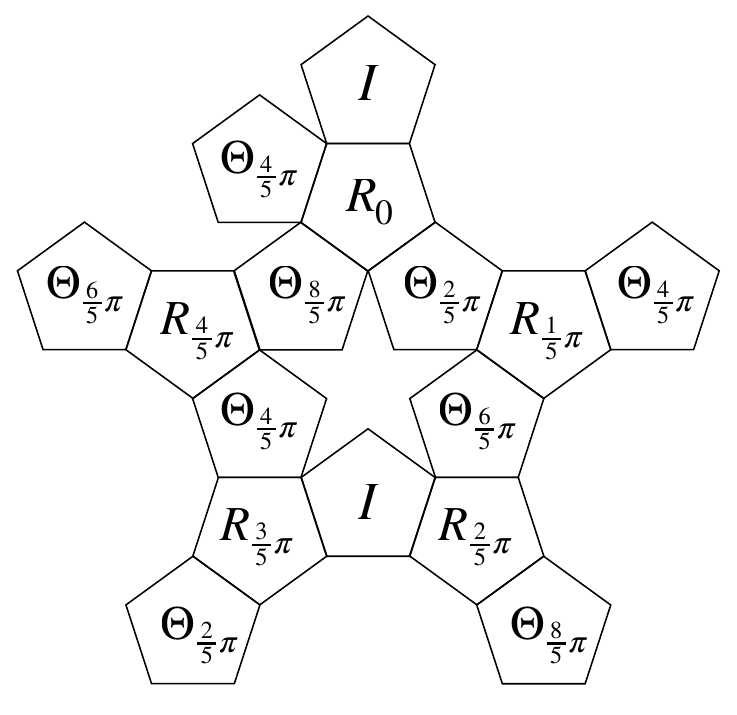}
\vspace{-20pt}
\end{minipage}
\hspace{10pt}
\begin{minipage}[b]{3cm}
\centering
\includegraphics[width = 100pt]{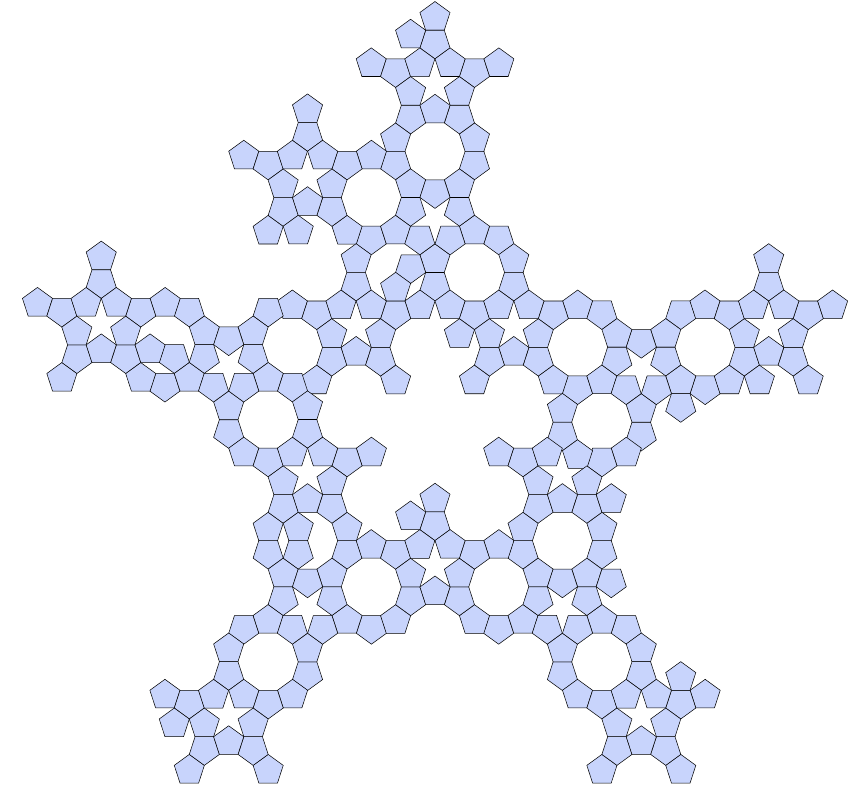}
\end{minipage}
\hspace{10pt}
\begin{minipage}[b]{3cm}
\centering
\includegraphics[width = 100pt]{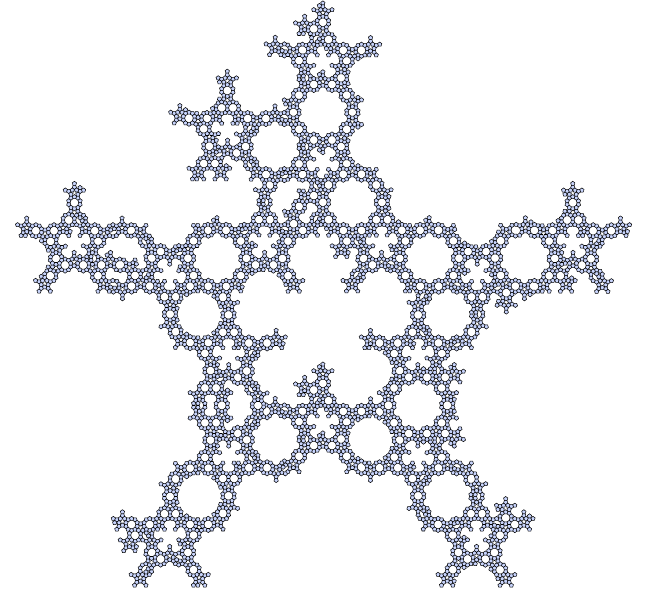}
\end{minipage}

\caption{$J = 5$, $G = \{I\}$}\label{NSpenta1}
\end{figure}

\begin{figure}[ht]
\hspace*{0pt}
\begin{center}
\includegraphics[width = 300pt]{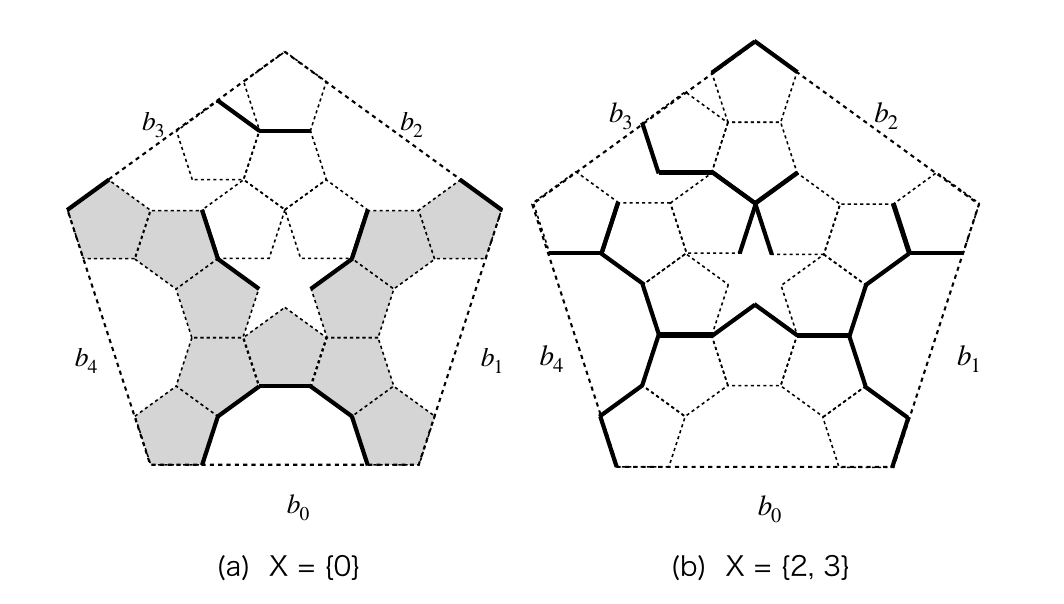}
\end{center}
\vspace{-20pt}
\caption{$\Con_T(1, X)$}\label{NSpentaConn}
\end{figure}

\example\label{WNG.ex10}
Let $S$ be the collection of pentagons in the most left figure of Figure~\ref{NSpenta1}. Define $\vp_s \in D_5^*$ as indicated in the corresponding pentagon $s \in S$ in the figure. Moreover, let $c_s$ be the center of the corresponding pentagon $s \in S$. Then we see that $(S, \{f_s\}_{s \in S}, G)$ is a pentagon-based self-similar system with $G = \{I\}$. Since $J = 5$, it follows that $\B^H = \B_1 \cup \B_2$. Connected components of $(T_1, E_1^{\ell}(X^c))$, which is $\Con_T(1, X)$, for $X = \{0\}$ is illustrated in Figure~\ref{NSpentaConn}-(a). Let $X = \{0\}$ for the moment. The thick lines are images of $b_0$ under the family of maps $\{f_s\}_{s \in S}$ and hence they are disconnected in $(T_1, E_1^{\ell}(X^c))$. The grey region in Figure~\ref{NSpentaConn}-(a) is one of the elements of $\Con_T(1, X)$. Denoting this connected component by $C$, we see that $\partial_{X^c}C = \{0, 1, 4\}$ and hence $\{0, 1, 4\}^c = \{2, 3\} \in F_{\partial}(1, X)$. In the same manner, we have
\[
F_{\partial}(1, X) = \{\{2, 3\}, \{3\}^c, \{0, 1, 4\}\}.
\]
To show the conductive homogeneity, we try to verify $(F_{\partial}2)$ of Theorem~\ref{WNG.thm20}. Since $\{3\}^c$ and $\{0, 1, 4\}$ belong to $\B^L$, it is enough to take care of $\{2, 3\}$ to show $(F_{\partial}2)$ for $X = \{0\}$. Actually, Figure~\ref{NSpentaConn}-(b) shows $\Con_{T}(1, \{2, 3\})$ and it follows that
\[
F_{\partial}(1, \{2, 3\}) = \{\{0\}^c, \{0, 1, 2\}, \{0, 3, 4\}, \{3\}^c\} \subseteq \B^L.
\]
Thus we have obtained $F_{\partial}(2, \{0\}) \subseteq \B^L$. Consequently Proposition~\ref{WNG.prop100}-(2) implies that $F_{\partial}(2, \{0, k\}) \subseteq \B^L$ for any $k \in \BbZ_5$. Entirely similar arguments yields $(F_{\partial}2)$ for other $X \in \B^H$ and hence the corresponding self-similar set $K$ is $p$-conductively homogeneous for any $p > \dim_{AR}(K, d_*)$.
\endexample

\begin{figure}[ht]
\hspace*{0pt}
\begin{minipage}[b]{3cm}
\centering
\includegraphics[width = 100pt]{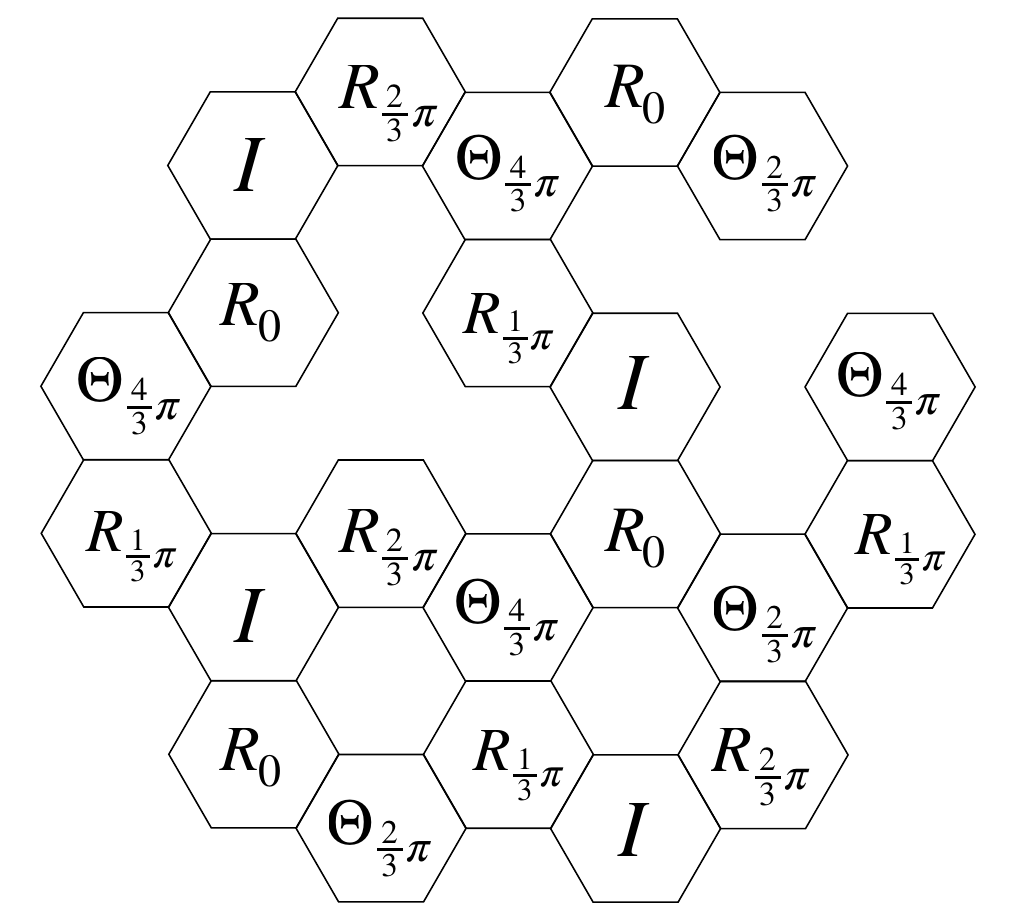}
\end{minipage}
\hspace{10pt}
\begin{minipage}[b]{3cm}
\centering
\includegraphics[width = 100pt]{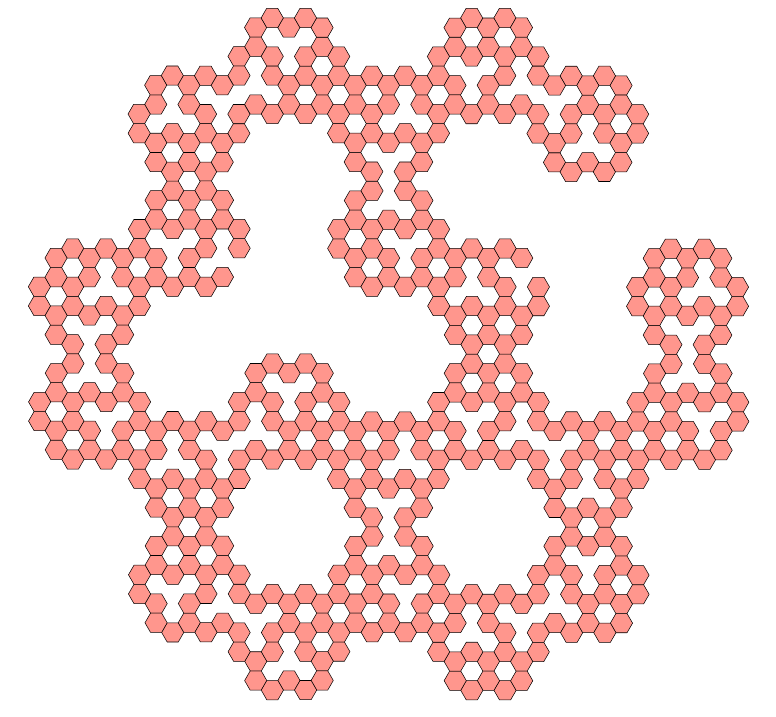}
\end{minipage}
\hspace{10pt}
\begin{minipage}[b]{3cm}
\centering
\includegraphics[width = 100pt]{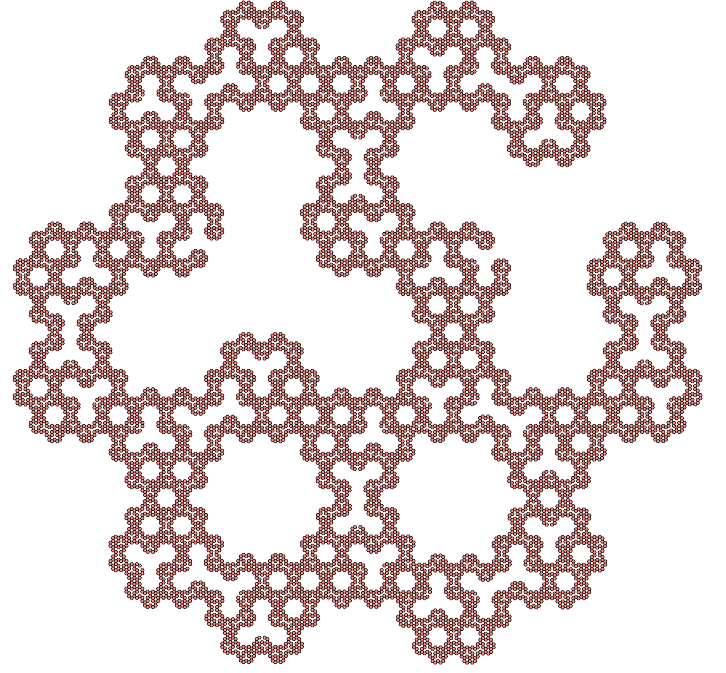}
\end{minipage}

\caption{$J = 6$, $G = \{I\}$}\label{NShexa1}
\end{figure}

\example\label{WNG.ex20}
Let $S$ be the collection of hexagons in the most left figure of Figure~\ref{NShexa1}. Define $\vp_s$ as indicated in the corresponding hexagon $s \in S$. Moreover, let $c_s$ be the centre of the corresponding hexagon $s \in S$. The $(S, \{f_s\}_{s \in S}, G)$ is a $(6, G)$-self-similar system with $G = \{I\}$. As the last example, we are going to show the condition $(F_{\partial}2)$ in Theorem~\ref{WNG.thm20}. First, examining the first step, i.e. the most left figure of Figure~\ref{NShexa1}, we see that
\begin{align*}
F_{\partial}(1, \{0\}) &= \{\BbZ_6, \{1, 2, 3, 4\}, \{0, 3, 4, 5\}, \{0, 1, 2, 5\}\} \subseteq \B^L\\
F_{\partial}(1, \{2\}) &= \{\{k\}^c| k \in \BbZ_6\} \cup \{\BbZ_6\} \subseteq \B^L\\
F_{\partial}(1, \{4\}) &= \{\BbZ_6, \{0, 1, 2, 3\}, \{2, 3, 4, 5\}, \{0, 1, 4, 5\}\} \subseteq \B^L.
\end{align*}
Hence by Proposition~\ref{WNG.prop100}-(2), if $X \cap \BbZ_6^0 \neq \emptyset$, then $F_{\partial}(1, X) \subseteq \B^L$. Consequently,  the remaining cases are when $X \subseteq \BbZ_6^1$. In such a case, we have the following two facts:\\
Fact 1:\,\,$(s_1, s_2) \in E_1^{\ell}$ if and only if there exists $i \in \BbZ_6^0$ such that $b_i(s_1) = b_i(s_2)$,\\
Fact 2:\,\,If $b_s(i) \subseteq b_j$ for some $i, j \in \BbZ_6$ and $s \in S$, then $i \in \BbZ_6^0$.\par
Suppose $X \subseteq \BbZ_6^1$. Then $\Con_T(1, X) = \{S\}$ by Fact 1 and Fact 2 implies $F_{\partial}(1, X) = \{\emptyset\}$. 
Thus we have verified the condition $(F_{\partial}2)$ and hence the corresponding self-similar set $K$ is $p$-conductively homogeneous for any $p > \dim_{AR}(K, d_*)$.
\endexample

\section{Proofs of Theorems~\ref{WNG.thm10} and \ref{WNG.thm20}}\label{POT}

This section is devoted to proofs of Theorems~\ref{WNG.thm10} and \ref{WNG.thm20}. Throughout this section, $(S, \{f_s\}_{s \in S}, G)$ is a $(J, G)$-s.s.\,system with $G = \{I\}$.

\lemma\label{WNG.lemma21}
Assume that $G = \{I\}$. Let $Y \subseteq \BbZ_J$. If $w, v \in T_{n + m}$ and $(w, v) \in E_{n + m}^{\ell}(Y)$, then 
\[
\begin{cases}
\s_n(w) = \s_n(v)\,\,&\text{if $[w]_n \neq [v]_n$,}\\
(\s_n(w), \s_n(v)) \in E_m^{\ell}\,\,&\text{otherwise.}
\end{cases}
\]
\endlemma

\demo
If $[w]_n \neq [v]_n$, then Lemma~\ref{WNG.lemma10}-(3) shows that $\s_n(w) = \s_n(v)$. Suppose $[w]_n = [v]_n$. Since $(w, v) \in E_{n + m}^{\ell}(Y)$, there exists $s \in Y$ such that $f_w(b_s) = f_v(b_s)$, and hence $f_{\s_n(w)}(b_s) = f_{\s_n(v)}(b_s)$, so that $(\s_n(w), \s_n(v)) \in E_n^{\ell}(Y)$.
\enddemo

\lemma\label{WNG.lemma40}
Assume that $G = \{I\}$. Let $n, m \ge 1$, $X \subseteq \BbZ_J$ and $A \in \Con_T(n + m, X)$. Then $\s_n(A) \in \Con_T(m, X)$. Moreover, $A = \cup_{w \in [A]_n} w\s_n(A)$ and $[A]_n \in \Con_T(n, (\partial_{X^c}(\s_n(A)))^c)$. Furthermore, 
\begin{equation}\label{WNG.eq10}
\partial_{\partial_{X^c}(\s_n(A))}[A]_n = \partial_{X^c}A
\end{equation}
\endlemma

\demo
First we show $\s_n(A) \in \Con_T(m, X)$. Let $w, v \in \s_n(A)$. Then there exists a path $(w(1), \ldots, w(k))$ of $(T_{n + m}, E_{n + m}^{\ell}(X^c))$ such that $\s_n(w(1)) = w$ and $\s_n(w(k)) = v$. Lemma~\ref{WNG.lemma21} implies that either $\s_n(w(i)) = \s_n(w(i + 1))$ or $(\s_n(w(i)), \s_n(w(i + 1))) \in E_m^{\ell}(X^c)$. This shows that $\s_n(A)$ is connected in $(T_m, E_m^{\ell}(X^c))$. Suppose that $(w, u) \in E_m^{\ell}(X^c)$ for some $w \in \s_n(A)$ and $u \in T_m$. Let $w' \in A$ satisfying $\s_n(w') = w$. Set $\tau = [w']_n$. Then $w' = \tau{w}$ and $(\tau{w}, \tau{u}) \in E_{n + m}^{\ell}(X^c)$. Since $A$ is a connected component of $(T_{n + m}, E_{n + m}^{\ell}(X^c))$, we see that $\tau{u} \in A$ and hence $u \in \s_n(A)$. Thus we have shown that $\s_n(A) \in \Con_T(n + m, X)$. \par
Second, let $\widetilde{A} = \bigcup_{w \in [A]_n} w\s_n(A)$. Obviously $A \subseteq \widetilde{A}$. Since $A \cap w\s_n(A) \neq \emptyset$ and $w\s_n(A)$ is connected in $(T_{n + m}, E_{n + m}^{\ell}(X^c))$, we see that $\widetilde{A}$ is connected. This yields $A = \widetilde{A}$.\par
Next, we show $[A]_n \in \Con_T(n, (\partial_{X^c}{(\s_n(A))})^c)$. Let $\c = (w(1), \ldots, w(k))$ be a path in $(T_{n + m}, E_{n + m}^{\ell}(X^c))$. Let $w_i = [w(i)]_n$ for $i = 1, \ldots, k$. Suppose that $w_i \neq w_{i + 1}$. Then Lemma~\ref{WNG.lemma21} implies that there exists $u \in T_m$ such that $w(i) = w_iu$ and $w(i + 1) = w_{i + 1}u$ for some $u \in \s_n(A)$. Moreover, there exists $t \in \BbZ_J$ such that $Q_{w_i} \cap Q_{w_{i + 1}} = b_t(w_i) = b_t(w_{i + 1})$. On the other hand, since $(w(i), w(i + 1)) \in E_{n + m}^{\ell}(X^c)$. there exists $s \in \BbZ_J$ such that $Q_{w(i)} \cap Q_{w(i + 1)} = b_s(w_iu) = b_s(w_{i + 1}u)$. Hence $b_s(u) \subseteq b_t$ and hence $t \in \partial_{X^c}(\s_n(A))$. Thus we conclude that $(w_i, w_{i + 1}) \in E_n^{\ell}(\partial_{X^c}(\s_n(A)))$. This shows that $[\c]_n$ is a path in $(T_n, E_n^{\ell}(\partial_{X^c}(\s_n(A)))$. Thus $[A]_n$ is connected in $(T_n, E_n^{\ell}(\partial_{X^c}(\s_n(A))))$. Next assume that $(w, v) \in E_n^{\ell}(\partial_{X^c}(\s_n(A)))$ for some $w \in [A]_n$ and $v \in T_n$.Then there exists $t \in \partial_{X^c}(\s_n(A))$ such that $Q_w \cap Q_v = b_t(w) = b_t(v)$. Moreover since $t \in \partial_{X^c}(\s_n(A))$, there exists $u \in \s_n(A)$ and $s \in X^c$ such that $b_s(u) \subseteq b_t$. Using Lemma~\ref{WNG.lemma10}-(2), it follows that $Q_{wu} \cap Q_{vu} = b_s(wu) = b_s(vu)$. Therefore $(wu,vu) \in E_{n + m}^{\ell}(X^c)$. Since $wu \in A$, we see $vu \in A$, so that $v \in [A]_n$. Consequently, we have shown that $[A]_n$ is a connected component of $(T_n, E_n^{\ell}(\partial_{X^c}(\s_n(A))))$. This immediately yields the desired statement.\par
Finally, we show \eqref{WNG.eq10}. Assume that $t \in \partial_{\partial_{X^c}(\s_n(A))}[A]_n$. Then there exists $s_1 \in \partial_{X^c}(\s_n(A))$ and $w \in [A]_n$ such that $b_{s_1}(w) \subseteq b_t$. Moreover, since $s_1 \in \partial_{X^c}(\s_n(A))$, there exists $s \in X^c$ and $u \in \s_n(A)$ such that$b_s(u) \subseteq b_{s_1}$. Combining the preceding two relations, we see that $b _s(wu) \subseteq b_{s_1}(w) \subseteq b_t$. This shows $t \in \partial_{X^c}A$. Conversely let $t \in \partial_{X^c}A$. Then there exists $s \in X^c$, $w \in [A]_n$ and $u \in \s_n(A)$ such that $b_s(wu) \subseteq b_t$. Since one of the boundaries $\{b_{s_1}(w)\}_{s \in \BbZ_J}$ of $Q_w$ must include $b_s(wu)$ and be included in $b_t$, there exists $s_1 \in \BbZ_J$ such that $b_s(wu) \subseteq b_{s_1}(w) \subseteq b_t$. Since $b_s(u) \subseteq b_{s_1}$, it follows that $s_1 \in \partial_{X^c}(\s_n(A))$. This along with the fact that $b_{s_1}(w) \subseteq b_t$ implies that $\partial_{\partial_{X^c}(\s_n(A))}[A]_n$. 
\enddemo

\demo[Proof of Theorem~\ref{WNG.thm10}]
For any $Z \in F_{\partial}(n + m, X)$, there exists $A \in \Con_T(n + m, X)$ such that $Z = (\partial_{X^c}{A})^c$. Lemma~\ref{WNG.lemma40} shows that if $Y = (\partial_{X^c}{\s_n(A)})^c$, then $Y \in F_{\partial}(m, X)$ and $Z \in F_{\partial}(n, Y)$. \par
Conversely, let $Y \in F_{\partial}(m, X)$ and let $Z \in F_{\partial}(n, Y)$.  Then there exists $D \in \Con_T(m, X)$ and $F \in \Con_T(n, Y)$ such that $Y = (\partial_{X^c}D)^c$ and $Z = (\partial_{Y^c}F)^c$. Define $A = \{wu| w \in F, u \in D\}$. \\
{\bf Claim} $A \in \Con_T(n + m, X)$.\\
Proof of Claim: Let $w_1, w_2 \in F$ and $u_1, u_2 \in D$. Then there is a $(T_n, E_n^{\ell}(Y^c))$-path $(w(1), \ldots, w(k))$ in $F$ satisfying $w(1) = w_1$ and $w(k) = w_2$. By the fact that $(w(i), w(i + 1)) \in E_n^{\ell}(Y^c)$, we see that $Q_{w(i)} \cap Q_{w(i + 1)} = b_t(w(i)) = b_t(w(i + 1))$ for some $t \in Y^c = \partial_{X^c}D$, which implies $b_{s_i}(v_i) \subseteq b_t$ for some $s_i \in X^c$ and $v_i \in D$. Since $b_s(w(i)v_i) \subseteq b_t(w(i))$, Lemma~\ref{WNG.lemma10}-(1) shows that
\[
Q_{w(i)v_i} \cap Q_{w(i + 1)v_i} = b_s(w(i)v_i) = b_s(w(i + 1)v_i)
\]
and hence $(w(i)v_i, w(i + 1)v_i) \in E_{n + m}^{\ell}(X^c)$. Set $v_0 = u_1$ and $v_k = u_2$. Let $i = 1, \ldots, k$. Since $v_{i - 1}, v_i \in D$, there exists a $(T_m, E_m^{\ell}(X^c))$-path $(v_i(1), \ldots, v_i(n_i))$ such that $v_i(1) = v_{i - 1}$ and $v_i(n_i) = v_i$. Define $\c_i = (w(i)v_i(1), \ldots, w(i)v_i(n_i))$. Then $\c_1$ is a $(T_{n + m}, E_{n + m}^{\ell}(X^c)$-path. Concatenating $\c_1, \ldots, \c_k$, we obtain a $(T_{n + m}, E_{n + m}^{\ell}(X^c))$-path in $A$ between $w_1u_1$ and $w_2u_2$. Thus we have shown that $A$ is connected in $(T_{n + m}, E_{n + m}^{\ell}(X^c))$. Suppose that $(w, v) \in E_{n + m}^{\ell}(X^c)$ and $w \in A$. If $[w]_n = [v]_n$, then $([\s_n(w), \s_n(v)) \in E_{m}^{\ell}(X^c)$. Since $\s_n(w) \in D$ and $D \in \Con_T(m, X)$, $\s_n(v) \in D$ and hence $v \in A$. In the case $[w]_n \neq [v]_n$, then there exist $t \in \BbZ_J$ and $s \in X^c$ such that
\[
Q_{w} \cap Q_v = b_s(w) = b_s(v) \subseteq Q_{[w]_n} \cap Q_{[v]_n} = b_t([w]_n) = b_t([v]_n).
\]
Using Lemma~\ref{WNG.lemma10}-(2), we see that $\s_n(v) = \s_n(w) \in D$. Moreover, since $b_s(\s_n(w)) \subseteq b_t$, it follows that $t \in \partial_{X^c}D$. Thus $([w]_n, [v]_n) \in (T_n, E_n^{\ell}(\partial_{X^c}D))$ and hence $[v]_n \in F$. This shows $v \in A$. This completes a proof of the claim. \qed\\
Now that $A \in \Con_T(n + m, X)$, Lemma~\ref{WNG.lemma40} yields that $\partial_{X^c}A = \partial_{\partial_{X^c}D}[A]_n = \partial_{Y^c}F$. Therefore we have $Z = \partial_{X^c}A \in F_{\partial}(n + m, X)$.
\enddemo

Next, we start to prove Theorem~\ref{WNG.thm20}.

\lemma\label{WNG.lemma50}
\[
\sup\{\#(A)| n \ge 1, A \in \Con_T(n, X), X \in \B^L\} < +\infty
\]
\endlemma

\demo
Since $\B^L$ is a finite set, it is enough to show that
\[
\sup\{\#(A)| n \ge 1, A \in \Con_T(n, X)\} < \infty
\]
for any $X \in  \B^L$. First if $X = \BbZ_J$, then $\#(A) = 1$ for any $n \ge 1$ and $A \in \Con_T(n, X)$. Next let $X \in \B_{J - 1}$, say $X = \sd{\BbZ_J}{\{s\}}$ and let $A \in \Con_T(n, X)$. Suppose that $A$ contains two distinct elements $w$ and $v$ and that $(w, v) \in E_n^{\ell}$. Then $Q_w \cap Q_v = b_s(w) = b_s(v)$. This implies that there exists no $u \in A$ other than $v$ such that $Q_w \cap Q_u$ is a line segment. Thus $A = \{w, v\}$. Hence $\$(A) \le 2$ for any $n \ge 1$ and $A \in \Con_T(n, X)$. The remaining case is when $X \in \sd{\B_{J - 2}}{\B_{J - 2}^{os}}$. Let $X = \sd{\BbZ_J}{\{s, t\}}$ with $s \neq t$.  Then there exists $j_* \in \{1, \ldots, j_0\}$ such that $s = t \pm j_*$, where  $j_0$ is the integer part of $(J - 1)/2$. Let $A \in \Con_T(n, X)$.  For $w \in A$, define $\ell_s(w)$ and $\ell_t(w)$ as the straight lines containing $b_s(w)$ and $b_t(w)$ respectively and define $p(w)$ as the intersection of $\ell_s(w)$ and $\ell_t(w)$. Note that the angle between $\ell_s(w)$ and $\ell_t(w)$ is
\[
\theta_* = \Big(1 - 2\frac{j_*}J\Big)\pi.
\]
Suppose that $w, v \in A$ and $Q_w \cap Q_v = b_s(w) = b_s(v)$. Then, since $R_{w, v}\circ{f_w} = f_v$, it follows that $R_{w, v}(b_t(w)) = b_t(v)$. Hence $R_{w, v}(\ell_t(w)) = \ell_t(v)$ and $p(w) = p(v)$. Now let $w_1, w_2 \in A$. Then there exists a $(T_n, E_n^{\ell}(\{s, t\})$-path between $w_1$ and $w_2$. Using the above argument inductively, we see that $p(w_1) = p(w_2)$. Thus, we have shown that $p(w)$ does not depend on $w \in A$. Therefore 
\[
\#(A) \le \text{the integer part of $2\pi/\theta_*$.}
\]
\enddemo

\lemma\label{WNG.lemma60}
Let $X \subseteq \BbZ_J$. If there exists $n_0$ such that $F_{\partial}^{n_0}(X) \subseteq \B^L$, then
\[
\sup\{\#(A)| n \ge 1,  A \in \Con_T(n, X)\} < \infty.
\]
\endlemma

\demo
Let $A \in \Con_T(n, X)$. Assume that $n \ge n_0$. By Lemma~\ref{WNG.lemma40}, there exists $Y \in  F_{\partial}^{n_0}(X)$ such that $[A]_{n - n_0} \in \Con_T(n - n_0, Y)$, $A = \bigcup_{w \in [A]_n} w\s_{n - n_0}(A)$ and $\s_{n - n_0}(A) \in \Con_T(n_0, X)$. Therefore, 
\begin{multline*}
\#(A) \le \sup\{\#(A')| n' \ge 1, A' \in \Con_T(n', Y), Y \in \B^L\} \times\\
\sup\{\#(B)| B \in \Con_T(n_0, X)\},
\end{multline*}
where the right-hand side is independent of $n$. Next, assume that $1 \le n < n_0$. Since $\Con_T(n, X)$ is a finite set, we immediately have
\[
\sup\{\#(A)| A \in \Con_T(n, X), 1 \le n <  n_0\} < \infty.
\]
Thus, we have shown the desired statement.
\enddemo

\lemma\label{WNG.lemma70}
Suppose that $M \ge M_J$, $m \ge 1$ and $\c \in \C^{\ell}_{M, m}(w)$. Set $n = |w|$ and let $Y = (\partial(\s_n(\c)))^e$. Then $Y^c \in \B^H$ and there exists $A \in \Con_T(n, Y^c)$ such that $[\c]_n \subseteq A$.
\endlemma

\demo
Let $l_* = l_n(\c)$. Moreover, let $[\c]_n = (w(1), \ldots, w(l_*))$ and let $\s_n(\c)_i = (u_i(1), \ldots, u_i(k_i))$ for $i = 1, \ldots, l_*$. Let $u_i = u_i(k_i)$. Then Lemma~\ref{WNG.lemma21} implies that $u_i = u_{i + 1}(1)$ for any $i = 1, \ldots, l_* -1$. Moreover 
\begin{multline*}
Q_{w(i)u_i} \cap Q_{w(i + 1)u_{i}} = b_{t_i}(w(i)u_i) = b_{t_i}(w(i + 1)u_i) \\
                               \subseteq Q_{w(i)} \cap Q_{w(i + 1)} = b_{s_i}(w(i)) = b_{s_i}(w(i + 1))
\end{multline*}
for some $s_i, t_i \in \BbZ_J$. This shows that $s_i  \in Y$. Hence $[\c]_n$ is connected in $(T_n, E_n^{\ell}(Y))$, so that there exists $A \in \Con_T(n, Y^c)$ such that $[\c]_n \subseteq A$.  Since $G = \{I\}$ and $\c \in \C^{\ell}_{M, m}(w)$, we see that $\H_*(\c) = \s_n(\c)$ and $\s_n(\c)$ is connected in $(T_m, E_m^{\ell})$. Now Theorem~\ref{BAS.thm20} shows that $Y^c \in \B^H$.
\enddemo

\lemma\label{WNG.lemma80}
Let $X \subseteq \BbZ_J$. If $\emptyset \in F_{\partial}^n(X)$, then $F_{\partial}^n(X) \subseteq \{\emptyset\} \cup \B^L$.
\endlemma

\demo
Since $\emptyset \in F_{\partial}^n(X)$, there exists $A \in \Con_T(n, X)$ such that $\partial_{X^c}A = \BbZ_J$. Let $B \in \Con(n, X)$ with $B \neq A$. This means $A \cap B = \emptyset$. Suppose that $s, t \in \partial_{X^c}(B)$ and $\d(s, t) \ge 2$. Then there exists a path $(v_1, \ldots, v_k)$ of $(T_n, E_n^{\ell}(X^c))$ in $B$ such that $b_{s'}(v_1) \subseteq b_s$ and $b_{t'}(v_k) \subseteq b_t$ for some $s', t' \in X^c$. On the other hand, since $s - 1, s + 1 \in \partial_{X^c}A$, there exists a path $(u_1, \ldots, u_l)$ of $(T_n, E_n^{\ell}(X^c))$ in $A$ such that $b_{s''}(u_1) \subseteq b_{s - 1}$ and $b_{t''}(u_l) \subseteq b_{s + 1}$ for some $s'', t'' \in X^c$. Then two paths $(v_1, \ldots, v_k)$ and $(u_1, \ldots, u_l)$ must have intersection but this contradicts the fact that $A \cap B = \emptyset$. Thus $\partial_{X^c}B \subseteq \{s_1, s_1 + 1\}$ for some $s_1 \in \BbZ_J$ and hence $(\partial_{X^c}B)^c \in \B^L$.
\enddemo

\lemma\label{WNG.lemma90}
$F_{\partial}(\B^L) \subseteq \B^L$.
\endlemma

\demo
Let $X = \sd{\BbZ_J}{\{s, t\}} \in \sd{\B_{J - 2}}{\B_{J - 2}^{os}}$ and let $A \in \Con_T(1, X^c)$. Then for any $w \in A$, it follows that
\[
\#(\{v| v \in A, (w, v) \in E_1^{\ell}(X^c)\}) \le 2.
\]
This implies that $A = \{v_1, \ldots, v_k\}$ for some $v_1, \ldots, v_k \in A$ and that $(v_j, v_{j + 1}) \in E_1^{\ell}(X^c)$ for any $j = 1, \ldots, k - 1$. Note that if $w, v, u \in A$ and $(w, u), (w, v) \in E_1^{\ell}(X^c)$, then both $b_s(w)$ and $b_t(w)$ can not be included in $b_p$ for any $p \in \BbZ_J$. Thus, we see $\#(\partial_{X^c}(A)) \le 2$. Assume that $\#(\partial_{X^c}(A)) = 2$. Let $\ell_j$ be the straight line including $K_{v_j} \cap K_{v_{j + 1}}$ for $j = 1, \ldots, k - 1$. Then the lines $\{\ell_j\}_{j = 1, \ldots, k - 1}$ have a common point and the angle between $\ell_j$ and $\ell_{j + 1}$ are the same for any $j = 1, \ldots, k - 2$. This shows that $\partial_{X^c}(A)^c \notin \B_{J - 2}^{os}$. Thus we have shown that $F_{\partial}(X) \in \B^L$ for any $X \in \sd{\B_{J - 2}}{\B_{J - 2}^{os}}$.\par
Next assume that $X = \sd{\BbZ_J}{\{s\}} \in \B_{J - 1}$. In this case, we immediately see that $\#(\partial_{X^c}(A)) \le 1$ for any $A \in E_1^{\ell}(X^c)$. Hence $F_{\partial}(X) \subseteq \B_{J - 1} \cup \B_J$. So, we have the desired consequence for every case.
\enddemo

\demo[Proof of Theorem~\ref{WNG.thm20}]
First we show the equivalence of the three conditions $(F_{\partial}1)$, $(F_{\partial}2)$ and $(F_{\partial}3)$. The equivalence of $(F_{\partial}2)$ and $(F_{\partial}3)$ is immediate by Lemma~\ref{WNG.lemma80}. \par
Assume $(F_{\partial}3)$ and suppose that there exists $X \in \B^H$ and $m \ge 1$ such that $X \in F_{\partial}^m(X)$. Then using inductive argument, we see that $X \in F_{\partial}^{mk}(X)$ for any $k \ge 1$. However, by $(F_{\partial}3)$, there exists $n \ge 1$ such that $F_{\partial}^n(X) \subseteq \{\emptyset\} \cup \B^{L}$. Since $F_{\partial}(\{\emptyset\} \cup \B^L) \subseteq \{\emptyset\} \cup \B^L$ by Lemma~\ref{WNG.lemma90}, it follows that $F_{\partial}^l(X) \subseteq \{\emptyset\} \cup \B^L$ for any $l \ge n$. This contradiction shows that $(F_{\partial}3)$ implies $(F_{\partial}1)$. \par
Next assume that $(F_{\partial}1)$ holds but $(F_{\partial}3)$ does not, i.e. there exists $X \in \B^H$ such that $F_{\partial}^n(X) \cap \B^H \neq \emptyset$ for any $n \ge 1$. Define
\[
\B_* = \{Y| Y \in \B^H, F_{\partial}^n(Y) \cap \B^H \neq \emptyset\,\,\text{for any $n \ge 1$}\}.
\]
{\bf Claim 0.}\,\,$\B_* \cap F_{\partial}(Y) \neq \emptyset$ for any $Y \in \B_*$.\\
Proof.\,\,Suppose that $\B_* \cap F_{\partial}(Y) = \emptyset$. Then for any $Z \in F_{\partial}(Y)$, there exists $n(Z) \ge 1$ such that $F_{\partial}^{n(Z)}(Z) \subseteq \{\emptyset\} \cup \B^L$. Let $n_* = \max_{Z \in F_{\partial}(Y)}n(Z)$. Since $F_{\partial}(\{\emptyset\} \cup \B^L) \subseteq \{\emptyset\} \cup \B^L$, we see that $F_{\partial}^{n_*}(Y) \subseteq \{\emptyset\} \cup \B^L$. This contradicts the fact that $Y \in \B_*$. Thus, we have shown the claim.\qed\par
Now since $X \in \B_*$, using the above claim inductively, we obtain a sequence $X_1, X_2, \ldots \in \B_*$ such that $X_1 = X$ and $X_{j + 1} \in F_{\partial}(X_j)$ for any $j \ge 1$. Since $\B_*$ is a finite set, there exist $k, l \ge 1$ such that $X_k = X_{k + l}$. Let $X_* = X_k$. Then $X_* \in \B_* \subseteq \B^H$ and $X_* \in F_{\partial}l(X_*)$. This contradicts $(F_{\partial}1)$. Thus we have show that $(F_{\partial}1)$ implies $(F_{\partial}3)$ and this concludes the first part of the proof.\par
Second we assume $(F_{\partial}3)$ and show that $(K, d_*)$ is $p$-conductively homogeneous for any $p > \dim_{AR}(K, d_*)$.
Define
\[
\B^1 = \{X | X  \in \B^H, \text{there exists $n_0 \ge 1$ such that}\,\,F_{\partial}^{n_0}(X) \subseteq  \B^L\}
\]
and
\[
\B^2 = \{X | X  \in \B^H, \emptyset \in F_{\partial}^{n_0}(X)\,\,\text{for some $n_0 \ge 0$}\}.
\]
Since $\B^1$ is a finite set, Lemma~\ref{WNG.lemma60} shows 
\[
\sup\{\#(A)| n \ge 1, X \in \B^1, A \in \Con_T(n, X)\} < \infty.
\]
Let $N_1$ be the above supremum. Also since $\B^2$ is finite, there exists $N_2 \ge 1$ such that $\emptyset \in \Con_T(N_2, X)$ for any  $X \in \B_2$.\par
Set $M_1 = \max\{M_J, N_1 + 1\}$. Let $M \ge M_1$, $w \in T$ and $\c \in \C^{\ell}_{M, m}(w)$.  Set $n = |w|$. Define $Y = (\partial(\s_n(\c)))^e$ and $X = Y^c$. \\
{\bf Claim 1.}\,\, $X \in \B^2$.\\
Proof. \,\,Combining  Lemma~\ref{WNG.lemma70} and the assumption of Theorem~\ref{WNG.thm20}, we see that $X \in \B^1 \cup \B^2$. Suppose that $X \in \B^1$. Then by Lemma~\ref{WNG.lemma70}, it follows that $\#([\c]_n) \le N_1 < M_1$. On the other hand, since $\c \in \C^{\ell}_{M, m}(w)$, we see that $\#([\c]_n) \ge M_1$.  This contradiction shows that $X \in \B^2$. \\
{\bf Claim 2.}\,\,There exists $A \subseteq \H_{n, N_2, m}(\c)$ such that $K(A)$ is connected and $\partial{A} = \BbZ_J$.\\
Proof.\,\,Note that $\H_{n, N_2, m}(\c) = \{\pi^{N_2}(v\s_n(u))| v \in T_{N_2}, u \in \c\}$. In the case $N_2 \ge m$, then $\H_{n, N_2, m}(\c) = T_m$ and hence the claim is obvious. Suppose that $N_2 < m$. By Claim 1, there exists $B \in \Con_T(N_2, X)$ such that $\partial_{X^c}B = \BbZ_J$. Define
\[
A_* = \{vu| v \in B, u \in  \s_n(\c)\}.
\]
First we show that $A_*$ is $(T_{N_2 + m}, E_{N_2 + m}^{\ell})$-connected. Let $v_1, v_2 \in B$ and let $\tau_1, \tau_2 \in \s_n(\c)$. Since $B$ is a connected component of $(T_{N_2}, E_{N_2}^{\ell}(Y))$, there exists a path $(v(1), \ldots, v(l))$ of $(T_{N_2}, E_{N_2}^{\ell})$ such that $v(1) = v_1$ and $v_2 = v(l)$. For each $i = 1, \ldots, l - 1$, there exists $s_i \in (\partial(\s_n(\c)))^e$ such that $Q_{v(i)} \cap Q_{v(i)} = b_{s_i}(v(i)) = b_{s_i}(v(i + 1))$. Since $s_i \in (\partial(\s_n(\c)))^e$ and there exists no isolated contact point of cells, Lemma~\ref{PTE.lemma10} implies that there exists $u_i \in \s_n(\c)$ and $t_i \in \BbZ_J$ such that $b_{t_i}(u) \subseteq b_{s_i}$. Then by Lemma~\ref{WNG.lemma10}, it follows that
\begin{multline*}
b_{t_i}(v(i)u_i) = b_{t_i}(v(i + 1)u_i) = Q_{v(i)u_i} \cap Q_{v(i + 1)u_i} \\
\subseteq Q_{v(i)} \cap Q_{v(i + 1)} = b_{s_i}(v(i)) = b_{s_i}(v(i + 1)).
\end{multline*}
Set $u_0 = \tau_1$ and $u_l = \tau_2$. Since $\s_n(\c)$ is $(T_m, E_m^{\ell})$-connected, there exists a $(T_m, E_m^{\ell})$-path $\c_i$ that is contained in $\s_n(\c)$ and connect $u_{i - 1}$ and $u_i$. Concatenating $\{v(i)\c_i\}_{i = 1, \ldots, l}$, we obtain a $(T_{N_2 + m}, E_{N_2 + m}^{\ell})$-path between $v_i\tau_1$ and $v_2\tau_2$. Thus we have shown that $A_*$ is $(T_{N_2 + m}, E_{N_2 + m}^{\ell})$-connected. Furthermore, since $\partial_{Y}B = \BbZ_J$, for any $t \in \BbZ_J$, there exist $s \in (\partial(\s_n(\c)))^e$, $s' \in \BbZ_J$, $v \in B$ and $u \in \s_n(\c)$ such that $b_s(v) \subseteq b_t$ and $b_{s'}(u) \subseteq b_s$. Thus we see that $b_{s'}(uv) \subseteq b_t$. Consequently $t \in \partial{A_*}$, so that $\partial{A_*} = \BbZ_J$. Now letting $A = \pi_{N_2}(A_*)$, we see that $A \subseteq \H_{n, N_2, m}(\c)$, $K(A)$ is connected and $\partial{A} = \BbZ_J$. Thus, Claim 2 has been obtained. \qed\par
Now that Claim 2 is verified, Theorem~\ref{BAS.thm10} yields the desired conclusion.
\enddemo

\section{List of definitions and notations}

\noindent{\bf Definitions}\\
alternated --- Definition~\ref{BAS.def30}\\
boundary of the regular $J$-gon --- Lemma~\ref{RPB.lemma20}\\
conductance constant --- Definition~\ref{CHC.def10}\\
conductively homogeneous --- Definition~\ref{CHC.def30}\\
connect $0$ and $1$ --- Definition~\ref{PTE.def10}\\
connected, connected component (of a graph) --- Definition~\ref{RPB.def00}\\
essential boundary segment -- Definition~\ref{COP.def20}\\
folding map -- Definition~\ref{WNG.def10}\\
folding of a path --- Lemma~\ref{BAS.lemma200}\\
$G$-symmetric $J$-gon-based self-similar system --- Definition~\ref{GPS.def50}\\
invariant --- Definition~\ref{COP.def30}\\
isolated contact point of cells --- Definition~\ref{COP.def10}\\
$(J, G)$-s.s.\,system --- Definition~\ref{GPS.def50}\\
knight move --- After Theorem~\ref{CHC.thm40}\\
length of a path -- Definition~\ref{RPB.def00}\\
modulus of a family of curves --- Definition~\ref{BAS.def40}\\
neighbor disparity constant -- Definition~\ref{CHC.def20}\\
path of a graph -- Definition~\ref{RPB.def00}\\
self-similar set --- Proposition~\ref{GPS.prop10}\\
self-similarity of energy --- Theorem~\ref{CHC.thm30}\\
simple (path) --- Definition~\ref{RPB.def00}\\
transitive -- Definition~\ref{COP.def30}\\
Regular $J$-gon --- Definition~\ref{GPS.def10}\\

\noindent{\bf Notations}\\
$A_j(\c)$ --- Definition~\ref{MIS.lemma20}\\
$\A_m^{\#}(A_1, A_2)$, $A_m(A_1, A_2)$ --- Definition~\ref{BAS.def40}\\
$B_*(x, r)$ --- Theorem~\ref{RPB.thm10}\\
$\B_l$, $\B_{J - 2}^{os}, \B^H, \B^L$ --- Definition~\ref{WNG.def40}\\
$b_i$ --- Definition~\ref{GPS.def10}\\
$b_i(w)$ --- Definition~\ref{GPS.def60}\\
$\C_m^{\#}(A_1, A_2), \C_m(A_1, A_2)$ --- Definition~\ref{BAS.def40}\\
$\C_{M, m}(w)$, $\C^{\ell}_{M, m}(w)$ --- Definition~\ref{BAS.def10}\\
$\C^{\ell}_{m, \partial}, \C^{\ell}_{m, \partial}(X)$ --- Definition~\ref{BCM.def20}\\
$\C^i_{m, \partial}$ --- Lemma~\ref{MIS.lemma40}\\
$\Con_T(n, X), \Con_T(X)$ --- Definition~\ref{MIS.def200}\\
$D_J$ --- Definition~\ref{GPS.def30}\\
$D_J^*$ --- Lemma~\ref{GPS.lemma00}\\
$D_{J/2}^V$ --- Definition~\ref{MTH.def10}\\
$d_*(\cdot, \cdot)$ --- Proposition~\ref{GPS.prop10}, Lemma~\ref{RPB.lemma10}\\
$d_{Q_*}$ --- Definition~\ref{PTE.def40}\\
$\dim_H(K, d_*)$ --- Proposition~\ref{GPS.prop10}\\
$\diam{A, d}$ --- Definition~\ref{GPS.def60}\\
$E_n^*$ --- \eqref{GPS.eq200}, Definition~\ref{RPB.def10}\\
$E_n^A$ for $A \subseteq T_n$ --- Definition~\ref{RPB.def10}\\
$E_1^{\ell}$ --- Definition~\ref{GPS.def50}\\
$E_n^{\ell}$ --- Lemma~\ref{BAS.thm00}\\
$E_n^{\ell}(Y)$ --- Definition~\ref{MIS.def200}\\
$\E_{p, A}^n(\cdot)$ --- Definition~\ref{CHC.def10}\\
$\E_{p, m}(A_1, A_2, A)$ --- Definition~\ref{CHC.def10}\\
$\E_{M, p, m}(w)$ --- Definition~\ref{CHC.def10}\\
$\hE_p(\cdot)$ --- Theorem~\ref{CHC.thm30}\\
$F_{\partial}(X)$ --- Definition~\ref{MIS.def200}\\
$F_{\partial}(n, X)$ --- Definition~\ref{MIS.def200}\\
$G(A)$ for $A \subseteq T_n$ --- Definition~\ref{BAS.def30}\\
$G(X)$ for $X \subseteq \BbZ_J$ --- Definition~\ref{COP.def20}\\
$G_*(S, \{f\}_{s \in S})$, $G_*$ --- Definition~\ref{GPS.def70}\\
$g_*$ --- Proposition~\ref{GPS.prop20}, Lemma~\ref{BAS.lemma10}, Definition~\ref{COP.def20}\\
$g_{w, v}$ --- Definition~\ref{PTE.def21}\\
$\H_{n_1, n_2, m}(\cdot)$ --- Definition~\ref{BAS.def10}\\
$\H_*(\c)$ --- Definition~\ref{BAS.def20}\\
$K(A)$ -- Definition~\ref{BAS.def10}\\
$K_w$ --- Definition~\ref{GPS.def60}\\
$L_*$ --- Proposition~\ref{RPB.prop05}\\
$l(\cdot)$ --- Definition~\ref{RPB.def00}\\
$\ell(\cdot)$ --- Definition~\ref{CHC.def10}\\
$l_k(\c)$ --- Definition~\ref{PTE.def20}\\
$M_J$ --- Theorem~\ref{BAS.thm20}\\
$\M_{p, m}^{\#}(A_1, A_2), \M_{p, m}(A_1, A_2), \M_{M, p, m}^{\#}(w), \M_{M, p, m}(w)$ --- Definition~\ref{BAS.def40}\\
$N_*$ --- Theorem~\ref{RPB.thm10}\\
$\N_p(\cdot)$ --- Definition~\ref{CHC.def40}\\
$O_w$ --- Lemma~\ref{GPS.lemma20}\\
$P_m(\cdot)$ --- Definition~\ref{CHC.def40}\\
$\P_{\c}$ --- Theorem~\ref{BAS.thm10}\\
$p_i$ --- Definition~\ref{GPS.def10}\\
$Q(A)$, $Q_0(A)$ --- Definition~\ref{PTE.def40}\\
$Q_*^{(J)}, Q_*$ --- Definition~\ref{GPS.def10}\\
$Q^{(m)}$ --- Definition~\ref{PTE.def40}\\
$Q_n(x)$ --- Definition~\ref{COP.def10}\\
$Q_w$ --- Definition~\ref{GPS.def60}\\
$Q_w^{(m)}$ --- Lemma~\ref{GPS.lemma30}\\
$Rot_J$ --- Definition~\ref{GPS.def30}\\
$R_{w, v}$ --- Lemma~\ref{BAS.thm00}\\
$R_{w, v}^*$ --- Lemma~\ref{BAS.lemma200}\\
$R_{\theta}$ --- Definition~\ref{GPS.def20}\\
$(S, \{f_s\}_{s \in S}, G)$ --- Definition~\ref{GPS.def50}\\
$S^m(\cdot)$ --- Definition~\ref{GPS.def60}\\
$T$ --- Definition~\ref{GPS.def60}\\
$T^{(k)}$ --- Corollary~\ref{RPB.cor10}\\
$T_m$ --- Definition~\ref{GPS.def60}\\
$T(w)$ --- Lemma~\ref{BAS.lemma201}\\
$U(x : n)$ --- Definition~\ref{RPB.def10}\\
$\W^p$ --- Definition~\ref{CHC.def40}\\
$\BbZ_J$ --- Definition~\ref{GPS.def10}\\
$(\BbZ_J)^e$ --- Definition~\ref{COP.def20}\\
$\BbZ_J^i$ --- \eqref{ENJ.eq10}\\
$\d_J(\cdot, \cdot), \d(\cdot, \cdot)$ --- Definition~\ref{GPS.def10}\\
$\GG_M^A(w)$,$\GG_M(w)$ --- Definition~\ref{RPB.def10}\\
$\GG^n(x)$ --- Definition~\ref{RPB.def10}\\
$\Theta_{\theta}$ --- Definition~\ref{GPS.def20}\\
$\mu_*$ --- Theorem~\ref{RPB.thm10}\\
$\pi$ --- Definition~\ref{RPB.def10}\\
$\rho(i)$ --- Definition~\ref{GPS.def40}\\
$\Sigma$ --- Proposition~\ref{GPS.prop15}\\
$\s$ -- Theorem~\ref{CHC.thm20}\\
$\s_n$ -- Definition~\ref{BAS.def10}\\
$\sigma_w^*$ --- Proposition~\ref{GPS.prop15}\\
$\sigma_{p, m, n}$ --- Definition~\ref{CHC.def20}\\
$\s_k(\c)_i$ --- Definition~\ref{PTE.def20}\\
$\s_k^F(\c)$ --- Lemma~\ref{BAS.lemma200}\\
$\chi$ --- Proposition~\ref{GPS.prop15}\\
$\partial{A}$, where $A \subseteq T$ --- Definition~\ref{BAS.def10}\\
$\partial_YA$ --- Definition~\ref{MIS.def200}\\
$\partial{Q_*}$ --- Lemma~\ref{RPB.lemma20}\\
$(\cdot)^e$ --- Definition~\ref{COP.def20}\\
$[w]_n$, where $w \in T$ --- Definition~\ref{BAS.def53}\\
$[\c]_n$, where $\c$ is a path --- Definition~\ref{PTE.def20}\\

\providecommand{\bysame}{\leavevmode\hbox to3em{\hrulefill}\thinspace}
\providecommand{\MR}{\relax\ifhmode\unskip\space\fi MR }
\providecommand{\MRhref}[2]{%
  \href{http://www.ams.org/mathscinet-getitem?mr=#1}{#2}
}
\providecommand{\href}[2]{#2}


\begin{thebibliography}{10}

\bibitem{BB1}
M.~T. Barlow and R.~F. Bass, \emph{The construction of {Brownian} motion on the
  {Sierpinski} carpet}, Ann. Inst. Henri Poincar{\'e} \textbf{25} (1989),
  225--257.

\bibitem{BeHeStr}
T.~Berry, S.~Hellman, and R.~Strichartz, \emph{Outer approximation of the
  spectrum of a fractal {Laplacian}}, Experimental Math. \textbf{18} (2009),
  449--480.

\bibitem{CaHaHu}
C.~Canner, C.~Hayes, R.~Huang, M.~Orwin, and L.~G. Rogers, \emph{{Resistance
  scaling on 4N-carpets}}, Forum Mathematicum \textbf{34} (2022), 61--75.

\bibitem{CaoQiu}
S.~Cao and H.~Qiu, \emph{Dirichlet forms on unconstrained {Sierpinski}
  carpets}, Probab. Theory Relat. Fields \textbf{189} (2024), 613--657.

\bibitem{Cheeger}
J.~Cheeger, \emph{Differentiability of {Lipschitz} functions on metric measure
  spaces}, Geom. Funct. Anal. \textbf{3} (1999), 428--517.

\bibitem{Haj1}
P.~Haj{\l}asz, \emph{Sobolev spaces on an arbitrary metric spaces}, Pontential
  Analysis \textbf{5} (1996), 403--415.

\bibitem{HeiKoShTy}
J.~Heinonen, P.~Koskela, N.~Shanmugalingam, and J.~Tyson, \emph{{Sobolev Spaces
  on Metric Measure Spaces - An Approach Based on Upper Gradients}}, New
  Mathematical Monographs, Cambridge University Press, 2015.

\bibitem{AndrewsIV}
U.~A.~Andrews IV, \emph{{Existence of diffusions on 4N-carpets}}, Doctoral
  Dissertation, University of Connecticut, 2017.

\bibitem{KajMur3}
N.~Kajino and M.~Murugan, \emph{On the conformal walk dimension {II}:
  Non-attainment for some {Sierpi\'{n}ski} carpets}, in preparation.

\bibitem{KajMur2}
\bysame, \emph{On the conformal walk dimension: Quasisymmetric uniformization
  for symmetric diffusions}, Invent. Math. \textbf{231} (2023), 263--405.

\bibitem{AOF}
J.~Kigami, \emph{{Analysis on Fractals}}, Cambridge Tracts in Math. vol. 143,
  Cambridge University Press, 2001.

\bibitem{GAMS}
\bysame, \emph{{Geometry and Analysis of Metric Spaces via Weighted
  Partitions}}, Lecture Notes in Math. vol. 2265, Springer, 2020.

\bibitem{Ki22}
\bysame, \emph{Conductive homogeneity of compact metric spaces and construction
  of $p$-energy}, Mem. European Math. Soc. \textbf{5} (2023).

\bibitem{Ota1}
Y.~Ota, \emph{Conductive homogeneity of self-similar sets; a new class of
  examples}, Doctoral Dissertation, Kyoto University, 2025.

\bibitem{Sasaya1}
K.~Sasaya, private communication.

\bibitem{Shanm1}
N.~Shanmugalingam, \emph{Newtonian spaces: an extension of {Sobolev} spaces to
  metric measure spaces}, Rev. Mat. Iberoamer. \textbf{16} (2000), 243--279.

\end{thebibliography}

\end{document}